\documentclass[12pt,a4paper]{article}
\usepackage[margin=2.5 cm]{geometry}

\usepackage[T1]{fontenc}    
\usepackage{hyperref}       
\usepackage{url}            
\usepackage{nicefrac}       
\usepackage[numbers]{natbib}
\usepackage{xcolor}         
\usepackage{bbm}            
\usepackage{amsmath, amssymb, amsfonts, amssymb, amsthm}
\usepackage{mathtools}
\usepackage{soul} 
\usepackage{comment}
\usepackage[shortlabels]{enumitem}
\usepackage[leqno,fleqn,intlimits]{empheq}

\newcommand{\N}{\mathbb{N}} 
\newcommand{\Z}{\mathbb{Z}}

\newcommand{\R}{\mathbb{R}}

\newcommand\restr[2]{\left.\kern-\nulldelimiterspace #1 \vphantom{|} \right|_{#2}}
\newcommand{\E}{{\mathbb E}}
\newcommand{\cB}{{\mathcal B}}
\newcommand{\cG}{{\mathcal G}}

\renewcommand{\epsilon}{\varepsilon}
\renewcommand{\phi}{\varphi}
\newcommand*\diff{\mathop{}\!\mathrm{d}}
\newcommand{\SMALL}{\textstyle}
\newcommand{\wt}{\widetilde}

\DeclareMathOperator*{\Id}{Id}
\DeclareMathOperator*{\Var}{Var}
\DeclareMathOperator*{\diam}{diam}

\DeclareMathOperator*{\Lip}{Lip}
\DeclareMathOperator*{\diag}{diag}

\numberwithin{equation}{section}

\theoremstyle{definition}

\newtheorem{theorem}{Theorem}[section]
\newtheorem{lemma}[theorem]{Lemma}
\newtheorem{corollary}[theorem]{Corollary}
\newtheorem{proposition}[theorem]{Proposition}
\newtheorem{remark}[theorem]{Remark}
\newtheorem{definition}[theorem]{Definition}

\newtheorem*{definition*}{Definition}
\newtheorem*{theorem*}{Theorem}
\newtheorem*{corollary*}{Corollary}
\newtheorem*{proposition*}{Proposition}
\newtheorem*{lem*}{Lemma}
\newtheorem*{remark*}{Remark}
\newtheorem*{notation*}{Notation}

\title{Rates for maps and flows in a deterministic multidimensional weak invariance principle}

\author{Nicolò Paviato\footnote{Einstein Institute of Mathematics, The Hebrew University of Jerusalem, Jerusalem 91904, Israel. e-mail: {\bf nicolo.paviato@mail.huji.ac.il} - 
Research supported by ISF research grant 3056/21.}}

\begin{document}
\maketitle

\begin{abstract}
We present the first rates of convergence to an $N$-dimensional Brownian motion when $N\ge2$ for discrete and continuous time dynamical systems. Additionally, we provide the first rates for continuous time in any dimension.
Our results hold for nonuniformly hyperbolic and expanding systems, such as Axiom A flows,
suspensions over a 
Young tower with exponential tails, and some classes of intermittent solenoids.    
\end{abstract}

\section{Introduction}

The study of statistical laws for dynamical systems is a widely explored area that attracted the attention of a rich research community, starting from the work of Bowen, Ruelle and Sinai in the 1970s. 
This paper provides quantitative estimates on the number of iterations of a hyperbolic/expanding\footnote{For the remainder of the introduction, we write "hyperbolic" instead of "hyperbolic/expanding".} system to accurately approximate Brownian motion. 
In particular, we focus on the rates of convergence for a deterministic version of the classical functional central limit theorem of Donsker~\cite{Don51}. 
In the following, we will refer to this as the weak invariance principle (WIP). 
Such a result was proved for numerous nonuniformly hyperbolic maps in~\cite{HofKel82} and for uniformly hyperbolic flows in~\cite{DenPhil84}. 
More recent developments in this direction for a nonuniformly hyperbolic setting are~\cite{BalMel18,Gou04,MelNic05,Melvar16,MelZwe15}. 

A natural consequence of the WIP is the central limit theorem, where the rates of convergence are commonly referred to as Berry-Esseen estimates.
Sharp results in this direction for uniformly and nonuniformly hyperbolic diffeomorphisms were found respectively in~\cite{CoePar90} and~\cite{Gou05}. 
In the continuous time literature,~\cite{Pen07} was the first to provide Berry-Esseen estimates; these are of the order of $\mathcal{O}(n^{-1/4+\epsilon})$, $\epsilon>0$, in the Prokhorov metric for 
a billiard flow with finite horizon.

For the WIP, Antoniou and Melbourne~\cite{AntMel19} proved a convergence rate of $\mathcal{O}(n^{-1/4+\delta})$ in the Prokhorov metric for nonuniformly hyperbolic maps, while Liu and Wang~\cite{LiuWan23} proved the same rate in the $q$-Wasserstein  metric for $q>1$. 
Here, $\delta>0$ depends on the degree of nonuniformity and it gets smaller for greater~$q$; it can be chosen arbitrarily small if the system has exponential tails.
The methods of~\cite{AntMel19,LiuWan23} are based on a generalisation~\cite{KorKosMel18} of the martingale-coboundary decomposition technique of Gordin~\cite{Gor69}, which allows to apply a martingale version of the Skorokhod embedding theorem. 
It is known~\cite{Bor73,SawSta72} that such a method cannot yield better rates than $\mathcal{O}(n^{-1/4})$. 
It is important to mention that~\cite{AntMel19} and~\cite{LiuWan23} deal exclusively with discrete time systems and real-valued observables and, to our knowledge, the literature does not have any further results on this topic. 
This paper is the first one that gives rates of convergence in $N$ dimension for $N\ge2$; furthermore, it is the first work that addresses flows in this context.  

When $N=1$, we get a rate of $\mathcal{O}(n^{-1/4}(\log n)^{3/4})$ in Prokhorov for uniformly hyperbolic maps and flows, improving the one of~\cite{AntMel19} in discrete time. 
For nonuniformly hyperbolic flows, we recover the same rate that~\cite{AntMel19} proved for maps. For $N\ge2$, we are able to achieve a rate of $\mathcal{O}(n^{-1/6+\delta})$ in the $1$-Wasserstein  metric independently of the dimension. 

To illustrate how the nonuniformity affects the rates, we consider the LSV map studied in~\cite{LivSauVai99}. 
For~$\gamma>0$ we define $T\colon[0,1]\to[0,1]$ as
\[
Tx=\begin{cases}
x(1+2^\gamma x^\gamma)&x\in[0,1/2)\\
2x-1&x\in[1/2,1].
\end{cases}
\]
This system is a type of Pomeau-Maneville intermittent map~\cite{PomMan80} and when $\gamma\in(0,1)$ there is a unique ergodic invariant probability measure that is absolutely continuous to Lebesgue measure. 
The WIP holds for Hölder observables $v\colon[0,1]\to\R^d$ if~$\gamma\in(0,\frac12)$. 
A first example of a nonuniformly expanding flow is obtained as a suspension of the map $T$ with a Hölder continuous roof function. 
The WIP is also valid on the suspension by~\cite{KelMel16,MelTor04,MelZwe15}. 
When $N=1$, we obtain for both map and flow the rates in Prokhorov of $\mathcal{O}(n^{-(1-2\gamma)/4+\epsilon})$ for $\gamma\in(0,\frac12)$ and $\epsilon>0$, similarly to~\cite{AntMel19}. 
If $N\ge2$, we get in $1$-Wasserstein  $\mathcal{O}(n^{-1/6+\epsilon})$ for $\gamma\in(0,\frac13]$, and $\mathcal{O}(n^{-(1-2\gamma)/2+\epsilon})$ if $\gamma\in(\frac13,\frac12)$.

Our proofs utilize results from general martingale theory \cite{Cou99, CunDedMer20, Kub94}. 
To apply the latter to discrete time dynamical systems, we  follow the same strategy of~\cite{AntMel19,LiuWan23} and rely on an advanced adaptation~\cite{KorKosMel18} of the martingale-coboundary decomposition introduced by Gordin~\cite{Gor69}. 
We proceed similarly in the flow case, where first we generalise~\cite{KorKosMel18} to continuous time; this original work is found in Section~\ref{sec:new_dec}.

The remainder of the paper is organized as follows. 
Section~\ref{sec:setup_and_main_results} presents our main results for flows and maps. 
In Section~\ref{sec:results_for_maps}, we recall some techniques from~\cite{KorKosMel18} and prove the rates for maps. 
Section~\ref{sec:new_dec} presents two new decompositions for regular observables with estimates in continuous time, extending the work of~\cite{KorKosMel18}. 
In Section~\ref{sec:continuous_time_rates}, we use the new estimates from Section~\ref{sec:new_dec} to prove the rates for semiflows. 
Finally, Section~\ref{sec:sinai_flow} shows that our rates for (non invertible) expanding systems are still valid for a family of (invertible) nonuniformly hyperbolic flows that satisfy an exponential contraction along stable leaves.
Apart from Section~\ref{sec:sinai_flow}, this work can be found in the author's PhD thesis~\cite{Pav23}, completed at the University of Warwick~(UK) under the supervision of Prof.\@ Ian Melbourne.

\begin{notation*} 
We write interchangeably  $a_n=\mathcal{O}(b_n)$ or  $a_n\ll b_n$ for two sequences $a_n,b_n\ge0$, if there exists a constant $C>0$ and an integer $n_0\ge0$ such that $a_n\le Cb_n$ for all $n\ge n_0$. 

\noindent For $x\in\R^m$ and $J\in\R^{m\times n}$, we write $|x|=(\sum_{i=1}^mx_i^2)^{1/2}$ and $|J|=(\sum_{i=1}^m\sum_{j=1}^nJ^2_{i,j})^{1/2}$.
\end{notation*}

\section{Setup and main results}\label{sec:setup_and_main_results}
We introduce here the metrics used to describe the rates in the WIP, and the class of dynamical systems under consideration. 
This section presents our new rates of convergence in discrete and continuous time.

\subsection{Metrics for probability measures}\label{subsec:metrics}

Let $(S,d_S)$ be a separable metric space with Borel $\sigma$-algebra $\cB$, and denote with $\mathcal{M}_1(S)$ the set of Borel probability measures on $S$. 
Let $\mu,\nu\in\mathcal{M}_1(S)$; following~\cite{GibSu02}, we have the following metrics on $\mathcal{M}_1(S)$,
\vspace{1ex}

\begin{itemize}
    \item \textbf{1-Wasserstein} (or Kantorovich)
    \begin{equation*}\SMALL
\mathcal{W}(\mu,\nu)=
\sup_{f\in \Lip_1}\bigl|\int_S f\diff\mu-\int_S f\diff\nu\bigr|,
\end{equation*}
where $\Lip_1=\{f\colon S\to \R:|f(x)-f(y)|\le d_S(x,y)\}$.\vspace{1ex}

\item \textbf{ Prokhorov} (or Lévy-Prokhorov)
\begin{equation*}\SMALL
\Pi(\mu,\nu)=
\inf\{\epsilon>0:\mu(B)\le\nu(B^\epsilon)+\epsilon\text{ for all }B\in\cB\},
\end{equation*}
where
$B^\epsilon=\bigcup_{x\in B}\{y\in S\colon d_S(x,y)<\epsilon\}.$
\end{itemize}\vspace{1ex}

If $X$ and $Y$ are $S$-valued random elements with respectively laws $\mu$ and $\nu$, we write $\Pi(X,Y)=\Pi(\mu,\nu)$ and $\mathcal{W}(X,Y)=\mathcal{W}(\mu,\nu)$. 
Denoting with $X_n\to_wX$ the weak convergence of the sequence of laws of $(X_n)_{n\ge1}$ to the law of $X$, we have that $\mathcal{W}(X_n,X)\to0$ implies $X_n\to_w X$.
The distance $\Pi$ metrizes weak convergence on $\mathcal{M}_1(S)$ and the same is true for $\mathcal{W}$ under the extra assumption $\diam(S)<\infty$.

We recall that  $\Pi(X,Y)\le\sqrt{\mathcal{W}(X,Y)}$ (see \cite[Theorem~2]{GibSu02}) and that, if $X$ and~$Y$ are defined on a common probability space, then
$\Pi(X,Y)\le|d_S(X,Y)|_\infty$. 
This estimate follows from the definition of $\Pi$,
noting that 
$\mathbb{P}(d_S(X,Y)>\epsilon)\le\epsilon$ for some $\epsilon>0$ implies that
$\Pi(X,Y)\le\epsilon$.

\subsection{Nonuniformly expanding maps}\label{subsec:nonunif_exp_maps}
Let $(X,d)$ be a bounded metric space with a Borel probability measure $\rho$ and suppose that 
$T\colon X\to X$ 
is a nonsingular map 
($\rho(T^{-1}E)=0$ if and only if $\rho(E)=0$ for all Borel sets
$E\subset X$). Assume that $\rho$ is ergodic.

We suppose that there exists a measurable $Y\subset X$ with $\rho(Y)>0$ and $\{Y_j\}_{j\ge1}$ an at most countable measurable partition of $Y$. Let $\tau\colon Y\to\Z^+$ be an integrable function with constant values $\tau_j\ge1$ on partition elements $Y_j$. We assume that $T^{\tau(y)}y\in Y$ for all $y\in Y$ and define $F\colon Y\to Y$ as $Fy=T^{\tau(y)}y$.

The dynamical system $(X,T,\rho)$ is said to be a \textit{nonuniformly expanding map} if there are constants $ \lambda>1$, $\eta\in(0,1]$, $C\ge1$, such that for each $j\ge1$ and $y,y'\in Y_j$,

\begin{enumerate}[(a)]
\item\label{item:full_branch} $\restr{F}{Y_j}
\colon Y_j\to Y$ is a measure-theoretic bijection;
\item\label{item:expansion} $d(Fy,Fy')\ge \lambda d(y,y')$;
\item\label{item:non_expansion} $d(T^\ell y,T^\ell y')\le Cd(Fy,Fy')$ for all $0\le\ell\le\tau_j-1$;
\item\label{item:zeta} $\zeta=\restr{d\rho}{Y}/(\restr{d\rho}{Y}\circ F)$ satisfies
$ |\log\zeta(y)-\log\zeta(y')|\le Cd(Fy,Fy')^\eta$.
\end{enumerate}
We say that $T$ is nonuniformly expanding of \textit{order} $p\in[1,\infty]$ if the \textit{return time} $\tau$ lies in $ L^p(Y)$. 
A map $F$ satisfying \ref{item:full_branch}, \ref{item:expansion}, and \ref{item:zeta} is a (full-branched) Gibbs-Markov map as in~\cite{AarDen01}.
It is standard that there exists a unique $\rho$-absolutely continuous ergodic (and mixing) $T$-invariant probability measure $\mu_X$ on $X$ (see for example~\cite[Theorem 1]{You99}).

The LSV map mentioned in the introduction is an example of nonuniformly expanding map of order $p$ for every $p\in[1,1/\gamma)$ (see~\cite[Subsection 2.5.2]{Alv20}).

\begin{definition}[Hölder functions]\label{def:Hölder}
Let $N\ge1$, $\eta\in(0,1]$, and $v: X\to \R^N$. Define
\[\SMALL
\|v\|_\eta=|v|_\infty+|v|_\eta\qquad\text{and}
\qquad
|v|_\eta=\sup_{x\neq x'}|v(x)-v(x')|/d(x,x')^\eta.
\]
Let $ C^\eta( X, \R^N)$ consist of observables $v: X\to \R^N$ with
$\|v\|_\eta<\infty$.
\end{definition}

\subsection{Rates in the WIP for maps}\label{subsec:rates_maps}
Let $T\colon X\to X$ be nonuniformly expanding of order $p\in[1,\infty]$ with ergodic invariant measure $\mu_X$. 
For $N\ge1$, let $v\in C^\eta( X, \R^N)$ with mean zero and define $B_n\colon[0,1]\to \R^N$, $n\ge1$, as
\begin{equation*}\label{eq:Bn}
B_{n}(k/n)=\frac{1}{\sqrt{n}}\sum_{j=0}^{k-1}v\circ T^j,
\end{equation*}
for $0\le k\le n$, 
and using linear interpolation in $[0,1]$. The process $B_n$ is a random element of $ C([0,1], \R^N)$ defined on the probability space $( X,\mu_X)$.
Note that the randomness of~$B_n$ comes exclusively from the initial point $x_0\in  X$, chosen according to $\mu_X$.

Here follows a standard result (see for example~\cite{Gou04,KorKosMel18,MelNic05}).

\begin{theorem}\label{theorem:CLT_WIP_discrete} Let $p\in[2,\infty]$ and
$v\in C^\eta( X, \R^N)$ with mean zero. Then,
\begin{enumerate}[(i)] 
\item\label{item:Sigma_discrete} The matrix $\Sigma=\lim_{n\to\infty}n^{-1}\int_{ X}(\sum_{j=0}^{n-1}v\circ T^j)
(\sum_{j=0}^{n-1}v\circ T^j  )^T\diff\mu_ X$ exists and is positive semidefinite. Typically $\Sigma$ is positive definite: there exists a closed subspace $ C_{deg}$ of $ C^\eta( X, \R^N)$ with infinite codimension, such that $\det(\Sigma)\ne 0$ if $v\notin C_{deg}$.
\item\label{item:WIP_discrete} The WIP holds: $B_{n}\to_wW$ in $ C([0,1], \R^N)$ on the probability space $( X,\mu_X)$, where~$W$ is a centred $N$-dimensional Brownian motion with covariance $\Sigma$.\qed
\end{enumerate}
\end{theorem}

The following theorems display rates in the WIP, where the order~$p\in(2,\infty]$ influences the speed of convergence. These rates are stated in the 1-Wasserstein and Prokhorov metrics on $\mathcal{M}_1(S)$, where $S= C([0,1], \R^N)$ with the uniform distance.

\begin{theorem}\label{theorem:rate_wass_maps}
Let $p\in(2,3)$ and $v\in C^\eta( X, \R^N)$ with mean zero. There is~$C>0$ such that
$\mathcal{W}(B_n,W)\le Cn^{-\frac{p-2}{2p}}(\log n)^\frac{p-1}{2p}$ 
for all integers $n>1$.
\end{theorem}

To our knowledge, the rates of Theorem~\ref{theorem:rate_wass_maps} are the first for multidimensional observables in the dynamical systems literature. They are likely not optimal, as one expects an improvement  when $p$ increases (as it happens for $N=1$ in Theorem~\ref{theorem:rates_N=1_maps}). Yet, our proof of Theorem~\ref{theorem:rate_wass_maps} uses modern techniques by~\cite{CunDedMer20}, which do not improve for $p>3$. In such cases, our rates become $\mathcal{O}(n^{-1/6+\epsilon})$ for any $\epsilon>0$.

\begin{remark}
For $N=1$ and $p\ge4$,~\cite[Theorem 3.4]{LiuWan23} of Liu and Wang implies the rate
$\mathcal{W}(B_n,W)\ll n^{-(p-2)/(4(p-1))}$. 
Hence, our Theorem~\ref{theorem:rate_wass_maps} provides a new estimate in $\mathcal{W}$ in one dimension when $p\in(2,3)$. By the inclusion of $L^p$ spaces
\end{remark}

\begin{theorem}\label{theorem:rates_N=1_maps}
Let $p\in(2,\infty]$ and $v\in C^\eta( X,\R)$ with mean zero. 
Then there is~$C>0$ such that
\begin{empheq}[
left=
{\Pi(B_n,W)} \le C
\empheqlbrace
]{align}
&n^{-\frac{p-2}{4p}}
&& p\in(2,\infty), \label{eq:rate_Bn_p}\\
&n^{-1/4}(\log n)^{3/4} 
&& p=\infty \label{eq:rate_Bn_infty}
\end{empheq}
for all integers $n>1$.
\end{theorem}

While the rates in \eqref{eq:rate_Bn_p} were established by Antoniou and Melbourne \cite[Theorem~3.2]{AntMel19}, the estimate in \eqref{eq:rate_Bn_infty} are derived in Section~\ref{sec:results_for_maps} and is a novel contribution of this paper.

\begin{remark}\label{rem:MultidimMap}
Using that $\Pi(X,Y)\le\sqrt{\mathcal{W}(X,Y)}$, Theorem~\ref{theorem:rate_wass_maps} yields for $p\in(2,3)$ and every $N\ge1$ that
$\Pi(B_n,W)\ll n^{-\frac{p-2}{4p}}(\log n)^\frac{p-1}{4p}$. This result is only relevant for $N\ge2$, as Theorem~\ref{theorem:rates_N=1_maps} gives better rates when $N=1$.
\end{remark}

\subsection{Nonuniformly expanding semiflows}\label{subsec:nonunif_exp_semiflow}

Let $( M,d)$ be a bounded metric space and let $\{\Psi_t\colon M\to M\}_{t\ge0}$ be a semiflow, so that $\Psi_0=\Id$ and $\Psi_{s+t}=\Psi_s\circ \Psi_t$ for all $s,t\ge0$. 
Assume Lipschitz continuous dependence  on initial conditions, that is there exists $K>0$ such that, for all $t\in[0,1]$ and $x,y\in M$,
\begin{equation}\label{eq:continuous_dependence}\SMALL
d(\Psi_tx,\Psi_ty)\le K d(x,y).
\end{equation}
We suppose Lipschitz continuity in time: there exists $L>0$ such that, for all $t,s\ge0$ and $x\in M$
\begin{equation}\label{eq:Lipschitz}\SMALL
d(\Psi_tx,\Psi_sx)\le L|t-s|.
\end{equation}

Let $\eta\in(0,1]$. Suppose that there exist a Borel subset $X\subset M$ and a first return time function $r\colon X\to[1,\infty)$, $r(x)=\inf\{t>0:\Psi_tx\in X\}$ and let $T\colon X\to X$ be $Tx=\Psi_{r(x)}x$. Assume that $r\in C^\eta(X,\R)$.
We say that the continuous time system $(M,\Psi_t)$ is a \textit{ nonuniformly expanding semiflow} if $T$ is a nonuniformly expanding map as described in Subsection~\ref{subsec:nonunif_exp_maps}. Moreover, we say that $\Psi_t$ is a nonuniformly expanding flow of \textit{order} $p\in[1,\infty]$ if $(X,T)$ is nonuniformly expanding map of order~$p$.

We define the suspension 
$X^r=\{(x,u)\in X\times[0,\infty):u\in [0,r(x)]\}/\sim$,
where $(x,r(x))\sim (Tx,0)$.
The suspension semiflow $T_t:X^r\to X^r$ is given by
$T_t(x,u)=(x,u+t)$ computed modulo identifications. 
We have the semiconjugacy
$\pi_r\colon X^r\to M$ defined as $\pi_r(x,u)=\Psi_ux$ between $T_t$ and $\Psi_t$, and can now define the ergodic $T_t$-invariant probability measure
$\mu^r=(\mu_X\times{\rm Lebesgue})/\overline{r}$, where $\mu_X$ is the ergodic invariant measure for $T$ and $\overline{r}=\int_Xr\,d\mu_X$. Hence, $\mu_M=(\pi_r)_*\mu^r$ is an ergodic $\Psi_t$-invariant probability measure on $ M$.

We define the space $C^\eta(M,\R^N)$  with norm $\|\cdot\|_\eta$ similarly to Definition~\ref{def:Hölder}.

\subsection{Rates in the WIP for semiflows}
Let $\Psi_t\colon M\to M$ be a nonuniformly expanding semiflow of order $p\in[1,\infty]$ with ergodic invariant measure $\mu_M$. 
Let $T\colon X\to X$ be the corresponding nonuniformly expanding map with ergodic invariant measure $\mu_X$, defined via the first return function $r\in C^\eta(X,\R)$. 

Let $v\in C^\eta( M, \R^N)$ with mean zero and define the sequence $W_n$ as
\begin{equation}\label{eq:Wn}
W_{n}(t)=\frac{1}{\sqrt{n}}\int_0^{nt}v\circ \Psi_s\diff s,
\end{equation}
for $n\ge1$ and $t\in[0,1]$. For fixed $n$, the process $W_n$  is a random element of $ C([0,1], \R^N)$ defined on the probability space $( M,\mu_ M)$.

Let $v_X\colon X\to\R^N$ be $v_X(x)=\int_0^{r(x)}v(\Psi_sx)\diff s$.     
Writing $\overline{r}=\int_Xr\,d\mu_X$, we get
\[\SMALL
\int_Xv_X\diff\mu_X=\overline{r}\int_{X^r}v(\Psi_sx)\diff\mu^r(x,s)=\overline{r}\int_Mv\diff\mu_M=0.
\]
For $x,x'\in X$ with $r(x)\le r(x')$,
\[\begin{split}\SMALL
|v_X(x)-v_X(x')|&\SMALL
\le\int_0^{r(x)}|v(\Psi_sx)-v(\Psi_sx')|\diff s
+|v|_\infty|r(x')-r(x)| \\
&\SMALL\le |r|_\infty|v|_\eta d(x,x')^\eta+|v|_\infty|r|_\eta d(x,x')^\eta.
\end{split}
\]
Hence, $v_X\in C^{\eta}(X,\R^N)$ with mean zero. When $p\in[2,\infty]$, we can apply Theorem~\ref{theorem:CLT_WIP_discrete} to $v_X$, obtaining a limiting centred Brownian motion $W_X\colon[0,1]\to\R^N$. Define the new Brownian motion $W=W_X/\sqrt{\overline{r}}$.

The following result is a consequence of Theorem~\ref{theorem:CLT_WIP_discrete} applied to $v_X$ and passed to the suspension~\cite{KelMel16,MelTor04,MelZwe15}.

\begin{theorem}[WIP]\label{theorem:CLT_WIP}
Let $p\in[2,\infty]$ and  
 $v\in C^\eta( M, \R^N)$ with mean zero.
We have that
$W_n\to_wW$ as $n\to\infty$. \qed
\end{theorem}

The following theorems are the continuous time versions of Theorems~\ref{theorem:rate_wass_maps} and~\ref{theorem:rates_N=1_maps}. They display for semiflows the same rates that we have for maps.

\begin{theorem}\label{theorem:rate_wass_flow}
Let $p\in(2,3)$ and $v\in C^\eta( M, \R^N)$ with mean zero. There is $C>0$ such that
$\mathcal{W}(W_n,W)\le Cn^{-\frac{p-2}{2p}}(\log n)^\frac{p-1}{2p}
$for all integers $n>1$.
\end{theorem}

\begin{theorem}\label{theorem:rates_N=1_flow}
Let $p\in(2,\infty]$ and
$v\in C^\eta( M,\R)$ with mean zero.
There exists $C>0$ such that
\begin{empheq}[
left=
{\Pi(W_n,W)} \le C
\empheqlbrace
]{align}
&n^{-\frac{p-2}{4p}}
&& p\in(2,\infty), \label{eq:rate_Wn_p}\\
&n^{-1/4}(\log n)^{3/4} 
&& p=\infty \label{eq:rate_Wn_infty}
\end{empheq}
for all integers $n>1$.
\end{theorem}

To our knowledge, Theorems~\ref{theorem:rate_wass_flow} and~\ref{theorem:rates_N=1_flow} provide the first rates for the WIP in the dynamical systems literature for continuous time. Note that  Theorem~\ref{theorem:rate_wass_flow} implies rates in the Prokhorov metric $\Pi$ by the same argument of Remark~\ref{rem:MultidimMap}.

\section{Discrete time rates}\label{sec:results_for_maps}
In this section we provide the proofs of Theorems~\ref{theorem:rate_wass_maps} and~\ref{theorem:rates_N=1_maps}. To achieve this, we firstly recall results from~\cite{KorKosMel18} and derive some estimates in $L^\infty$. These estimates are then employed respectively with~\cite[Theorem2.3(2)]{CunDedMer20} and~\cite[Lemma 3]{Cou99} to derive our proofs.

\subsection{Approximation via martingales}\label{subsec:martingale_discrete}
We present here the relevant results from~\cite{KorKosMel18} to obtain a Gordin-type~\cite{Gor69} reversed martingale differences sequence with a control over the sum of its squares.

Let $T\colon X\to X$ be nonuniformly expanding of order~$p\in[2,\infty]$ with ergodic invariant measure $\mu_ X$. We call an \textit{extension} of $( X,T,\cB,\mu_ X)$ any measure-preserving system $(\Delta,f,\mathcal{A},\mu_\Delta)$ with a measure-preserving
$\pi_\Delta\colon\Delta\to  X$, such that 
$T\circ\pi_\Delta=\pi_\Delta\circ f$.
Denote by $P:L^1(\Delta)\to L^1(\Delta)$ the transfer operator for $f$ with respect to $\mu_\Delta$, which is characterised by $\int (Pv) w\diff\mu_\Delta=\int v(w\circ f)\diff\mu_\Delta$ for all $v\in L^1$, $w\in L^\infty$. It is well known that 
$P(v\circ f)=v $ and
$(Pv)\circ f=\mathbb{E}[v|f^{-1}\mathcal{A}]$,
for any integrable $v$.

\begin{proposition}\label{proposition:martingale_coboundary}

There is an extension 
$f\colon \Delta\to\Delta$ of 
$T\colon  X\to  X$
 such that for any
$v\in C^\eta( X, \R^N)$ with mean zero there exist 
$m\in L^p(\Delta, \R^N)$ and $\chi\in L^{p-1}(\Delta, \R^N)$, with the convention that $\infty-1=\infty$, satisfying
\begin{equation}\label{eq:martingale_coboundary}
v\circ \pi_\Delta=m+\chi\circ f-\chi,\qquad
Pm=0.
\end{equation}
If $p\in[2,\infty)$, then there exists $C>0$ such that for all $n\ge1$
\begin{equation}\label{eq:estimate_chi} \SMALL
|m|_{p}\le C\|v\|_\eta\qquad\text{and}\qquad
\bigl|\max_{1\le k\le n}|\chi\circ f^{k}-\chi|\bigr|_{p}\le C
\|v\|_\eta n^{\frac{1}{p}}.
\end{equation}
If $p=\infty$, then 
\begin{equation}\label{eq:estimate_chi_unif_exp} \SMALL
|m|_\infty\le C\|v\|_\eta\qquad\text{and}\qquad|\chi|_\infty\le C\|v\|_\eta.
\end{equation}

\begin{proof}
Equations~\eqref{eq:martingale_coboundary} and~\eqref{eq:estimate_chi} are proven in~\cite[Propositions 2.4, 2.5, 2.7]{KorKosMel18}. The estimates in~\eqref{eq:estimate_chi_unif_exp} follows by the arguments displayed before~\cite[Proposition~2.4]{KorKosMel18}. 
\end{proof}
\end{proposition}

We call $m$ the \textit{martingale part} of $v$ and $\chi$ its \textit{coboundary part}. It is relevant to cite~\cite[Corollary~2.12]{KorKosMel18} that provides the identity
$\Sigma=\int_\Delta mm^T\diff\mu_\Delta$, where
$\Sigma$ is the matrix defined in Theorem~\ref{theorem:CLT_WIP_discrete}\ref{item:Sigma_discrete}.

\begin{proposition}\label{proposition:control_sums}
Let $p\in[2,\infty)$. There exists $C>0$ such that
\begin{equation*}\label{eq:control_squares}\SMALL
\bigl|\max_{1\le k\le n}|\sum_{j=0}^{k-1}
(\mathbb{E}[mm^T-
\Sigma|f^{-1}\mathcal{A}])
\circ f^j|\bigr|_{p}
\le C\|v\|_\eta^2n^\frac12,
\end{equation*}
for every $n\ge1$.
\begin{proof}
Let 
$\breve{\Phi}=(P(mm^T))\circ f-\int_\Delta mm^T\diff\mu_\Delta$. Hence,
$\breve{\Phi}=
\E[mm^T-\Sigma|f^{-1}\mathcal{A}]$ and the result follows by~\cite[Corollary~3.2]{KorKosMel18}.
\end{proof}
\end{proposition}

\begin{proposition}\label{proposition:orthogonality}
Let $n\ge1$ and let $k,\ell\in\{0,\dots,n-1\}$ such that $k\neq\ell$. Then 
\[
\mathbb{E}[(m\circ f^k)(m\circ f^\ell)^T|f^{-n}\mathcal{A}]=0.
\]
\begin{proof}
Without loss suppose $k<\ell$. 
Hence,
\begin{equation*}\label{eq:orthogonality}
\begin{split}
\mathbb{E}[(m\circ f^k)(m\circ f^\ell)^T|f^{-n}\mathcal{A}]&=
(P^{n}[(m\circ f^k)(m\circ f^\ell)^T])\circ f^{n}\\
&=(P^{n-k}P^k[(m(m\circ f^{\ell-k})^T)\circ f^k])\circ f^{n}\\
&=(P^{n-k} [m(m\circ f^{\ell-k})^T]
)\circ f^{n}.
\end{split}
\end{equation*}
The proof is finished by 
$P[m(m\circ f^{\ell-k})^T]=(Pm)(m\circ f^{\ell-k-1})^T=0$.
\end{proof}
\end{proposition}

\begin{definition}\label{def:reversed_mg_diff_seq}
A sequence of integrable $ \R^N$-valued random variables $(d_n)_{n\ge0}$, together with the $\sigma$-algebras $(\cG_n)_{n\ge0}$, is called a \textit{ reversed martingale differences sequence} (in brief RMDS) if for all $n\ge0$ the variable $d_n$ is $\cG_n$-measurable,
$\cG_{n+1}\subseteq\cG_n$, and 
$\mathbb{E}[d_n|\cG_{n+1}]=0$.
\end{definition}

\begin{remark}\label{remark:martingale}
If $d_k$ is a RMDS with $\sigma$-algebras $\cG_k$, then for any $n\ge1$ the sequence $M_n(k)=\sum_{j=1}^{k}d_{n-j}$, $1\le k\le n$, with filtration $\mathcal{A}_k=\cG_{n-k}$, is a martingale. Namely, $M_n(k)$ is $\mathcal{A}_k$ measurable and $\E[M_n(k+1)|\mathcal{A}_k]=M_n(k)$.
\end{remark}

We now recall a classical result connecting martingale theory with measure-preserving systems. It justifies the use of the terminology "martingale part" for $m$, and it will be employed in the proofs of the rates for both maps and semiflows.

\begin{proposition}\label{proposition:mds_intro}
Let $(\Lambda,T,\cG,\mu)$ be measure-preserving system with transfer operator~$P$.
If $v\in \ker P$, then the sequence $(v\circ T^n)_{n\ge0}$ with $\sigma$-algebras $(T^{-n}\cG)_{n\ge0}$ is an RMDS.
\begin{proof}
We have that 
$T^{-(n+1)}\cG\subseteq T^{-n}\cG$ and $v\circ T^n$ is
$T^{-n}\cG$-measurable for all $n\ge0$. Using
\[\SMALL
\mathbb{E}[v\circ T^n|T^{-(n+1)}\cG]=\E[v|T^{-1}\cG]\circ T^n
\quad\text{and}\quad 
\mathbb{E}[v|T^{-1}\cG]=(Pv)\circ T,
\]
we conclude that 
$\mathbb{E}[v\circ T^n|T^{-(n+1)}\cG]=(Pv)\circ T^{n+1}=0$.
\end{proof}
\end{proposition}

The next theorem follows from~\cite[Theorem~2.3(2)]{CunDedMer20} and it is our main tool to prove multidimensional rates for the WIP.

\begin{theorem}[Cuny, Dedecker, Merlevède]\label{theorem:Cuny_et_al}
Let $p\in(2,3)$ and $N\ge1$. Suppose that $(d_n)_{n\ge0}$ is a $ \R^N$-valued stationary RMDS lying in $L^p$ with $\sigma$-algebras $(\cG_n)_{n\ge0}$. Let $M_n=\sum_{k=0}^{n-1}d_k$, $n\ge1$, and assume that
\begin{equation}\label{eq:Cuny_condition}
\sum_{n=1}^\infty\frac{1}{n^{3-p/2}}
\bigl|\mathbb{E}[M_nM_n^T|\cG_{n}]-
\mathbb{E}[M_nM_n^T]\bigr|_{p/2}<\infty.
\end{equation}
Then, there is $C>0$ and there exists a probability space supporting a sequence of random variables
$(M_n^*)_{n\ge1}$ with the same joint distributions as $(M_n)_{n\ge1}$ and a sequence $(\mathcal{N}_n)_{n\ge0}$ of i.i.d.\@ $ \R^N$-valued centered Gaussians with 
$\Var(\mathcal{N}_0)=\mathbb{E}[d_0d_0^T]$, 
such that for every integer $n>1$,
\begin{equation}\label{eq:Cuny_L1_ASIP}\SMALL
\bigl|\max_{1\le k\le n}|M_k^*-\sum_{\ell=0}^{k-1}\mathcal{N}_\ell|\bigr|_1\le C
n^{\frac{1}{p}}(\log n)^{\frac{p-1}{2p}}.
\end{equation}
\begin{proof}
This proposition is a version of~\cite[Theorem~2.3(2)]{CunDedMer20} for $p\in(2,3)$. Such a theorem is stated for a martingale differences sequence, however~\cite[Remark 2.7]{CunDedMer20} affirms that its thesis is true for \textit{reversed} martingale differences sequences, as well. To prove the sufficiency of condition~\eqref{eq:Cuny_condition}, reason as in~\cite[Remark 2.4]{CunDedMer20}.
\end{proof}
\end{theorem}

The last theorem of this subsection is taken from~\cite{Cou99} and applied to a bounded RMDS. It will be used to prove one-dimensional rates in the WIP.

\begin{theorem}[Courbot]\label{theorem:Courbot_applied}
Let $(d_n)_{n\ge0}$ with $\sigma$-algebras $(\cG_n)_{n\ge0}$ be a real bounded stationary RMDS. Consider $W\colon[0,1]\to\R$ a centred Brownian motion with variance $\sigma^2=\mathbb{E}[d_0^2]$.
Define for $1\le k\le n$ the process $\mathcal{M}_n\colon[0,1]\to\R$ as
$\mathcal{M}_n(k/n)=n^{-\frac12}\sum_{j=1}^{k}d_{n-j}$, using linear interpolation in $[0,1]$, and 
$V_{n}(k)=n^{-1}\sum_{j=1}^k\mathbb{E}[d_{n-j}^2|\cG_{n-(j-1)}]$. Let
\begin{align}
\SMALL
\kappa_n&\SMALL=
\inf\bigl\{\epsilon>0:\mathbb{P}\bigl( \max_{0\le k\le n}|V_{n}(k)-(k/n)\sigma^2|>\epsilon \bigr)\le
\epsilon\bigr\}\label{eq:kappa_n},\\
\widetilde{\kappa}_n&\SMALL=\max\bigl\{\kappa_n|\log\kappa_n|^{-\frac12}, n^{-\frac12}\bigr\}.\label{eq:tilde_kappa_n}
\end{align}
If $\lim_{n\to\infty}\kappa_n=0$, then there exists $C>0$ such that for all integers $n>1$
$$\SMALL \Pi(\mathcal{M}_n,W) \le 
C\widetilde{\kappa}_n^{1/2} |\log\widetilde{\kappa}_n|^{3/4}.$$

\begin{proof}
This theorem derives from~\cite[Lemma 3]{Cou99}, which utilizes~\cite{CoqMemVos94} to embed a continuous time martingale into a Brownian motion. It is important to note that here~$\Pi$ denotes a metric on the space of probability measures over $C([0,1])$, while in~\cite{Cou99} it is defined via the set of càdlàg functions. However, this difference does not create an obstruction, as the proof of~\cite[Lemma 3]{Cou99} can be adapted to stochastic processes with continuous sample paths, such as $\mathcal{M}_n$, thus establishing our thesis. (See also~\cite[Appendix A]{Pav23}).
\end{proof}
\end{theorem}

\subsection{Proof of Theorem~\ref{theorem:rate_wass_maps}}\label{subsec:proof}
For fixed $v\in C^\eta( X, \R^N)$  with martingale part $m\in L^p(\Delta, \R^N)$, $p\in(2,\infty)$, define the sequence of processes $X_n\colon[0,1]\to \R^N$, $n\ge1$,

\begin{equation}\label{eq:X_n}
X_{n}(k/n)=\frac{1}{\sqrt{n}}\sum_{j=0}^{k-1}m\circ f^j,
\end{equation}
for $0\le k\le n$, 
and using linear interpolation in  $[0,1]$. 
Recall that the sequence $B_n$ is defined as 
$B_{n}(k/n)=n^{-1/2}\sum_{j=0}^{k-1}v\circ T^j$ plus linear interpolation.

\begin{remark}\label{rem:p_general}
In spite of Theorem~\ref{theorem:rate_wass_maps} being valid only for $p\in(2,3)$, we work with $p\in(2,\infty)$ where possible and restrict the range only when we apply Theorem~\ref{theorem:Cuny_et_al}.
\end{remark}

\begin{lemma}\label{lem:estimate_Xn}
There exists $C>0$ such that 
$\mathcal{W}(B_n,X_n)\le Cn^{-\frac{p-2}{2p}}$ for all $n\ge1$.
\begin{proof}
By Proposition~\ref{proposition:martingale_coboundary},
\begin{align*}\SMALL
B_n(k/n)\circ\pi_\Delta-X_n(k/n)&\SMALL=
n^{-\frac12}\sum_{j=0}^{ k-1}(v\circ\pi_\Delta-m)\circ f^j
=n^{-\frac12}(\chi\circ f^{ k}-\chi)
\end{align*}
for $0\le k\le n$.
Since $B_n$ and $X_n$ are piecewise linear with the same interpolation nodes, equation~\eqref{eq:estimate_chi} yields
\begin{align*}\SMALL
\bigl|\sup_{t\in[0,1]}|B_n(t)\circ\pi_\Delta-X_n(t)|\bigr|_{p}
&\SMALL=
\bigl|\sup_{t\in\{0,\frac1n,\dots,1\}}|B_n(t)\circ\pi_\Delta-X_n(t)|\bigr|_{p}\\
&\SMALL=
n^{-\frac12}\bigl|\max_{1\le k\le n}|\chi\circ f^{k}-\chi|\bigr|_{p}\ll
n^{-\frac{p-2}{2p}}.
\end{align*}
We use that $\pi_\Delta$ is a semiconjugacy and get for any $ f\in \Lip_1$
\begin{align*}\SMALL 
\bigl|\int_ X  f(B_n)\diff\mu_ X-\int_\Delta  f(X_n)\diff\mu_\Delta\bigr|
&\SMALL\le\int_\Delta| f(B_n\circ \pi_\Delta)- f(X_n)|\diff\mu_\Delta\\
&\SMALL \le
\bigl|\sup_{t\in[0,1]}|B_n(t)\circ\pi_\Delta-X_n(t)|\bigr|_{p}
\ll n^{-\frac{p-2}{2p}},
\end{align*}
which completes the proof.
\end{proof}
\end{lemma}

We now show some estimates for i.i.d.\@ random variables and Brownian motion.

\begin{proposition}\label{proposition:max_by_exponential_tails}
Let $\{\xi_n\}_{n\ge1}$ be a sequence of identically distributed real random variables, defined on the same probability space. 
If $a=\mathbb{E}[e^{\xi_1}]<\infty$, then 
$\mathbb{E}[\max_{1\le k\le n}\xi_k]\le\log(na)$ for all $n\ge1$.
\begin{proof}
We have that 
$e^{\max_{1\le k\le n}\xi_k}=\max_{1\le k\le n}e^{\xi_k}\le\sum_{k=1}^ne^{\xi_k}$. Since all $\xi_k$ share the same distribution,
$\mathbb{E}[e^{\max_{1\le k\le n}\xi_k}]\le
\mathbb{E}[\sum_{k=1}^ne^{\xi_k}]=na$. By Jensen's inequality,
\[\SMALL
\mathbb{E}[\max_{1\le k\le n}\xi_k]\le
\log\mathbb{E}[e^{\max_{1\le k\le n}\xi_k}]\le\log(na).\tag*{\qedhere}
\]
\end{proof}
\end{proposition}

\begin{proposition}\label{proposition:exponential_tails}
Let $W\colon[0,1]\to\R^N$ be a centred Brownian motion with covariance~$\Sigma$. 
Then $\mathbb{E}[e^{\sup_{t\in[0,1]}|W(t)|}]<\infty$.

\begin{proof}
Since $\Sigma$ is symmetric and positive semidefinite, there exists an  orthogonal $N\times N$ matrix $O$ such that 
$O\Sigma O^T=\diag(\sigma_1^2,\dots,\sigma_N^2)$, with 
$\sigma_i^2\ge0$. 
Then, $OW$ is a centred Brownian motion with covariance 
$O\Sigma O^T$, and for all $1\le i\le N$ the real-valued processes $(OW)_i$ are independent centred Brownian motions with variances $\sigma_i^2$. Let $\bar{\sigma}=\max_{1\le i\le N}\sigma^2_i$. If $\bar{\sigma}=0$, then both $OW$ and $W$ are the constant zero process and the proof is finished. 
If~$\bar{\sigma}>0$, we use standard Gaussian estimates and get that for every $1\le i\le N$ there exists $C_i>0$ such that for all $s>1$
\begin{align*}\SMALL
\mathbb{P}\bigl(\sup_{t\in[0,1]}|(OW(t))_i|>s\bigr)
\le C_i\exp\bigl(-s^2/(2\bar{\sigma})\bigr).
\end{align*} 
Writing 
$\xi=\sup_{t\in[0,1]}|OW(t)|$, 
$\hat{C}=\sum_{i=1}^N C_i$, 
and $c=2N^2\bar{\sigma}$, we get
\[\SMALL
\mathbb{P}(\xi>s)\le
\sum_{i=1}^N\mathbb{P}\bigl(\sup_{t\in[0,1]}|OW_i(t)|>s/N\bigr)\le
\hat{C}\exp(-s^2/c)
\]
and
\[
\mathbb{P}(e^\xi>s)=
\mathbb{P}(\xi>\log s)\le
\hat{C}\exp\bigl(-(\log s)^2/c\bigr).
\]
By a change of variable $x=\log s$,
$$\SMALL
\mathbb{E}[e^{\xi}]
=\int_0^1 \mathbb{P}(e^{\xi}>s)\diff s+\int_1^\infty \mathbb{P}(e^{\xi}>s)\diff s
\le 1+\hat{C}\int_0^\infty e^{-x^2/c}e^{x}\diff x<\infty.$$
By orthogonality,
$|O^Tx|=|x|$ for all $x\in \R^N$. Hence,
$$|W(t)|=|O^TOW(t)|=|OW(t)|$$
for every $t\in[0,1]$. Therefore,
$\mathbb{E}[e^{\sup_{t\in[0,1]}|W(t)|}]
=\mathbb{E}[e^\xi]<\infty$.
\end{proof}
\end{proposition}

\begin{lemma}\label{lem:W1_piecewise_Brownian_motion}
Let $W\colon[0,1]\to\R^N$ be a centred Brownian motion, and let $(\mathcal{N}_n)_{n\ge0}$ be a sequence of i.i.d.\@ $ \R^N$-valued centered Gaussians with $\Var(\mathcal{N}_0)=\Var(W(1))$. Define the sequence of processes $Y_n\colon [0,1]\to \R^N$ as $Y_n(k/n)=n^{-1/2}\sum_{j=0}^{k-1}\mathcal{N}_j$ for $0\le k\le n$, with linear interpolation. Then, there exists $C>0$ such that
$\mathcal{W}(Y_{n},W)\le C n^{-\frac12}\log n$ for all integers $n>1$.
\begin{proof}
Define the sequence $Y_n^*\colon [0,1]\to \R^N$ as 
$Y_n^*(k/n)=W(k/n)$ for $0\le k\le n$, plus linear interpolation. 
We have that $Y_n=_dY_n^*$ as continuous processes for all $n\ge1$. So, for $ f\in\Lip_1$,
\[\SMALL
|\mathbb{E}[ f(Y_n)]-\mathbb{E}[ f(W)]|=
|\mathbb{E}[ f(Y_n^*)- f(W)]|\le \mathbb{E}[\sup_{t\in[0,1]}|Y_n^*(t)-W(t)|]\le A_1+A_2,
\] 
where
\[\SMALL
A_1=\mathbb{E}[\sup_{t\in[0,1]}|Y_n^*(t)-W(\lfloor nt\rfloor/n)|]
\quad\text{and}\quad
A_2=\mathbb{E}[\sup_{t\in[0,1]}|W(\lfloor nt\rfloor/n)-W(t)|].
\]
Since
\[\SMALL\begin{split}\SMALL
A_1&\SMALL=\mathbb{E}[\max_{1\le k\le n}
|W(k/n)-W((k-1)/n)|]\\
&\SMALL\le \mathbb{E}[\max_{1\le k\le n}\sup_{t\in(\frac{k-1}{n},\frac{k}{n})}
|W(t)-W((k-1)/n)|]=A_2,
\end{split}
\]
it is sufficient to estimate $A_2$.
By the rescaling property, $\widehat{W}_n(t)=n^\frac12W(t/n)$, $t\in[0,n]$ is a Brownian motion for every $n\ge1$, with the same parameters of $W$. Let $(\xi_k)_{k\ge1}$ be a identically distributed sequence of random variables with $\xi_1=_d\sup_{t\in[0,1]}|W(t)|$. Then, for every $1\le k\le n$,
\[\begin{split}\SMALL
\sup_{t\in(\frac{k-1}{n},\frac{k}{n})}
|W(t)-W((k-1)/n)|&\SMALL=
n^{-\frac12}\sup_{t\in(\frac{k-1}{n},\frac{k}{n})}
 |\widehat{W}_n(nt)-\widehat{W}_n(n(k-1))|\\
&\SMALL=n^{-\frac12}\sup_{t\in(k-1,k)}
 |\widehat{W}_n(t)-\widehat{W}_n(k-1)|=_dn^{-\frac12}\xi_k.
\end{split}
\]
Proposition~\ref{proposition:exponential_tails} yields 
$\mathbb{E}[e^{\xi_1}]<\infty$, hence we can apply Proposition~\ref{proposition:max_by_exponential_tails} getting that
$\SMALL
A_2= n^{-\frac12}\mathbb{E}[
\max_{1\le k\le n}\xi_k]\ll n^{-\frac12}\log n$,
which completes the proof.
\end{proof}
\end{lemma}

We present here the proof of Theorem~\ref{theorem:rate_wass_maps}, which we obtain by Theorem~\ref{theorem:Cuny_et_al} in combination with Lemmas~\ref{lem:estimate_Xn} and~\ref{lem:W1_piecewise_Brownian_motion}.

\begin{proof}[\textbf{Proof of Theorem~\ref{theorem:rate_wass_maps}}]
Let $p\in(2,3)$, $X_n$ be from~\eqref{eq:X_n}, and let $Y_n$ be as in Lemma~\ref{lem:W1_piecewise_Brownian_motion}. Consider $W$ an $N$-dimensional Brownian motion with mean $0$ and covariance $\Sigma$ from Theorem~\ref{theorem:CLT_WIP_discrete}\ref{item:WIP_discrete}. Recall that $\Sigma=\int_\Delta mm^T\diff\mu_\Delta$.
By Lemmas~\ref{lem:estimate_Xn} and~\ref{lem:W1_piecewise_Brownian_motion}, to prove the rate on $\mathcal{W}(B_n,W)$ it suffices to estimate $\mathcal{W}(X_n,Y_n)$.

Let us check the hypotheses of Theorem~\ref{theorem:Cuny_et_al} for $d_n=m\circ f^n$.
By~\eqref{eq:martingale_coboundary} we have that $Pm=0$ , so Proposition~\ref{proposition:mds_intro} yields that $d_n$  with $\sigma$-algebras $f^{-n}\mathcal{A}$ is an RMDS on the probability space~$(\Delta,\mathcal{A}, \mu_\Delta)$. It lies in $L^p$ by Proposition~\ref{proposition:martingale_coboundary}, and it is  stationary because $f$ is measure-preserving. 
Let us show that $M_n=\sum_{j=0}^{n-1}m\circ f^j$, $n\ge1$, satisfies condition~\eqref{eq:Cuny_condition}.
In the following equation, the off-diagonal terms are zero by Proposition~\ref{proposition:orthogonality}:
\begin{align*}\SMALL 
\mathbb{E}[M_nM_n^T|f^{-n}\mathcal{A} ]-&
\mathbb{E}[M_nM_n^T]\\
=&\SMALL \sum_{k,\ell=0}^{n-1}\bigl( 
\mathbb{E}[(m\circ f^k)(m\circ f^\ell)^T|f^{-n}\mathcal{A}]-
\mathbb{E}[(m\circ f^k)(m\circ f^\ell)^T]
\bigr)\\
=&\SMALL \sum_{k=0}^{n-1} \bigl( 
 \mathbb{E}[(mm^T)\circ f^k|f^{-n}\mathcal{A}]-
\mathbb{E}[(mm^T)\circ f^k]\bigr)\\
=& \SMALL \mathbb{E}\bigl[\sum_{k=0}^{n-1}(mm^T-\Sigma)\circ f^k
|f^{-n}\mathcal{A}\bigr].
\end{align*}
Using Proposition~\ref{proposition:control_sums},
\begin{align*}
\SMALL
\bigl|\mathbb{E}\big[\sum_{k=0}^{n-1}(mm^T-\Sigma)\circ f^k
|f^{-n}\mathcal{A}\bigr]\bigr|_{p/2}
&\SMALL
= \bigl|\mathbb{E}\bigl[ \sum_{k=0}^{n-1}
\mathbb{E}[(mm^T-\Sigma)\circ f^k|f^{-k-1}\mathcal{A}]
\bigl|f^{-n}\mathcal{A}\bigr]\bigr|_{p/2}\\
&\SMALL 
\le \bigl| \sum_{k=0}^{n-1}
\mathbb{E}[(mm^T-\Sigma)\circ f^k|f^{-k-1}\mathcal{A}]\bigr|_{p/2}\\
&\SMALL 
=\bigl|\sum_{k=0}^{n-1}
\mathbb{E}[(mm^T-\Sigma)|f^{-1}\mathcal{A}]\circ f^k
\bigr|_{p/2}\ll n^{\frac12}.
\end{align*}
Since $p\in(2,3)$, the series in~\eqref{eq:Cuny_condition} converges.
 
By Theorem~\ref{theorem:Cuny_et_al}, there exists a probability space supporting a sequence $(M^*_n)_{n\ge1}$ with the same joint distributions as $(M_n)_{n\ge1}$ and a sequence $(\mathcal{N}_n)_{n\ge0}$ of i.i.d.\@ $ \R^N$-valued centered Gaussians with 
$\Var(\mathcal{N}_0)=\mathbb{E}[mm^T]=\Sigma$, such that~\eqref{eq:Cuny_L1_ASIP} holds.

Let $M^*_0=0$ and define for $n\ge1$ the process $X_n^*\colon[0,1]\to \R^N$ as 
$X^*_n(k/n)=n^{-\frac12}M^*_{k}$
for $0\le k\le n$, with linear interpolation.
So, $X^*_n=_dX_n$ as continuous processes.
Using~\eqref{eq:Cuny_L1_ASIP}, we have that for all $ f\in\Lip_1$,
\[
\begin{split}\SMALL
\mathcal{W}(X_n,Y_n)&\SMALL
\le \mathbb{E}[ f(X^*_n)- f(Y_n)]
\le \mathbb{E}[\sup_{t\in[0,1]}|X_n^*(t)-
Y_n(t)|]\\
&\SMALL=n^{-\frac12}
\bigl|\max_{1\le k\le n}|M^*_{k}-\sum_{\ell=0}^{k-1}\mathcal{N}_\ell|\bigr|_1
\ll n^{-\frac{p-2}{2p}}(\log n)^\frac{p-1}{2p}.
\end{split}
\]
Hence,
$\mathcal{W}(B_n,W)
\ll n^{-\frac{p-2}{2p}}(\log n)^\frac{p-1}{2p}$ 
and the proof is complete.
\end{proof}

\subsection{Using bounded martingales}\label{subsec:using_bounded_martingales}

Let $T$ be nonuniformly expanding of order~$p=\infty$. 
For $v\in C^\eta( X,\R)$ with mean zero, we consider its martingale part $m\in L^\infty(\Delta)$  and write  
$\breve{\Phi}=\mathbb{E}[m^2|f^{-1}\mathcal{A}]-\sigma^2$, where
$\sigma^2=\int_\Delta m^2\diff\mu_\Delta$. As pointed out before~\cite[Corollary~3.2]{KorKosMel18}, there exist $\breve{m},\breve{\chi}\colon\Delta\to\R$ with $P\breve{m}=0$ such that 
$\breve{\Phi}=\breve{m}+\breve{\chi}\circ f-\breve{\chi}$, which is called the \textit{secondary martingale-coboundary decomposition} of $v$.
Since $\tau\in\ L^\infty$, \cite[Proposition~3.1]{KorKosMel18} and the arguments displayed before~\cite[Proposition~2.4]{KorKosMel18} yield that there is $C>0$ such that 
\begin{equation}\label{eq:estimate_chi_unif_exp_secondary} \SMALL
|\breve{m}|_\infty\le C\|v\|^2_\eta\qquad\text{and}\qquad
|\breve{\chi}|_\infty\le C\|v\|^2_\eta.
\end{equation}

\begin{proposition}[Azuma-Hoeffding inequality~{\cite[pg 237]{Wil91}}]\label{proposition:azuma-hoeff}
Let $M(n)= \sum_{j=1}^{n} X_{j}$, $n\ge1$, 
be a real-valued martingale with $X_j\in L^\infty$ for $j\ge1$. Then
\[
\mathbb{P}\Bigl(\max_{1\le k\le n}|M(k)|\ge x\Bigr)\le
\exp\biggl\{\frac{-x^2/2}{\sum_{j=1}^n|X_j|_\infty^2}\biggr\},
\]
for every $x\ge0$ and $n\ge1$.\qed
\end{proposition}

For the rest of this section, for any $k\ge1$ and $g\colon \Delta\to\R$ we write $g_k=\sum_{j=0}^{k-1}g\circ f^j$.
\begin{proposition}\label{proposition:control_sums_unif_exp}
Let $v\in C^\eta( X,\R)$ with mean zero. 
There exist $a,C>0$ such that
\begin{equation*}\label{eq:control_squares_unif_exp}
\mu_\Delta\Bigl(\max_{1\le k\le n}|{\SMALL\sum_{j=0}^{n-1}
\breve{\Phi}
\circ f^j}|\ge x\Bigr)
\le C\exp\biggl\{-\frac{ax^2}{n}\biggr\},
\end{equation*}
for every $x>0$ and $n\ge1$.
\begin{proof} 
Let $\breve{\Phi}=\breve{m}+\breve{\chi}\circ f-\breve{\chi}$. 
Hence,
$\breve{\Phi}_k=\breve{m}_k+\breve{\chi}\circ f^k-\breve{\chi}$, and  
by~\eqref{eq:estimate_chi_unif_exp_secondary} 
there is $K>0$ such that for any $n\ge1$ we have 
$\max_{1\le k\le n}|\breve{\Phi}_k|\le \max_{1\le k\le n}|\breve{m}_k|+K$. 
So,
\begin{equation}\label{eq:estimate_bounded}\begin{split}\SMALL
\mu_\Delta(\max_{1\le k\le n}|\breve{\Phi}_k|\ge x)&\SMALL\le
\mu_\Delta(\max_{1\le k\le n}|\breve{m}_k|+K \ge x)\\
&\SMALL\le \mu_\Delta(\max_{1\le k\le n}|\breve{m}_k|\ge x/2) +
\mu_\Delta(K\ge x/2).
\end{split}
\end{equation}

If $m=0$, we have automatically  
$\mu_\Delta(\max_{1\le k\le n}|\breve{m}_k|\ge x/2)=0$. If $m\ne 0$, 
we use that $P\breve{m}$=0 to get from Proposition~\ref{proposition:mds_intro} that
$(\breve{m}\circ f^n)_{n\ge0}$ is an RMDS on the probability space $(\Delta,\mu_\Delta)$. 
As in Remark~\ref{remark:martingale}, for every $n\ge1$ the process 
$\breve{M}_n(k)=\sum_{j=1}^{k}\breve{m}\circ f^{n-j}$, $1\le k\le n$,
is a martingale.
Since 
$\breve{m}_k=\breve{M}_n(n)-\breve{M}_n(n-k)$,
using Proposition~\ref{proposition:azuma-hoeff} and~\eqref{eq:estimate_chi_unif_exp_secondary}, there is $c>0$ such that
\begin{equation*}\label{eq:Azuma-H_sec_dec}
\begin{split}
\mu_\Delta(\max_{1\le k\le n}|\breve{m}_k|\ge x/2)&\le
\mu_\Delta(\max_{1\le k\le n}|\breve{M}_n(k)|\ge x/4)\\
&\le \exp\biggl\{\frac{-x^2/32}{\sum_{j=1}^n|\breve{m}|_\infty^2}\biggr\}=
\exp\biggl\{-\frac{cx^2}{n}\biggr\}.
\end{split}
\end{equation*}
Since 
$\mu_\Delta(K\ge x/2)=1$ for $x\le 2K$ and $0$ otherwise,
\begin{equation*}\label{eq:azuma_coboundary}
\mu_\Delta(K\ge x/2)\le
\exp\{4K^2-x^2\}\le \exp\{4K^2\}\exp\{-x^2/n\}.
\end{equation*}
Conclude by applying these estimates to~\eqref{eq:estimate_bounded}.
\end{proof}
\end{proposition}

\subsection[]{Proof of Theorem~\ref{theorem:rates_N=1_maps} ($p=\infty$)}\label{subsec:theorem:rates_N=1_maps}

For fixed $v\in C^\eta( X,\R)$ with mean zero and martingale part $m\in L^\infty(\Delta)$, define the sequence of processes $Y_n\colon[0,1]\to\R$, $n\ge1$
\begin{equation*}
Y_n(k/n)=\frac{1}{\sqrt{n}}\sum_{j=1}^{k}m\circ f^{n-j},
\end{equation*}
for $1\le k\le n$, using linear interpolation in  $[0,1]$. 
Following~\cite[Lemma 4.8]{KelMel16}, let $h\colon C([0,1],\R)\to C([0,1],\R)$ be the linear operator $(hf)(t)= f(1)-f(1-t)$.

\begin{lemma}\label{lem:estimate_Bn_Yn}
There exists $C>0$ such that 
$\Pi(h\circ B_n\circ\pi_\Delta,Y_n)\le C n^{-\frac12}$ for all $n\ge1$.
\begin{proof}
The process $h\circ B_n$ is piecewise linear on $[0,1]$ with interpolation nodes $k/n$ for $0\le k\le n$, and
$h\circ B_n(k/n)\circ\pi_\Delta=n^{-\frac12}\sum_{j=n-k}^{n-1} v\circ f^j$.
By~\eqref{eq:martingale_coboundary},
\begin{align*}
h\circ B_n(k/n)\circ \pi_\Delta- Y_n(k/n)&=\SMALL
n^{-\frac12}\bigl(\sum_{j=n-k}^{n-1} v\circ\pi_\Delta\circ f^j -
 \sum_{j=1}^{k} m\circ f^{n-j}\bigr) \\
&= n^{-\frac12}((v\circ \pi_\Delta)_n-(v\circ\pi_\Delta)_{n-k}
-(m_n-m_{n-k}))\\
&=n^{-\frac12}(\chi\circ f^n-\chi\circ f^{n-k}).
\end{align*}
Since $h\circ B_n\circ \pi_\Delta$ and $Y_n$ have the same interpolation nodes, we have by~\eqref{eq:estimate_chi_unif_exp}, 
\[\begin{split}\SMALL
\bigl|\sup_{t\in[0,1]}|h\circ B_n(t)\circ \pi_\Delta- Y_n(t)|\bigr|_\infty
&\SMALL=\bigl|\max_{0\le k\le n}|h\circ B_n(k/n)\circ \pi_\Delta- Y_n(k/n)|\bigr|_\infty
\\
&\SMALL\le 2n^{-\frac12}|\chi|_\infty\ll  n^{-\frac12}.    
\end{split}
\]
Since the Prokhorov metric is bounded by the infinity norm, 
\[\SMALL
\Pi(h\circ B_n\circ \pi_\Delta,Y_n)\le
\bigl|\sup_{t\in[0,1]}|h\circ B_n(t)\circ \pi_\Delta- Y_n(t)|\bigr|_\infty
\ll n^{-\frac12}.\tag*{\qedhere}
\]
\end{proof}
\end{lemma}

\begin{lemma}\label{lem:Courbot}
There is $C>0$ such that 
$\Pi(Y_n,W)\le Cn^{-\frac14}(\log n)^{\frac34}$ for all integers $n>1$.
\begin{proof}
Following the proof of Theorem~\ref{theorem:rate_wass_maps}, the sequence $d_n=m\circ f^n$, $n\ge0$, with $\sigma$-algebras $f^{-n}\mathcal{A}$ is a stationary RMDS on the probability space~$(\Delta,\mathcal{A}, \mu_\Delta)$. 
Equation~\eqref{eq:estimate_chi_unif_exp} yields that $d_n$ is bounded.
We adopt the same notation of Theorem~\ref{theorem:Courbot_applied}, noting that $\sigma^2=\int_\Delta m^2\diff\mu_\Delta$ and that $Y_n$ coincides with $\mathcal{M}_n$. 
We have that
$$\SMALL V_n(k)=n^{-1}\sum_{j=1}^k\E[m^2\circ f^{n-j}|f^{-n-(j-1)}\mathcal{A}]=
n^{-1}\sum_{j=1}^k\E[m^2|f^{-1}\mathcal{A}]\circ f^{n-j}.$$ 

We claim that
\[\SMALL
\kappa_n\ll \sqrt{n^{-1}\log n}.
\]
Assuming the claim true, let us evaluate
$\widetilde{\kappa}_n$ from~\eqref{eq:tilde_kappa_n}. 
Note that $x\mapsto x^2(\log x)^{-1}$ is decreasing for $x\in(0,1)$. Hence
$x\mapsto x^2|\log x|^{-1}$ 
is increasing and so is
$x\mapsto x|\log x|^{-\frac12}$.
Since $\kappa_n\ll \sqrt{n^{-1}\log n}$, we get that
\begin{equation*}\label{eq:kappa_estimate}
\kappa_n|\log \kappa_n|^{-\frac12}
\ll\sqrt{\frac{\log n}{n|\log\log n-\log n|}}
\ll\frac{1}{\sqrt{n}}.
\end{equation*}
By definition, $\widetilde{\kappa}_n\ll n^{-\frac12}$ as well, and the statement follows from Theorem~\ref{theorem:Courbot_applied}.\vspace{1ex}

Let us now prove the claim. 
Writing 
$\breve{\Phi}=\mathbb{E}[m^2|f^{-1}\mathcal{A}]-\sigma^2$ and 
$\breve{\Phi}_k=\sum_{j=0}^{k-1}\breve{\Phi}\circ f^j$, 
\begin{equation*}\SMALL
V_n(k)-(k/n)\sigma^2=n^{-1}\sum_{j=1}^k\breve{\Phi}\circ f^{n-j}=
n^{-1}(\breve{\Phi}_n-\breve{\Phi}_{n-k}),
\end{equation*}
for every $n\ge1$. 
So, 
$\max_{0\le k\le n}|V_n(k)-(k/n)\sigma^2|\le
2 n^{-1} \max_{1\le k\le n}|\breve{\Phi}_k|$. 
By Proposition~\ref{proposition:control_sums_unif_exp}, there are $a,C>0$ such that
\[\SMALL
\mu_\Delta\bigl(\max_{0\le k\le n}|V_n(k)-(k/n)\sigma^2|\ge \epsilon\bigr)\le
\mu_\Delta\bigl(\max_{1\le k\le n}|\breve{\Phi}_k|\ge n\epsilon/2\bigr)
\le Ce^{-an\epsilon^2},
\]
for all $\epsilon>0$ and $n\ge1$. 
Let now 
$\epsilon_n=\sqrt{\log n/(an)}$. We have that 
$C\le n\epsilon_n$ for $n$ large enough and
\[\SMALL
\mu_\Delta\bigl(\max_{0\le k\le n}|V_n(k)-(k/n)\sigma^2|>\epsilon_n\bigr)\le
C\exp\{-an\epsilon_n^2\}=
C/n\le\epsilon_n.
\]
By definition~\eqref{eq:kappa_n}, we have that 
$\kappa_n\ll\epsilon_n\ll \sqrt{n^{-1}\log n}$ which proves the claim.
\end{proof}
\end{lemma}

The proof of the next proposition is part of the proof of~\cite[Theorem 2.2]{AntMel19}.
However, we present this result as a separate statement due to its multiple applications and to maintain our work self-contained.

\begin{proposition}\label{proposition:final_inequality}
Let $Z(t)$, $t\in[0,1]$, be a $\R^N$-valued continuous process with $Z(0)=0$ a.s.\@ and let $W(t)$, $t\in[0,1]$, be a $N$-dimensional Brownian motion. 
Then, 
$$\Pi(Z,W)\le 2\Pi(h\circ Z,W).$$
\begin{proof}
It is easy to see that 
$h\circ W=_dW$.
Note that $h(hf)=f$ if $f(0)=0$, and the map 
$h\colon  C([0,1],\R)\to C([0,1],\R)$ 
is Lipschitz with constant $\Lip(h)\le2$. We conclude by the Lipschitz mapping theorem~\cite[Theorem~3.2]{Whi74},
\[\SMALL
\Pi(Z,W)=
\Pi(h(h\circ Z),h(h\circ W))
\le 2\Pi(h\circ Z,h\circ W)=
2\Pi(h\circ Z,W).\tag*{\qedhere}
\]
\end{proof}
\end{proposition}

\begin{proof}[\textbf{Proof of Theorem~\ref{theorem:rates_N=1_maps} ($p=\infty$)}]
 Since $B_n(0)=0$ for all $n\ge1$, applying Proposition~\ref{proposition:final_inequality} with $N=1$ we get
\[
\Pi(B_n,W)\ll
\Pi(h\circ B_n,W)
\le \Pi(h\circ B_n, Y_n)+ \Pi(Y_n,W).\]
Since $\Pi(h\circ B_n, Y_n)=\Pi(h\circ B_n\circ\pi_\Delta, Y_n)$, we apply Lemmas~\ref{lem:estimate_Bn_Yn} and~\ref{lem:Courbot} to finish.
\end{proof}

\section{Martingale-coboundary decompositions for semiflows}\label{sec:new_dec}

In this section we recall how the semiflows defined in Subsection~\ref{subsec:nonunif_exp_semiflow} can be modelled by suspensions over a Gibbs-Markov map with an unbounded roof function.  
We show that Hölder observables on the ambient space lift to regular functions on the suspension, for which we get two new decompositions  in the style of Gordin~\cite{Gor69}. 
This follows and extends the approach of~\cite{KorKosMel18} to continuous time.

We inform the reader that sometimes in our arguments it may be necessary to diminish the parameter~$\eta\in(0,1]$;  such a change does not create any issue because of the inclusion of Hölder spaces.

\subsection{Gibbs-Markov semiflows}
Let $\Psi_t\colon  M\to M$ be a nonuniformly expanding semiflow of order~$p\in[2,\infty]$ as in Subsection~\ref{subsec:nonunif_exp_semiflow}.  Hence, there exist $\eta\in(0,1]$, a set $X\subset M$ with a Borel probability measure $\rho$ and a first return function $r\in C^\eta(X,\R)$ with $\inf r\ge1$.
 Since the map $T=\Psi_r\colon X\to X$, is nonuniformly expanding, following Subsection~\ref{subsec:nonunif_exp_maps} there are
a subset $Y\subset X$ with $\rho(Y)>0$ and  a measurable partition~$\{Y_j\}_{j\ge1}$, a return time $\tau\in L^p(Y,\rho)$, and a map $F=T^{\tau}\colon Y\to Y$ such that conditions~\ref{item:full_branch}-\ref{item:zeta} are satisfied. 
It is a standard result that there is a unique absolutely continuous ergodic $F$-invariant probability measure $\mu$ on $Y$ and 
 $d \mu/(d\restr{\rho}{Y})$ is bounded, and hence $\tau\in L^p(Y,\mu)$.

Let $\phi\colon Y\to[1,\infty)$ be defined as $\phi(y)=\sum_{j=0}^{\tau(y)-1}r(T^jy)$. It lies in $L^p(Y,\mu)$ because $\phi\le |r|_\infty \tau$. 
Define the suspension 
$Y^\phi=\{(y,u)\in Y\times[0,\infty):u\in [0,\phi(y)]\}/\sim$
where $(y,\phi(y))\sim (Fy,0)$.
The suspension semiflow $F_t:Y^\phi\to Y^\phi$ is given by
$F_t(y,u)=(y,u+t)$ computed modulo identifications.

There is a semiconjugacy
$\pi_M\colon Y^\phi\to M$ between $F_t$ and $\Psi_t$, defined as $\pi_M(y,u)=\Psi_uy$. We have the ergodic $F_t$-invariant probability measure
$\mu^\phi=(\mu\times{\rm Lebesgue})/\bar\phi$,
where $\bar\phi=\int_Y\phi\,d\mu$. Moreover, it is easy to see that $\mu_M=(\pi_M)_*\mu^\phi$ is exactly the ergodic $\Psi_t$-invariant probability measure on $M$, defined in Subsection~\ref{subsec:nonunif_exp_semiflow}.

\begin{proposition}\label{proposition:phi}
There is $C>0$ such that 
\begin{equation}\label{eq:phi}
\SMALL
|\phi(y)-\phi(y')|\le 
C(\inf_{Y_j}\phi) d(Fy,Fy')^{\eta},
\end{equation}
for every $j\ge1$ and $y,y'\in Y_j$. 
Moreover, if $p\neq\infty$ then
\begin{equation}\label{eq:polynomial_tails}\SMALL
\sum_{j}\mu(Y_j)(\SMALL\sup_{Y_j}\phi^p)<\infty.
\end{equation}
\begin{proof}

Recall that $\tau$ is constant on partition elements.
Using that $r\in C^\eta(X,\R)$, 
\begin{equation*}\SMALL
|\phi(y)-\phi(y')|\le
\sum_{\ell=0}^{\tau(y)-1}| r(T^\ell y)-r(T^\ell y')|
\le
|r|_{\eta}\sum_{\ell=0}^{\tau(y)-1} d(T^\ell y,T^\ell y')^{\eta}, 
\end{equation*}
for each $j\ge1$ and $y,y'\in Y_j$. Hence, point~\ref{item:non_expansion} from Subsection~\ref{subsec:nonunif_exp_maps} yields that there is $C>0$ such that
\[\SMALL
|\phi(y)-\phi(y')|\le C\tau(y) d(Fy,Fy')^{\eta}.
\]
By $\inf_{X} r\ge1$ and the definition of $\phi$, we get $\restr{\tau}{Y_j}\le\restr{\phi}{Y_j}$, which implies $\tau(y)\le (\inf_{Y_j}\phi)$. Equation~\eqref{eq:phi} follows.

By~\eqref{eq:phi}, we get
$\sup_{Y_j}\phi-\inf_{Y_j}\phi\le C\diam(Y)^\eta(\SMALL\inf_{Y_j}\phi)$. 
Hence, there exists $K>0$ such that 
$\sup_{Y_j}\phi\le K\inf_{Y_j}\phi$ 
for all $j\ge1$. So, 
\begin{align*}\SMALL
\sum_{j}\mu(Y_j)(\SMALL\sup_{Y_j}\phi^p)
\le K^p\sum_j\mu(Y_j)({\SMALL \inf_{Y_j}}\phi^p)
\le K^p|\phi|_p^p<\infty. \tag*{\qedhere} 
\end{align*}
\end{proof}
\end{proposition}

Since $F\colon Y\to Y$ is Gibbs-Markov and $\phi\colon Y\to[1,\infty)$ satisfies~\eqref{eq:phi}, we say as in~\cite[Definition~2.2]{BalButMel19} that $F_t\colon Y^\phi\to Y^\phi$ is a \textit{Gibbs-Markov semiflow}.
We define the space $C^\eta(Y,\R^N)$ with norm $\|\cdot\|_\eta$ similarly to Definition~\ref{def:Hölder}.

\begin{definition}[Function space on $Y^\phi$]\label{def:functions_on_Y^phi}
Let $\eta\in(0,1]$, $N\ge1$, and define for $j\ge1$ the set $Y^\phi_j=\{(y,u)\in Y^\phi:y\in Y_j\}$.
For $v:Y^\phi\to \R^N$, let
$|v|_{\infty}=\sup_{(y,u)\in Y^\phi} |v(y,u)|$
and
\[
\|v\|_{\eta}=|v|_{\infty}+|v|_{\eta},\qquad 
|v|_{\eta}=\sup_{j\ge1}\sup_{(y,u),(y',u)\in Y_j^\phi,\ y\neq y'} \frac{|v(y,u)-v(y',u)|}{(\inf_{Y_j}\phi)d(Fy,Fy')^\eta }.
\]
Let $\mathcal{F}^\eta(Y^\phi, \R^N)$ consist of observables $v:Y^\phi\to \R^N$ with
$\|v\|_{\eta}<\infty$. 
\end{definition}

Our next result shows that Hölder observables on the ambient space $M$ lift naturally into the just defined function space on $Y^\phi$.

\begin{proposition}\label{proposition:v_in_F}
Let $v\in C^\eta( M, \R^N)$. Then $w=v\circ\pi_M\in\mathcal{F}^{\eta^2}(Y^\phi, \R^N)$ and there exists $C>0$ such that 
$\|w\|_{\eta^2}\le C\|v\|_\eta$.

\begin{proof}
We clearly have that $|w|_\infty\le |v|_\infty$. 
Hence, we are left to show that $|w|_{\eta^2}\ll\|v\|_\eta$. 
Let $j\ge1$ and $(y,u), (y',u)\in Y_j^\phi$. We have
\begin{equation}\label{eq:v_in_F}\SMALL
|w(y,u)-w(y',u)|=|v(\Psi_uy)-v(\Psi_uy')|
\le|v|_\eta d(\Psi_uy,\Psi_uy')^\eta.
\end{equation}

For $r\in C^\eta(X,\R)$ and $k\ge1$, we write $r_k=\sum_{j=0}^{k-1}r\circ T^j$. 
Let $n\ge0$ be such that $r_n(y)\le u<r_{n+1}(y)$. 
So, $n\le\tau(y)$ and $u=r_n(y)+E(y)$, where 
$E(y)\le r(T^ny)\le|r|_\infty$. By~\eqref{eq:continuous_dependence}, there is $K>0$ (dependent on $|r|_\infty$) such that
\begin{equation*}\label{eq:control_flow}\SMALL
d(\Psi_uy,\Psi_uy')=d(\Psi_{E(y)}(\Psi_{r_n(y)}y),\Psi_{E(y)}(\Psi_{r_n(y)}y'))
\le K
d(\Psi_{r_n(y)}y,\Psi_{r_n(y)}y').
\end{equation*}
Using~\eqref{eq:Lipschitz}, 
\begin{align*}\SMALL
d(\Psi_{r_n(y)}y,\Psi_{r_n(y)}y')&\SMALL\le
d(\Psi_{r_n(y)}y,\Psi_{r_n(y')}y')+d(\Psi_{r_n(y')}y',\Psi_{r_n(y)}y')\\
&\SMALL=d(T^ny,T^ny')+d(\Psi_{r_n(y')}y',\Psi_{r_n(y)}y')\\
&\SMALL\ll d(T^ny,T^ny')+
L|r_n(y)-r_n(y')|.
\end{align*}
By point~\ref{item:non_expansion} of Subsection~\ref{subsec:nonunif_exp_maps},
\[\SMALL d(T^ny,T^ny')\ll d(Fy,Fy')
\le\diam(M)^{1-\eta}d(Fy,Fy')^{\eta}.
\]
Using again~\ref{item:non_expansion},
\[\SMALL
|r_n(y')-r_n(y)|\le |r|_{\eta}\sum_{j=0}^{n-1}d(T^jy,T^jy')^{\eta}\ll
n d(Fy,Fy')^{\eta}.
\]
By $n\le\tau(y)\le \inf_{Y_j}\phi$, we get 
$|r_n(y')-r_n(y)|\ll(\inf_{Y_j}\phi)d(Fy,Fy')^{\eta}$. 
Combining these estimates, there exists $C>0$ such that
\[\SMALL
d(\Psi_uy,\Psi_uy')\le C(\inf_{Y_j}\phi)d(Fy,Fy')^{\eta}.
\]
Combining the above equation with~\eqref{eq:v_in_F}, we get
\begin{equation*}\label{eq:general_function_space} \SMALL
|w(y,u)-w(y',u)|\le C|v|_\eta(\inf_{Y_j}\phi)^\eta d(Fy,Fy')^{\eta^2}.
\end{equation*}
By $\inf\phi\ge1$, we get 
$(\inf_{Y_j}\phi)^\eta\le\inf_{Y_j}\phi$ and so $|w|_{\eta^2}\le C|v|_\eta$.
\end{proof}
\end{proposition}

We conclude this subsection with two  estimates and introducing the transfer operators for both semiflow $F_t\colon Y^\phi\to Y^\phi$ and map $F\colon Y\to Y$. Let $g=\diff\mu/(\diff\mu\circ F)$ be the inverse Jacobian of $F$. 
Then, (see for example~\cite{AarDen01}) there is $C>0$  such that
\begin{equation}\label{eq:g}\SMALL
g(y)\le C\mu(Y_j)\qquad\text{and}\qquad|g(y)-g(y')|\le 
C\mu(Y_j) d(Fy,Fy')^{\eta},
\end{equation}
for all $j\ge1$ and $y,y'\in Y_j$. 

We denote with $L_t:L^1(Y^\phi)\to L^1(Y^\phi)$ the transfer operator for $F_t$, so $\int (L_tv) w\diff\mu^\phi=\int v (w\circ F_t)\diff\mu^\phi$ for all $v\in L^1$, $w\in L^\infty$, $t\ge0$. 
Let $P:L^1(Y)\to L^1(Y)$ be the transfer operator for $F$, so $\int (Pv) w\diff\mu=\int v(w\circ F)\diff\mu$ for all $v\in L^1$ and $w\in L^\infty$; recall that $|Pv|_q\le |v|_q$ for all $q\in[1,\infty]$.
The pointwise formula for $P$ is given by $(Pv)(y)=\sum_j g(y_j)v(y_j)$ where $y_j$ is the unique preimage of $y$ under $\restr{F}{Y_j}$.

\subsection{Primary decomposition}\label{subsec:primary_dec}
We start by describing the class of functions that will admit our new decomposition. 

\begin{definition}\label{def:star}
For $v\colon Y^\phi\to\R^N$, define $v'\colon Y\to \R^N$ as $v'(y)=\int_0^{\phi(y)}v(y,u)\diff u$.
We say that $v$ satisfies $(\star)$ if (i)~$v\in L^\infty(Y^\phi)$, (ii)~$\int_{Y^\phi} v\diff\mu^\phi=0$ and (iii)~$\|Pv'\|_\eta<\infty$.  
For such functions $v$, we write  $\langle v\rangle_\eta=|v|_\infty+\|Pv'\|_\eta$.
\end{definition}

Note that if $v$ satisfies $(\star)$ then  $Pv'$ has mean zero. 
This follows by
\[\SMALL
\int_Y Pv' \diff\mu=\int_Yv'\diff\mu
=\int_Y\int_0^{\phi(y)}v(y,u)\diff u\diff\mu
=\bar{\phi}\int_{Y^\phi}v\diff\mu^\phi=0.
\]

\begin{proposition}\label{proposition:Pholder} There exists $C>0$ such that
$\|Pv'\|_\eta\le C\|v\|_\eta$ for all 
$v\in\mathcal{F}^\eta(Y^\phi, \R^N)$. 
If in addition $\int_{Y^\phi}v\diff\mu^\phi=0$, then $v$ satisfies $(\star)$ and $\langle v\rangle_\eta\le C\|v\|_\eta$.

\begin{proof} 		
Let $y,y'\in Y_j$ and suppose without loss that $\phi(y)\le\phi(y')$. By~\eqref{eq:phi},
\begin{equation}\label{eq:preimages1}
\begin{split}
|v'(y)-v'(y')|
\le&\SMALL\int_0^{\phi(y)}|v(y,u)-v(y',u)|\diff u
+\int_{\phi(y)}^{\phi(y')}|v(y',u)|\diff u\\
\ll&\SMALL \bigl(|v|_\eta(\inf_{Y_j}\phi) (\sup_{Y_j}\phi)
+|v|_\infty(\inf_{Y_j}\phi)\bigr) d(Fy,Fy')^\eta\\
\le& \|v\|_\eta (\SMALL\sup_{Y_j}\phi^2)  d(Fy,Fy')^\eta.
\end{split}
\end{equation}
	
Let now $y,y'\in Y$, with preimages $y_j,y'_j\in Y_j$ under $F$. Since $|v'|\le \phi|v|_\infty$, we have that 
$|v'(y_j)|\le |v|_\infty(\sup_{Y_j}\phi)$.
Using~\eqref{eq:g},~\eqref{eq:preimages1}, and~\eqref{eq:polynomial_tails} with $p=2$	
\begin{align*}\SMALL
|(Pv')(y)-(Pv')(y')|
\le&\SMALL\sum_j|g(y_j)
-g(y_j)||v'(y'_j)|
+\sum_j g(y'_j)|v'(y_j)-v'(y'_j)|\\
\ll&\SMALL \|v\|_\eta 
\bigl(\sum_j\mu(Y_j)(\sup_{Y_j}\phi^{2})\bigr)d(Fy_j,Fy'_j)^\eta
\ll\|v\|_\eta d(y,y')^\eta.
\end{align*}
Similarly,~\eqref{eq:polynomial_tails} yields also that 
$|Pv'|_\infty\ll|v|_\infty$, concluding that
$\|Pv'\|_\eta\ll\|v\|_\eta$.
\end{proof}
\end{proposition}

Let $v\colon Y^\phi\to\R^N$ satisfy $(\star)$. We define $\chi',m'\colon Y\to \R^N$ as follows:
\begin{equation*}\label{eq:decompY}\SMALL
\chi'=\sum_{k=1}^\infty P^kv',
\qquad m'=v'-\chi'\circ F+\chi'.
\end{equation*}
It is well known for Gibbs-Markov maps (see~\cite[Theorem~1.6]{AarDen01}), that for every function $w\in C^\eta(Y,\R^d)$ with mean zero, there are $a,C>0$ such that 
$\|P^kw\|_\eta\le Ce^{-ak}$ for all $k\ge1$.
Since $Pv'\in  C^\eta(Y,\R^d)$ has mean zero,  the series 
$$\SMALL\sum_{k=1}^\infty \|P^kv'\|_\eta=\sum_{k=0}^\infty \|P^kPv'\|_\eta$$ converges.
By completeness,  
$\chi'\in C^\eta(Y,\R^d)$ and  
$Pm'=Pv'-\chi'+\sum_{k=2}^\infty P^kv'=0$. 
It follows that
\begin{equation}\label{eq:def_m'_chi'}\SMALL
\|\chi'\|_\eta\le\sum_{k=0}^\infty\|P^kPv'\|_\eta
\ll\|Pv'\|_\eta,
\qquad
|m'|_p\le|\phi|_p|v|_\infty+2|\chi'|_\infty
\ll\langle v\rangle_\eta.
\end{equation}
Hence $m'\in L^p(Y, \R^N)$.

Such a construction of $\chi'$ and $m'$ is conducted in the same way in~\cite[Subsection 2.2]{KorKosMel18}; we decided to include it to make our argument self-contained.

Define $m,\chi\colon Y^\phi\to \R^N$ for $y\in Y$ and $u\in[0,\phi(y))$ by
\begin{equation}\label{eq:chiandm}
\chi(y,u) = \chi'(y) + \int_0^u v(y,s) \diff s,
\qquad m(y,u) =
\begin{cases}
m'(y) & u\in[\phi(y)-1,\phi(y)), \\
0 & u\in[0,\phi(y)-1).
\end{cases}
\end{equation}

\begin{proposition}\label{proposition:m_and_chi}
We have that   
$m\in L^p(Y^\phi, \R^N)$ and $\chi\in L^{p-1}(Y^\phi, \R^N)$, with the convention that $\infty-1=\infty$.
Moreover, there exists $C>0$ such that
\begin{equation*}
|m|_p\le C \langle v\rangle_\eta
\qquad\text{and}\qquad  
|\chi|_{p-1}\le \langle v\rangle_\eta,
\end{equation*}
for all functions $v$ that satisfy $(\star)$.

\begin{proof}
Firstly, suppose that $p=\infty$. Then, by~\eqref{eq:def_m'_chi'} and~\eqref{eq:chiandm},
\begin{equation*}
|\chi|_\infty\le|\chi'|_\infty+|\phi|_\infty|v|_\infty
\ll\langle v\rangle_\eta.
\qquad
|m|_\infty=|m'|_\infty\ll\langle v\rangle_\eta.
\end{equation*}
Secondly, suppose that $p\in[2,\infty)$. By~\eqref{eq:def_m'_chi'}, 
$|\chi(y,u)|\le |\chi'|_\infty+u|v|_\infty\ll\phi(y)\langle v\rangle_\eta$.
Hence,
\begin{equation*}\label{eq:p-1_norm_chi}\SMALL
|\chi|_{p-1}\ll
\langle v\rangle_\eta(\int_Y\int_0^\phi|\phi|^{p-1}\diff s\diff \mu)^{\frac{1}{p-1}}=
\langle v\rangle_\eta|\phi|_p^\frac{p}{p-1}<\infty.
\end{equation*}
Since $m'\in L^p(Y, \R^N)$,~\eqref{eq:chiandm} and~\eqref{eq:def_m'_chi'} yield
\begin{equation*}\label{eq:p_norm_m}\SMALL
|m|_{p}^p\ll\int_Y\int_0^\phi
|m'|^p\mathbbm{1}_{\{\phi-1\le u<\phi\}}\diff u \diff\mu=
|m'|_p^p\ll\langle v\rangle_\eta^p<\infty.
\end{equation*}
The statement follows.
\end{proof}

\end{proposition}

Our next proposition shows how the transfer operator $L_1$ acts pointwise.

\begin{proposition}\label{proposition:L1}
Let $v\in L^1(Y^\phi)$. Then
\[(L_1v)(y,u)=
\begin{cases}
v(y,u-1)&u\in[1,\phi(y))\\
\sum_jg(y_j)v(y_j,u-1+\phi(y_j))&u\in[0,1),
\end{cases}
\]
where $y_j$ is the unique preimage of $y$ under $\restr{F}{Y_j}$.

\begin{proof}
Let $w\in L^\infty(Y^\phi)$. 
By definition of $L_1$ and $\mu^\phi$, using the substitution $u\mapsto u+1$,
\begin{equation}\label{eq:prima}
\begin{split}\SMALL
\int_{Y^\phi}L_1(\mathbbm{1}_{\{0\le u<\phi-1\}}v)w\diff\mu^\phi&\SMALL
=\bar{\phi}^{-1}\int_Y\int_0^{\phi(y)} \mathbbm{1}_{\{0\le u<\phi(y)-1\}}v(y,u)w(y,u+1)\diff u\diff\mu\\
&\SMALL=\int_{Y^\phi}\mathbbm{1}_{\{1\le u<\phi(y)\}}v(y,u-1)w(y,u)\diff\mu^\phi.
\end{split}
\end{equation}
Next, let us focus on 
$\mathbbm{1}_{\{\phi-1\le u<\phi\}}v$. 
By the substitution $u\mapsto u+1-\phi(y)$,
\begin{equation*}
\begin{split}\SMALL
\int_{Y^\phi}L_1(\mathbbm{1}_{\{\phi-1
\le u<\phi\}}v)w\diff\mu^\phi&\SMALL
=\bar{\phi}^{-1}\int_Y\int_{\phi(y)-1}^{\phi(y)} v(y,u)w(Fy,u+1-\phi(y))\diff u\diff\mu\\
&\SMALL=\bar{\phi}^{-1}\int_Y\int_0^1 v(y,u-1+\phi(y))w(Fy,u)\diff u\diff\mu.
\end{split}
\end{equation*}
Write 
$\tilde{v}_u(y)=v(y,u-1+\phi(y))$ 
and $w^u(y)=w(y,u)$. Then,
\begin{equation}\label{eq:seconda}
\begin{split}\SMALL
\int_{Y^\phi}L_1(\mathbbm{1}_{\{\phi-1\le u<\phi\}}v)w\diff\mu^\phi&\SMALL
=\bar{\phi}^{-1}\int_0^1\int_Y \tilde{v}_u  (w^u\circ F)\diff\mu\diff u\\
&\SMALL=\bar{\phi}^{-1}\int_0^1\int_Y (P\tilde{v}_u)  w^u\diff\mu\diff u\\
&\SMALL=\int_{Y^\phi} 
\mathbbm{1}_{\{0\le u<1\}}(P\tilde{v}_u) w\diff\mu^\phi.
\end{split}
\end{equation}
Equations~\eqref{eq:prima} and~\eqref{eq:seconda} yield that
\begin{align*}\SMALL
(L_1v)(y,u)
&\SMALL=L_1(\mathbbm{1}_{\{0\le u<\phi-1\}}v
+ \mathbbm{1}_{\{\phi-1\le u<\phi\}}v)(y,u)\\
&\SMALL=\mathbbm{1}_{\{1\le u<\phi\}}v(y,u-1)
+\mathbbm{1}_{\{0\le u<1\}}(P\tilde{v}_u)(y).
\end{align*}
The proof is completed by the pointwise formula for $P$.
\end{proof}
\end{proposition}

We present now our new primary martingale-coboundary decomposition.

\begin{proposition}\label{proposition:primarydec}
Suppose that $v$ satisfies $(\star)$ and define $\psi\colon Y^\phi\to \R^N$ as $\psi=\int_0^1v\circ F_s\diff s$. Then
$\psi=m+\chi\circ F_1-\chi$ and $m\in\ker L_1$.

\begin{proof}
Let $(y,u)\in Y^\phi$ with 
$u\in[0,\phi(y)-1)$. 
Then $F_1(y,u)=(y,u+1)$ and 
$\psi(y,u)=\int_u^{u+1}v(y,s)\diff s$, so
\[\SMALL
\chi(F_1(y,u))-\chi(y,u)
=\int_0^{u+1}v(y,s)\diff s-\int_0^uv(y,s)\diff s
=\psi(y,u)=\psi(y,u)-m(y,u).
\]
If $u\in[\phi(y)-1,\phi(y))$, then 
\[\SMALL
\psi(y,u)=\int_0^{u+1-\phi(y)}v(Fy,s)\diff s+v'(y)
-\int_0^uv(y,s)\diff s.
\]
We have that $F_1(y,u)=(Fy,u+1-\phi(y))$. 
By definition, $v'-m'=\chi'\circ F-\chi'$ and 
$m(y,u)=m'(y)$, so
\begin{align*}\SMALL
\chi(F_1(y,u))-\chi(y,u)&\SMALL
=\chi'(Fy)-\chi'(y)+\int_0^{u+1-\phi(y)}v(Fy,s)\diff s
-\int_0^uv(y,s)\diff s\\
&\SMALL=v'(y)-m'(y)+\psi(y,u)-v'(y)=\psi(y,u)-m(y,u).
\end{align*}
Therefore $\psi=m+\chi\circ F_1-\chi$ 
on the whole of $Y^\phi$.

We are left to prove that $m\in\ker L_1$ using the formula of Proposition~\ref{proposition:L1}. 
Let $y\in Y$. If $u\in[1,\phi(y))$, then 
$u-1\in[0,\phi(y)-1)$ and by definition of $m$,
\[\SMALL
(L_1m)(y,u)=m(y,u-1)=0.
\]
If $u\in[0,1)$, then 
$u-1+\phi(y_j)\in[\phi(y_j)-1,\phi(y_j))$ 
for all preimages $y_j$ of $y$, and
\[\SMALL
(L_1 m)(y,u)
=\sum_jg(y_j)m(y_j,u-1+\phi(y_j))
=(Pm')(y)=0,
\]
because $m'\in\ker P$.
\end{proof}
\end{proposition}
Following the terminology of Section~\ref{sec:results_for_maps}, the new functions $m$ and $\chi$ are called respectively martingale and coboundary part of $v$.
In view of Proposition~\ref{proposition:primarydec}, to estimate the Birkhoff sums of $\psi$ in $p$-norm, it would be desirable to have $\chi\in L^p$. This is indeed true for $p=\infty$ by Proposition~\ref{proposition:m_and_chi}; however, in general $\chi$ lies in $L^{p-1}$. 
Nevertheless, the next result show that for $p\in[2,\infty)$ the function $\chi\circ F_n-\chi$ lies in $L^p$ for all $n\ge1$.

\begin{proposition}\label{proposition:coboundary_p_norm}
There exists $C>0$ such that
$|\max_{1\le k\le n}|\chi\circ F_k-\chi||_p
\le C\langle v\rangle_\eta n^{1/p}$ for all $n\ge1$. Moreover,
\begin{equation}\label{eq:coboundary}\SMALL
\bigl|\max_{1\le k\le n}|\chi\circ F_k-\chi|\bigr|_p\le 
C\langle v\rangle_\eta\bigl(n^{1/q}+n^{1/p}|\mathbbm{1}_{\{\phi\ge n^{1/q}\}}\phi|_p\bigr)
\end{equation}
for all $n\ge1$, $q\ge p$, and $v$ satisfying $(\star)$.

\begin{proof}
This proof is identical to the one of~\cite[Proposition~2.7]{KorKosMel18}, with the obvious notational changes. (See also~\cite[Proposition~3.37]{Pav23}).
\end{proof}
\end{proposition}

The next Corollary is found~\cite[Corollary~2.8]{KorKosMel18} and we prove it for completeness. 

\begin{corollary}\label{corollary:coboundary_p_norm}
$|\max_{1\le k\le n}|\chi\circ F_k-\chi||_p=o(n^{1/p})$.
\begin{proof} 
Using that $\phi\in L^p(Y)$, we have  
$|\mathbbm{1}_{\{\phi\ge n^{1/q}\}}\phi|_p\to0$ 
by the bounded convergence theorem. 
Let $q>p$, then equation~\eqref{eq:coboundary} yields for $n\to\infty$ that
\begin{equation*}\SMALL
n^{-1/p} \bigl|\max_{1\le k\le n}|\chi\circ F_k-\chi|\bigr|_p\ll
 n^{-\frac{q-p}{pq}}+|\mathbbm{1}_{\{\phi\ge n^{1/q}\}}\phi|_p\longrightarrow 0.\tag*{\qedhere}
\end{equation*}
\end{proof}
\end{corollary}

\subsection{Key estimates}
This section displays results from~\cite{KorKosMel19} to get crucial estimates for the martingale-coboundary decomposition of any $v\colon Y^\phi\to\R^N$ that satisfies $(\star)$. We denote $m$ and $\chi$ for respectively the martingale and the coboundary part of $v$, as defined in~\eqref{eq:chiandm}

\begin{proposition}\label{proposition:Rio_v}
Let $p\in[2,\infty)$. 
There exists $C>0$ such that
\begin{equation}\label{eq:Burk_applied}\SMALL
\bigl|\max_{1\le k\le n}|{\SMALL\sum_{j=0}^{k-1}m\circ F_j}|\bigr|_{p}\le 
C\langle v\rangle_\eta n^\frac12,
\end{equation}
and
\begin{equation}\label{eq:Rio_applied}\SMALL
\bigl|\max_{1\le k\le n}|{\SMALL\int_0^kv\circ F_s\diff s}|\bigr|_{2(p-1)}\le 
C\langle v\rangle_\eta n^\frac12,
\end{equation}
for all $n\ge1$ and any $v$ that satisfies $(\star)$.

\begin{proof}
This proof is carried out exactly as the one of~\cite[Corollary~2.10]{KorKosMel18} with the obvious notational changes. 
We remark that an essential ingredient for~\eqref{eq:Burk_applied} is Burkholder's inequality~\cite{Bur73}, while~\eqref{eq:Rio_applied} follows from  Rio's inequality~\cite{KuePhi80}. (See also~\cite[Proposition~3.42]{Pav23}).
\end{proof}
\end{proposition}

\begin{corollary}\label{corollary:variance}
The limit $\Sigma=\lim_{n\to\infty}n^{-1}\int_{Y^\phi}(\int_0^nv\circ F_s\diff s)
(\int_0^nv\circ F_s\diff s)^T\diff\mu^\phi\in\R^{N\times N}$ exists for any $v$ satisfying $(\star)$. 
Moreover,
$\Sigma=\int_{Y^\phi}mm^T\diff\mu^\phi$.
\begin{proof}
This proof is identical to the one of~\cite[Corollary~2.12]{KorKosMel18}, with the obvious notational changes for the suspension $Y^\phi$ and applying Proposition~\ref{proposition:Rio_v} and Corollary~\ref{corollary:coboundary_p_norm} to the martingale-coboundary decomposition of $v$.
(See also~\cite[Corollary~3.43]{Pav23}).
\end{proof}
\end{corollary}

\begin{corollary}[WIP]\label{corollary:WIP} 
Let $\widehat{W}_n(t)=n^{-1/2}\int_0^{nt}v\circ F_s\diff s$, $t\in[0,1]$, and let $W\colon[0,1]\to\R^N$ be a centered Brownian motion with covariance $\Sigma$ from Corollary~\ref{corollary:variance}. Then, $\widehat{W}_n\to_wW$ on the probability space $(Y^\phi,\mu^\phi)$.

\begin{proof}
This proof is carried by the same approach of~\cite[Corollary~2.13]{KorKosMel18}, applying Proposition~\ref{proposition:Rio_v} with the obvious notational changes.
\end{proof}
\end{corollary}

\begin{proposition}\label{proposition:azuma_applied}
Let $p=\infty$ and let $v$ satisfy $(\star)$. There exist $a,C>0$ such that
\[
\mu^\phi\Bigl(\max_{1\le k\le n}|{\SMALL\int_0^kv\circ F_j}|\ge x\Bigr)\le
C\exp\biggl\{-\frac{ax^2}{n}\biggr\},
\]
for all $n\ge1$ and $x\ge0$.
\begin{proof}
Proposition~\ref{proposition:primarydec} yields that
$\int_0^1v\circ F_s\diff s
=m+\chi\circ F_1-\chi$, with $m\in L^\infty\cap\ker L_1$. Then,  $(m\circ F_n)_{n\ge1}$ is an RMDS
by Proposition~\ref{proposition:mds_intro} and it is bounded. We have
$\int_0^{n}v\circ F_s\diff s= \sum_{j=0}^{n-1}m\circ F_j+\chi\circ F_n-\chi$ for $n\ge1$. To conclude, reason similarly to the proof of Proposition~\ref{proposition:control_sums_unif_exp}, replacing $\breve{\Phi}$, $\breve{m}$, $\breve{\chi}$, $f$, $\mathcal{A}$ with  respectively $\psi$, $m$, $\chi$, $F_1$, $\cB$, and using the estimates from Proposition~\ref{proposition:primarydec} instead of~\eqref{eq:estimate_chi_unif_exp_secondary}.
\end{proof}
\end{proposition}

\subsection{Secondary decomposition}\label{subsec:secondary_dec}

Let $v\in\mathcal{F}^\eta(Y^\phi,\R^N)$ with mean zero and let $m$ and $\chi$ be respectively its martingale and coboundary parts, as defined in~\eqref{eq:chiandm}. By Proposition~\ref{proposition:Pholder}, the observable $v$ satisfies $(\star)$ and $\langle v\rangle_\eta\ll\|v\|_\eta$.
Let $U_Ff=f\circ F$, $f\in L^1(Y)$, and $U_1g=g\circ F_1$, $g\in L^1(Y^\phi)$, be the Koopman operators respectively for $F$ and $F_1$.
\begin{proposition}\label{proposition:v_breve}
$(U_1L_1(mm^T))(y,u)=
\begin{cases}
(U_FP(m'm'^{T})(y)&u\in[\phi(y)-1,\phi(y))\\
0&u\in[0,\phi(y)-1)
\end{cases}$

\begin{proof}
Let $(y,u)\in Y^\phi$. 
By Proposition~\ref{proposition:L1} and definition of $m$, if $u\in[1,\phi(y))$,
\begin{equation}\label{eq:first}
(L_1(mm^T))(y,u)=mm^T(y,u-1)=0;
\end{equation}
and if $u\in[0,1)$
\begin{equation}\label{eq:second}\SMALL
(L_1(mm^T))(y,u)
=\sum_jg(y_j)mm^T(y_j,u-1+\phi(y_j))
=(P(m'm'^{T}))(y).
\end{equation}
Let us analyse $U_1L_1(mm^T)$. 
If $(y,u)\in Y^\phi$ is such that $u\in[0,\phi(y)-1)$, 
then $u+1\in[1,\phi(y))$ and by~\eqref{eq:first} we get 
\[
(U_1L_1(mm^T))(y,u)=(L_1(mm^T))(y,u+1)=0.
\]
If $u\in[\phi(y)-1,\phi(y))$, then 
$u+1-\phi(y)\in[0,1)$ 
and~\eqref{eq:second} yields that
\[
(U_1L_1(mm^T))(y,u)
=(L_1(mm^T))(Fy,u+1-\phi(y))
=(P(m'm'^{T}))(Fy),
\]
finishing the proof.
\end{proof}
\end{proposition}
As in Corollary~\ref{corollary:variance}, let $\Sigma=\int mm^T\diff\mu^\phi$ and define
\begin{equation}\label{eq:vbreve}
\breve{v}=U_1L_1(mm^T)-\Sigma=\E[mm^T-\Sigma|F^{-1}_1\cB].
\end{equation}
By $\int U_1L_1(mm^T)\diff\mu^\phi
=\int L_1(mm^T)\diff\mu^\phi=\Sigma$ it follows that
$\int_{Y^\phi}\breve{v}\diff\mu^\phi=0$. Following Definition~\ref{def:star}, we write 
$\breve{v}'(y)=\int_0^{\phi(y)}\breve{v}(y,u)\diff u$. 

\begin{proposition}\label{proposition:Pholderbis} There exists $C>0$ such that 
\[
|\breve{v}|_\infty
\le C\|v\|_\eta^2\qquad\text{and}\qquad\|P\breve{v}'\|_\eta\le C\|v\|^2_\eta,
\]
for all $v\in\mathcal{F}^\eta(Y^\phi, \R^N)$ with mean zero. 
So, the observable $\breve{v}$ satisfies $(\star)$ and $\langle\breve{v}\rangle_\eta\le\|v\|_\eta^2$.

\begin{proof}

By definition of $m'$, we see that 
$|m'|\le |v'| +2|\chi'|_\infty$. 
Using $|v'|\le\phi|v|_\infty$ and that
$\|\chi'\|_\eta\ll\|Pv'\|_\eta$ from~\eqref{eq:def_m'_chi'}, 
we get that
$|m'|\ll\phi\langle v\rangle_\eta$ and $|m'm'^T|\le\phi^2\langle v\rangle_\eta^2$. 
By~\eqref{eq:g} and~\eqref{eq:polynomial_tails} with $p=2$, we get that for all $y\in Y$
\begin{equation}\label{eq:Pm'm'_infinity}\SMALL
|P(m'm'^T)(y)|\le\sum_jg(y_j)|m'm'^T(y_j)|\ll
\sum_j\mu(Y_j)(\SMALL \sup_{Y_j}\phi^2)\ \langle v\rangle_\eta^2\ll\langle v\rangle_\eta^2.
\end{equation}
Recall that Proposition~\ref{proposition:m_and_chi} yields $|m|_2\ll\langle v\rangle_\eta$.  
By Proposition~\ref{proposition:v_breve} and~\eqref{eq:Pm'm'_infinity}, 
\[\SMALL
|\breve{v}|_\infty\le|U_FP(m'm'^T)|_\infty+\bigl|\int_{Y^\phi}mm^T\diff\mu^\phi\bigr|\le
|P(m'm'^T)|_\infty+|m|_2^2\ll\langle v\rangle_\eta^2.
\]
The first estimate follows by $\langle v\rangle_\eta\ll\|v\|_\eta$.

Let us now show the second estimate. Proposition~\ref{proposition:v_breve} yields
\[\SMALL
\breve{v}'(y)
=\int_0^{\phi(y)}(U_FP(m'm'^{T})(y)\mathbbm{1}_{\{\phi(y)-1<u<\phi(y)\}}
-\Sigma)\diff u
=(U_FP(m'm'^{T}))(y)-\phi(y)\Sigma.
\]
The identity $PU_F=\Id_{L^1(Y)}$ implies that 
$P\breve{v}'=P(m'm'^{T})-(P\phi)\Sigma$. 
Therefore, to complete the proof it suffices to show that 
$\|P(m'm'^{T})\|_\eta\ll\|v\|_\eta^2$ and 
$\|(P\phi)\Sigma\|_\eta\ll\|v\|_\eta^2$.\vspace{1ex}

Let us focus on $(P\phi)\Sigma$.
Apply Proposition~\ref{proposition:Pholder} with 
$v\equiv 1$ 
to get that $v'=\phi$ 
and $\|P\phi\|_\eta\ll1$. 
Hence,
$\|(P\phi)\Sigma\|_\eta
=\|P\phi\|_\eta|\Sigma|
\ll|m|_2^2\ll\langle v\rangle_\eta
\ll\|v\|_\eta^2$.

Next, let us focus on $P(m'm'^{T})$. 
We already know by~\eqref{eq:Pm'm'_infinity} that
$|P(m'm'^{T})|_\infty\ll\|v\|_\eta^2.$
Let $y,y'\in Y_j$. 
By definition of $m'$, equation~\eqref{eq:preimages1} and 
$\chi'\in C^\eta(Y,\R^N)$, we get
\begin{align*}\SMALL
|m'(y)-m'(y')|&\SMALL
\le|v'(y)-v'(y')|+|\chi'(Fy)-\chi'(Fy')|
+|\chi'(y)-\chi'(y')|\\
&\SMALL \ll \|v\|_\eta (\sup_{Y_j}\phi)\  d(Fy,Fy')^\eta
+\|\chi'\|_\eta  d(Fy,Fy')^\eta
+\|\chi'\|_\eta  d(y,y')^\eta.
\end{align*}
Since $\|\chi'\|_\eta\ll\|Pv'\|_\eta\ll\|v\|_\eta$, we use point~\ref{item:expansion} of Subsection~\ref{subsec:nonunif_exp_maps} to get
$|m'(y)-m'(y')|\ll
\|v\|_\eta (\SMALL\sup_{Y_j}\phi) d(Fy,Fy')^\eta$. 
Using that 
$|m'|\ll\phi\langle v\rangle_\eta\ll\phi\|v\|_\eta$,  we obtain
\begin{equation*}
\begin{split}
|m'(y)m'(y)^T-m'(y')m'(y')^T|
\le&(|m'(y)|+|m'(y')|) |m'(y)-m'(y')|\\
\ll&\|v\|_\eta^2 (\SMALL\sup_{Y_j}\phi^2) d(Fy,Fy')^\eta.
\end{split}
\end{equation*}
Fix $y,y'\in Y$ with preimages $y_j,y'_j\in Y_j$ under $F$. By~\eqref{eq:g} and~\eqref{eq:polynomial_tails} with $p=2$,
\begin{align*}\SMALL
|(P(m'm'^{T}))(y)-(P(m'm'^{T}))(y')|&\SMALL
\le\sum_j|g(y_j)
-g(y'_j)||(m'm'^{T})(y_j)|\\
&\SMALL\qquad+\sum_j g(y_j) |(m'm'^{T})(y_j)-(m'm'^{T})(y'_j)|\\
&\SMALL\ll\|v\|_\eta^2\sum_j\mu(Y_j)(\sup_{Y_j}\phi^2)\
d(Fy_j,Fy'_j)^\eta\\
&\SMALL\ll \|v\|^2_\eta d(y,y')^\eta.
\end{align*}
We conclude that $\|P(m'm'^{T})\|_\eta\ll\|v\|_\eta^2$.
\end{proof}
\end{proposition}

By Proposition~\ref{proposition:Pholderbis}, for any $v\in\mathcal{F}^\eta(Y^\phi,\R^N)$ we can apply Proposition~\ref{proposition:primarydec} to $\breve{v}$, obtaining the \textit{secondary martingale-decomposition} of $v$. 

We show now that the Birkhoff sum and integral of $\breve{v}$ are close. 
For $n\ge1$, define
$\phi_k=\sum_{j=0}^{k-1}\phi\circ F^j$.
For $(y,u)\in Y^\phi$ and $t>0$, define the \textit{lap number} $N_t(y,u)=n\ge0$ to be the unique integer such that 
$\phi_n(y)\le t+u<\phi_{n+1}(y)$.
\begin{proposition}\label{proposition:intandsum}
There exists $C>0$ such that
\begin{equation*}\label{eq:they_are_close}\SMALL
\bigl|\int_0^n\breve{v}\circ F_s\diff s-\sum_{j=0}^{n-1}\breve{v}\circ F_j\bigr|_\infty
\le C\|v\|_\eta^2,
\end{equation*}
for every $n\ge1$ and 
$v\in\mathcal{F}^\eta(Y^\phi, \R^N)$ with mean zero.

\begin{proof}
Define 
$\alpha=U_FP(m'm'^T)$. 
Proposition~\ref{proposition:v_breve} gives that for all $(y,u)\in Y^\phi$,
$$(U_1L_1(mm^T))(y,u)
=\alpha(y)\ \mathbbm{1}_{\{\phi(y)-1\le u<\phi(y)\}}.$$  
The integral 
$\int_0^n(U_1L_1(mm^T))\circ F_s\diff s$ 
sums $\alpha$ along an orbit under $F$, with an error given by 
\begin{equation}\label{eq:integralestimate}\SMALL
\bigl|\int_0^n(U_1L_1(mm^T))(F_s(y,u))\diff s
-\sum_{j=0}^{N_{n-1}(y,u)}\alpha(F_jy)\bigr|
\le |\alpha(y)|
+|\alpha(F^{N_n(y,u)}y)|\le2|\alpha|_\infty,
\end{equation}
for all $n\ge1$ and $(y,u)\in Y^\phi$.	
	
We find that every initial point $(y,u)\in Y^\phi$ enters the strip $[\phi-1,\phi)$ exactly once every lap. 
Still, the sum 
$\sum_{j=0}^{n-1}(U_1L_1(mm^T))\circ F_j$ 
 could miss the term 
$\alpha\circ F^{N_{n-1}}$, 
giving that for every $(y,u)\in Y^\phi$ and all $n\ge1$,
\begin{equation}\label{eq:sumestimate}\SMALL
\bigl|\sum_{j=0}^{n-1}(U_1L_1(mm^T))(F_j(y,u))
-\sum_{j=0}^{N_{n-1}(y,u)}\alpha(F_jy)\bigr|
\le |\alpha(F^{N_{n-1}(y,u)}y)|\le|\alpha|_\infty.
\end{equation}
	
Both~\eqref{eq:integralestimate} and~\eqref{eq:sumestimate} can be restated with infinity norms, because the estimates are uniform in $(y,u)$. 
Combine~\eqref{eq:integralestimate} and~\eqref{eq:sumestimate}, noticing that the two terms $n\Sigma$ cancel out:
\begin{align*}\SMALL
\bigl|\int_0^n\breve{v}\circ F_s\diff s
-\sum_{j=0}^{n-1}\breve{v}\circ F_j\bigr|_\infty&\SMALL
=\bigl|\int_0^n(U_1L_1(mm^T))\circ F_s\diff s
-\sum_{j=0}^{{n-1}}(U_1L_1(mm^T))\circ F_j\bigr|_\infty\\
&\SMALL\le 3|\alpha|_\infty\le3|P(m'm'^T)|_\infty\ll\langle v\rangle_\eta^2,
\end{align*}
where the last inequality is true by~\eqref{eq:Pm'm'_infinity}. Conclude by $\langle v\rangle_\eta\ll\|v\|_\eta$.
\end{proof}
\end{proposition}

Proposition~\ref{proposition:intandsum} enables us to connect discrete and continuous time estimates. 
The two corollaries that follow will be used in the proofs of our main theorems.

\begin{corollary}\label{corollary:Burkholder_v_prime}
Let $p\in[2,\infty)$. There exists $C>0$ such that
\[\SMALL
\bigl|\max_{1\le k\le n}|{\SMALL\sum_{j=0}^{k-1}\breve{v}\circ F_j}|\bigr|_{2(p-1)}\le 
C\|v\|^2_\eta n^\frac12,
\]
for all $n\ge1$ and 
$v\in\mathcal{F}^\eta(Y^\phi, \R^N)$ with mean zero.
\begin{proof}
Since $\breve{v}$ satisfies $(\star)$ by Proposition~\ref{proposition:Pholderbis}, we get from~\eqref{eq:Rio_applied} that
\[\SMALL|\max_{1\le k\le n}|{\SMALL\int_0^k\breve{v}\circ F_j}||_{2(p-1)}
\ll \langle \breve{v}\rangle_\eta n^\frac12
\ll \|v\|_\eta^2n^\frac12.\] 
The statement follows by Proposition~\ref{proposition:intandsum}.
\end{proof}
\end{corollary}

\begin{corollary}\label{corollary:azuma_v_breve}
Let $p=\infty$ and $v\in \mathcal{F}^\eta(Y^\phi,\R)$ with mean zero. There exist $a,C>0$ such that
\[
\mu^\phi\Bigl(\max_{1\le k\le n}|{\SMALL\sum_{j=0}^{k-1}\breve{v}\circ F_j}|\ge x\Bigr)\le
C\exp\biggl\{-\frac{a\epsilon^2}{n}\biggr\},
\]
for all $x\ge0$ and $n\ge1$.
\begin{proof}
By Proposition~\ref{proposition:intandsum}, there exists $K>0$ such that 
\[\SMALL
\max_{1\le k\le n}\bigl|{\SMALL\sum_{j=0}^{k-1}\breve{v}\circ F_j}\bigr|\le 
\max_{1\le k\le n}\bigl|{\SMALL\int_0^k\breve{v}\circ F_j}\bigr|+
K.
\]
Hence,
\begin{equation*}
\mu^\phi\Bigl(\max_{1\le k\le n}|{\SMALL\sum_{j=0}^{k-1}\breve{v}\circ F_j}|\ge x\Bigr)\le
\mu^\phi\Bigl(\max_{1\le k\le n}|{\SMALL\int_0^k\breve{v}\circ F_j}|\ge x/2\Bigr)+
\mu^\phi(K\ge x/2).
\end{equation*}
Since $\breve{v}$ satisfies $(\star)$ by Proposition~\ref{proposition:Pholderbis}, the first term of the right-hand side is sorted by Propositions~\ref{proposition:azuma_applied}. The second term is treated as in~\eqref{eq:azuma_coboundary}.
\end{proof}
\end{corollary}

\section{Continuous time rates}\label{sec:continuous_time_rates}
This section provides the proofs of Theorems~\ref{theorem:rate_wass_flow} and~\ref{theorem:rates_N=1_flow}. 
Let $\Psi_t\colon  M\to  M$ be a nonuniformly expanding semiflow of order~$p\in(2,\infty]$ with ergodic invariant measure~$\mu_M$ as in Subsection~\ref{subsec:nonunif_exp_semiflow}, and let $v\in C^\eta( M, \R^N)$ with mean zero. 
For $t\in[0,1]$ and $n\ge1$, let $W_n(t)=n^{-\frac12}\int_0^{nt}v\circ \Psi_s\diff s$ be as in~\eqref{eq:Wn}, which converges weakly to a centred $N$-dimensional Brownian motion $W$ by Theorem~\ref{theorem:CLT_WIP}. 
Let $F_t\colon Y^\phi\to Y^\phi$ be the respective Gibbs-Markov semiflow with ergodic invariant measure $\mu^\phi$, which is semiconjugated to $\Psi_t$ by the map $\pi_M\colon Y^\phi\to M$, $\pi_M(y,u)=\Psi_uy$. 
Then the observable $w=v\circ\pi_M$ has mean zero, and by Proposition~\ref{proposition:v_in_F} it  lies in $\mathcal{F}^{\eta^2}(Y^\phi, \R^N)$. 
We define the sequence of processes $\widehat{W}_n$ on the probability space $(Y^\phi,\mu^\phi)$ as 
$\widehat{W}_n=W_n\circ\pi_M$. 
We have that \begin{equation}\label{eq:Wn'}
\widehat{W}_n(t)=\frac{1}{\sqrt{n}}\int_0^{nt}w\circ F_s\diff s,\qquad t\in[0,1].
\end{equation}
Following Corollary~\ref{corollary:WIP}, we note that the limiting Brownian motion $W$ has covariance matrix 
$\Sigma=\int_{Y^\phi}mm^T\diff\mu^\phi$, 
where $m$ is the martingale part of $w$.
Since $\pi_M$ is measure-preserving, we have
$W_n=_d\widehat{W}_n$ for all $n\ge1$, so
$\mathcal{W}(W_n,W)=\mathcal{W}(\widehat{W}_n,W)$ and
$\Pi(W_n,W)=\Pi(\widehat{W}_n,W)$.
Therefore, in the remainder of this section we deal with observables $w$ that lie in 
$\mathcal{F}^\eta(Y^\phi, \R^N)$ for some $\eta\in(0,1]$, and prove rates for $\widehat{W}_n$.\vspace{1ex} 

We recall from Proposition~\ref{proposition:primarydec} that there exist 
$m,\chi\colon Y^\phi\to \R^N$ such that 
\begin{equation}\label{eq:primary_dec}\SMALL
\int_0^1w\circ F_s\diff s=m+\chi\circ F_1-\chi.
\end{equation}
Since $\langle w\rangle_\eta\ll\|w\|_\eta$, Proposition~\ref{proposition:m_and_chi} yields that there is $C>0$ such that, for $p\in(2,\infty]$
\begin{equation}\label{eq:estimate_m_chi}
|m|_p\le C\|w\|_\eta,\qquad  
|\chi|_{p-1}\le C\|w\|_\eta.
\end{equation}
Proposition~\ref{proposition:coboundary_p_norm} yields for $p\in(2,\infty)$ that
\begin{equation}\label{eq:estimate_chi_bis}\SMALL
\bigl|\max_{1\le k\le n}|\chi\circ F_k-\chi|\bigr|_p
\le C\|w\|_\eta n^{1/p}.
\end{equation}

\begin{notation*}
For $n\ge1$ and $g\colon Y^\phi\to \R^N$, 
we write $g_n=\sum_{j=0}^{n-1}g\circ F_j$. 
\end{notation*}

\subsection{Proof of Theorem~\ref{theorem:rate_wass_flow}}
For fixed $w\in\mathcal{F}^\eta(Y^\phi,\R)$ with mean zero and martingale part $m\in L^p(Y^\phi, \R^N)$, $p\in(2,\infty)$, we define the sequence of processes $X_n\colon[0,1]\to \R^N$, $n\ge1$, as
\begin{equation}\label{eq:X_n_flow}
X_{n}(k/n)=\frac{1}{\sqrt{n}}\sum_{j=0}^{k-1}m\circ F_j,
\end{equation}
for $0\le k\le n$
and using linear interpolation in  $[0,1]$.

\begin{proposition}\label{proposition:max_by_Lp}
Let $\{\xi_n\}_{n\ge1}$ be a sequence of identically distributed real random variables, defined on the same probability space. If $\xi_1\in L^q$ for some $q\in[1,\infty)$, then
$|\max_{1\le k\le n}|\xi_k||_q\le n^{1/q}|\xi_1|_q$ for all $n\ge1$.
\begin{proof}
We have that 
$(\max_{1\le k\le n}|\xi_k|)^q=\max_{1\le k\le n}|\xi_k|^q\le \sum_{k=1}^n|\xi_k|^q$. Since all the $\xi_k$ share the same distribution,
$\mathbb{E}[(\max_{1\le k\le n}|\xi_k|)^q]\le
\mathbb{E}[\sum_{k=1}^n|\xi_k|^q]=n\mathbb{E}[|\xi_1|^q]$. The statement follows.
\end{proof}
\end{proposition}

\begin{lemma}\label{lem:estimate_Xn_flow}
There exists $C>0$ such that 
$\mathcal{W}(\widehat{W}_n,X_n)\le Cn^{-\frac{p-2}{2p}}$ for all $n\ge1$.
\begin{proof}
Let $\psi=\int_0^1w\circ F_s\diff s$. By~\eqref{eq:primary_dec}, we have
$\psi_k-m_k=\chi\circ F_k-\chi$ for $k\ge1$. So,
{\small
\begin{align*}
\widehat{W}_n(t)-X_n(t)
&= n^{-1/2}(\psi_{\lfloor nt\rfloor/n}-m_{\lfloor nt\rfloor/n})+R_n(t)= n^{-1/2}(\chi\circ F_{\lfloor nt\rfloor/n}-\chi)+R_n(t)
\end{align*}
}
for all $t\in[0,1]$, where 
$R_n(t)=(\widehat{W}_n(t)-\widehat{W}_n(\lfloor nt\rfloor/n))-(X_n(t)-X_n(\lfloor nt\rfloor/n))$. So,
\begin{equation*}\label{eq:R_n}\SMALL
n^{\frac12}|R_n(t)|\le 
\bigl|\int_{\lfloor nt\rfloor}^{nt}w\circ F_s\diff s \bigr|+
|m\circ F_{\lfloor nt\rfloor-1}|\le
|w|_\infty+
\max_{1\le k\le n}|m\circ F_{k-1}|.
\end{equation*}
By Proposition~\ref{proposition:max_by_Lp} and~\eqref{eq:estimate_m_chi}, 
$$\SMALL n^{-\frac12}\bigl|\max_{1\le k\le n}|m\circ F_{k-1}|\bigr|_p\le
n^{-\frac12+\frac1p}|m|_p\ll n^{-\frac{p-2}{2p}}\|w\|_{\eta}.$$
Hence,
\[\SMALL
\bigl|\sup_{t\in[0,1]}|R_n(t)|\bigr|_p\le 
n^{-\frac12}(|w|_\infty+
\bigl|\max_{1\le k\le n}|m\circ F_{k-1}|\bigr|_p)
\ll  n^{-\frac{p-2}{2p}}\|w\|_{\eta}.
\]
By the estimate on $R_n$ and \eqref{eq:estimate_chi_bis},
\begin{align*}\SMALL
\bigl|\sup_{t\in[0,1]}|\widehat{W}_n(t)-X_n(t)|\bigr|_{p}
\ll
n^{-\frac12}\bigl|\max_{1\le k\le n}|\chi\circ f^{k}-\chi|\bigr|_{p}+n^{-\frac{p-2}{2p}}\ll
n^{-\frac{p-2}{2p}}.
\end{align*}

We finish the proof showing that for any $ f\in \Lip_1$,
\begin{align*}\SMALL
\bigl|\int_{Y^\phi}  f(\widehat{W}_n)\diff\mu^\phi-\int_{Y^\phi}  f(X_n)\diff\mu^\phi\bigr|\le
\bigl|\sup_{t\in[0,1]}|\widehat{W}_n(t)-X_n(t)|\bigr|_{p}
\ll n^{-\frac{p-2}{2p}}.\tag*{\qedhere}
\end{align*}
\end{proof}
\end{lemma}

\begin{proof}[\textbf{Proof of Theorem~\ref{theorem:rate_wass_flow}}]
Let $p\in(2,3)$.
Consider $X_n$ from~\eqref{eq:X_n_flow} and let $Y_n$ be as in Lemma~\ref{lem:W1_piecewise_Brownian_motion}. 
Let $W\colon[0,1]\to\R^N$ be the limiting Brownian motion that has covariance $\Sigma=\int_\mu^\phi mm^T\diff\mu_\Delta$.
By Lemmas~\ref{lem:estimate_Xn_flow} and~\ref{lem:W1_piecewise_Brownian_motion}, to prove the rate on $\mathcal{W}(\widehat{W}_n,W)$ it suffices to estimate $\mathcal{W}(X_n,Y_n)$.

Let us check the hypotheses of Theorem~\ref{theorem:Cuny_et_al} for $d_n=m\circ F_n$, $n\ge0$.
By Proposition~\ref{proposition:primarydec} we have that $L_1m=0$, so Proposition~\ref{proposition:mds_intro} yields that $d_n$  with $\sigma$-algebras $F_1^{-n}\cB$ is an RMDS on the probability space~$(Y^\phi,\cB, \mu^\phi)$. 
It lies in $L^p$ by~\eqref{eq:estimate_m_chi}, and it is  stationary because $F_n$ is measure-preserving. 

Since $m\in\ker L_1$, we follow the proof of Proposition~\ref{proposition:orthogonality} and get 
$\mathbb{E}[(m\circ F_k)(m\circ F_\ell)^T|F_n^{-1}\cB]=0$ for all 
$0\le k\ne \ell\le n-1$. 
Using the notation $\breve{v}=\E[mm^T-\Sigma|F^{-1}_1\cB]$ from~\eqref{eq:vbreve}, we apply~Corollary~\ref{corollary:Burkholder_v_prime} and reason as in the proof of Theorem~\ref{theorem:rate_wass_maps} (in Subsection~\ref{subsec:proof}) to prove that $M_n=\sum_{j=0}^{n-1}m\circ F_j$, $n\ge1$ satisfies condition~\eqref{eq:Cuny_condition}. 

We can now apply Theorem~\ref{theorem:Cuny_et_al} and follow the proof of Theorem~\ref{theorem:rate_wass_maps} to get that $\mathcal{W}(X_n,Y_n)
\ll n^{-\frac{p-2}{2p}}(\log n)^\frac{p-1}{2p}$, which concludes the proof.
\end{proof}

\subsection[]{Proof of Theorem~\ref{theorem:rates_N=1_flow} ($p=\infty$)}\label{subsec:theorem:rates_N=1_flow_infty}
Let $p=\infty$ and $w\in\mathcal{F}^\eta(Y^\phi,\R)$ with mean zero and martingale part $m\in L^\infty(Y^\phi,\R)$.
Define the sequence of processes $Y_n\colon[0,1]\to\R$, $n\ge1$
\begin{equation*}
Y_n(k/n)=\frac{1}{\sqrt{n}}\sum_{j=1}^{k}m\circ F_{n-j},
\end{equation*}
for $1\le k\le n$, using linear interpolation in  $[0,1]$. 
Let $h\colon C([0,1],\R)\to C([0,1],\R)$ be the linear operator $(hf)(t)= f(1)-f(1-t)$ as in~\cite[Lemma 4.8]{KelMel16}.

\begin{lemma}\label{lem:estimate_Wn_Yn}
There exists $C>0$ such that 
$\Pi(h\circ \widehat{W}_n,Y_n)\le C n^{-\frac12}$ for all $n\ge1$.
\begin{proof}
Let $\psi=\int_0^1w\circ F_s\diff s$. 
By equation~\eqref{eq:primary_dec},
\begin{align*}
h\circ \widehat{W}_n(t)- Y_n(t)&=\SMALL
n^{-\frac12}\bigl(\int_{n-\lfloor nt\rfloor}^{n} w\circ F_s\diff s -
 \sum_{j=1}^{\lfloor nt\rfloor-1} m\circ F_{n-j}\bigr)+ R_n(t) \\
&= n^{-\frac12}\bigl(\psi_n-\psi_{n-\lfloor nt\rfloor}
-\bigl(m_n-m_{n-\lfloor nt\rfloor}\bigr)\bigr)+R_n(t)\\
&=n^{-\frac12}\bigl(\chi\circ F_n-\chi\circ F_{n-\lfloor nt\rfloor}\bigr)+R_n(t)
\end{align*}
for every $t\in[0,1]$,  where 
\[\begin{split}\SMALL
R_n(t)&\SMALL=h\circ(\widehat{W}_n(t)-\widehat{W}_n(\lfloor nt\rfloor/n))-
(Y_n(t)-Y_n(\lfloor nt\rfloor/n))\\
&\SMALL=(\widehat{W}_n((1-\lfloor nt\rfloor)/n)-\widehat{W}_n(1-t))-(Y_n(t)-Y_n(\lfloor nt\rfloor/n)).
\end{split}
\]
So, 
$$\SMALL
n^\frac12|R_n(t)|\le
\bigl|\int_{1-t}^{1-\lfloor nt\rfloor}w\circ F_s\diff s\bigr|+
|m\circ F_{n-\lfloor nt\rfloor-1}|\le
|w|_\infty+|m|_\infty,$$ 
and by~\eqref{eq:estimate_m_chi},
$\bigl|\sup_{t\in[0,1]}|R_n(t)|\bigr|_\infty
\ll n^{-\frac12}\|w\|_\eta$. Hence,
\[\SMALL
\bigl|\sup_{t\in[0,1]}|h\circ \widehat{W}_n(t)- Y_n(t)|\bigr|_\infty\ll
n^{-\frac12}( 2|\chi|_\infty+\|w\|_\eta)\ll 
n^{-\frac12}\|w\|_\eta.
\]

Since the Prokhorov metric is bounded by the infinity norm, we conclude that
\[\SMALL
\Pi(h\circ \widehat{W}_n,Y_n)\le\bigl|\sup_{t\in[0,1]}|h\circ \widehat{W}_n(t)- Y_n(t)|\bigr|_\infty\ll n^{-\frac12}.\tag*{\qedhere}
\]
\end{proof}
\end{lemma}

\begin{lemma}\label{lem:Courbot_flow}
There exists $C>0$ such that 
$\Pi(Y_n,W)\le Cn^{-\frac14}(\log n)^{\frac34}$ for all integers $n>1$.
\begin{proof}

Following the proof of Theorem~\ref{theorem:rate_wass_flow}, the sequence $d_n=m\circ F_n$ with $\sigma$-algebras $F_1^{-n}\cB$, $n\ge0$, is a stationary RMDS on the probability space~$(Y^\phi,\cB, \mu^\phi)$.  
Equation~\eqref{eq:estimate_m_chi} yields that $d_n$ is bounded.
We adopt the same notation of Theorem~\ref{theorem:Courbot_applied}, noting that $\sigma^2=\int_{Y^\phi}m^2\diff\mu^\phi$ and that $Y_n$ coincides with $\mathcal{M}_n$. We have that
$$\SMALL V_n(k)=
n^{-1}\sum_{j=1}^k\E[m^2\circ F_{n-j}|F^{-1}_{n-(j-1)}\cB]=
n^{-1}\sum_{j=1}^k\E[m^2|F^{-1}\cB]\circ F_{n-j}.$$

As pointed out in the proof of Lemma~\ref{lem:Courbot}, to apply Theorem~\ref{theorem:Courbot_applied} and complete the current proof, it suffices to show that
$\kappa_n\ll \sqrt{n^{-1}\log n}$.  
Writing 
$\breve{v}=\mathbb{E}[m^2|F^{-1}_1\cB]-\sigma^2$ as in~\eqref{eq:vbreve}, we have that
\begin{equation*}\SMALL
V_n(k)-(k/n)\sigma^2=n^{-1}\sum_{j=1}^k\breve{v}\circ F_{n-j}=
n^{-1}(\breve{v}_n-\breve{v}_{n-k}),
\end{equation*}
for every $n\ge1$. 
So, 
$\max_{0\le k\le n}|V_n(k)-(k/n)\sigma^2|\le
2 n^{-1} \max_{1\le k\le n}|\breve{v}_k|$. 
By Corollary~\ref{corollary:azuma_v_breve}, there are $a,C>0$ such that
\[\SMALL
\mu_\Delta\bigl(\max_{0\le k\le n}|V_n(k)-(k/n)\sigma^2|\ge \epsilon\bigr)\le
\mu_\Delta\bigl(\max_{1\le k\le n}|\breve{v}_k|\ge n\epsilon/2\bigr)
\le Ce^{-an\epsilon^2},
\]
for all $\epsilon\ge0$ and $n\ge1$. Reasoning as in Lemma~\ref{lem:Courbot}, this estimate is sufficient to show that $\kappa_n\ll \sqrt{n^{-1}\log n}$, which concludes the proof.
\end{proof}
\end{lemma}

\begin{proof}[\textbf{Proof of Theorem~\ref{theorem:rates_N=1_flow}~($p=\infty$)}]
Let $w\in\mathcal{F}^\eta(Y^\phi, \R^N)$ with mean zero and let $\widehat{W}_n$ be from~\eqref{eq:Wn'}.
Since $\widehat{W}_n(0)=0$ for all $n\ge1$,  applying Proposition~\ref{proposition:final_inequality} with $N=1$ we get
\[
\Pi(\widehat{W}_n,W)\ll
\Pi(h\circ \widehat{W}_n,W)
\le \Pi(h\circ \widehat{W}_n, Y_n)+ \Pi(Y_n,W).\]
Conclude by Lemmas~\ref{lem:estimate_Wn_Yn} and~\ref{lem:Courbot_flow}.
\end{proof}

\subsection[]{Proof of Theorem~\ref{theorem:rates_N=1_flow} ($p\in(2,\infty)$)}\label{subsec:theorem:rates_N=1_flow_p}
In the current subsection, we use our new estimates from Section~\ref{sec:new_dec} to apply the method found in~\cite[Section 4]{AntMel19} to the semiflow case with real-valued observables. We remark that the following results are proven by the same techniques of~\cite{AntMel19}, and are shown here for completeness.
Let $p\in(2,\infty)$ and 
$w\in\mathcal{F}^\eta(Y^\phi,\R)$ with mean zero and martingale part 
$m\in L^p(Y^\phi,\R)$. 
Consider $\sigma^2=\int_{Y^\phi}m^2\diff\mu^\phi$ and define $d_n=(m\circ F_n)/(n^{1/2}\sigma)$ with $\sigma$-algebras 
$\cG_n=F^{-1}_n\cB$, $n\ge0$. Since $L_1m=0$, the sequence $d_n$ is a RMDS by Proposition~\ref{proposition:mds_intro}. 
Then $(d_{n-j})_{0\le j\le n}$ with filtration $(\cG_{n-j})_{0\le j\le n}$ is a martingale differences array. Define for $0\le k\le n$
\[\SMALL
V_{n}(k)=\sum_{j=1}^k\mathbb{E}[d_{n-j}^2|\cG_{n-(j-1)}].
\]
Define now a sequence of processes $X_n\colon [0,1]\to\R$, $n\ge1$, as
\begin{equation}\label{eq:Xn_Kubilius}
X_n\Biggl(\frac{V_n(k)}{V_n(n)}\Biggr)=\sum_{j=1}^kd_{n-j},
\end{equation}
for $0\le k\le n$, and linear interpolation in $[0,1]$.  
As stated in~\cite{AntMel19}, the integer $k$ in~\eqref{eq:Xn_Kubilius} is a random variable $k=k_{n}(t)\colon Y^\phi\to\{0,\dots,n\}$, such that 
$V_n(k)\le t V_n(n)<V_n(k+1)$.

\begin{proposition}\label{proposition:estimates_nodes}
There exists $C>0$ such that 
$\bigl|\sup_{t\in[0,1]}|k_{n}(t)-\lfloor nt\rfloor|\bigr|_{2(p-1)}\le Cn^\frac12$ for all $n\ge1$.
\begin{proof}
The proof is carried out as in~\cite[Proposition~4.4]{AntMel19}. The only fact left to show is that
\begin{equation}\label{eq:estimate_Kubilius}\SMALL
\bigl|\max_{1\le k\le n}|V_n(k)-k/n|\bigr|_{2(p-1)}\ll n^{-\frac12}.
\end{equation} 
We have that
\[\begin{split}
V_n(k)-\frac{k}{n}=&\frac{1}{n\sigma^2}
\sum_{j=1}^k\mathbb{E}[m^2\circ F_{n-j}|F^{-1}_{n-(j-1)}\cB]-\frac{k}{n}
=\frac{1}{n\sigma^2}
\sum_{j=1}^k\bigl(\mathbb{E}[m^2-\sigma^2|F^{-1}_1\cB]\circ F_{j}\bigr),
\end{split}
\]
and~\eqref{eq:estimate_Kubilius} follows by Corollary~\ref{corollary:Burkholder_v_prime}.
\end{proof}
\end{proposition}

Recall that for $n\ge1$ and $g\colon Y^\phi\to \R^N$, 
we write $g_n=\sum_{j=0}^{n-1}g\circ F_j$.

\begin{proposition}\label{proposition:Z_n}
For $n\ge1$ and $\psi=\int_0^1v\circ F_s\diff s$, define 
$Z_n=\max_{0\le i,\ell\le\sqrt{n}}|\psi_\ell|\circ F_{i\lfloor \sqrt{n}\rfloor}$.
\begin{enumerate}[(a)]
\item\label{item:a} $\bigl|\sum_{j=a}^{b-1}\psi\circ F_j\bigr|\le Z_n((b-a)(\sqrt{n}-1)^{-1}+3)$ 
for all $0\le a< b\le n$.
\item\label{item:b} $|Z_n|_{2(p-1)}\le C\|w\|_\eta n^{1/4+1/(4(p-1))}$ for all $n\ge1$.
\end{enumerate}
\begin{proof}
This proof is identical to the one of~\cite[Proposition~4.6]{AntMel19} with the obvious notational changes and applying equation~\eqref{eq:Rio_applied} to get 
$\bigl|\max_{1\le k\le n}|\psi_k|\bigr|_{2(p-1)}\ll n^{1/2}$.
\end{proof}
\end{proposition}

Let $h\colon C([0,1],\R)\to C([0,1],\R)$ be the linear operator $(hf)(t)= f(1)-f(1-t)$.

\begin{lemma}\label{lem:Xn_Kubilius}
There exists $C>0$ such that 
$\Pi( h\circ \widehat{W}_n,\sigma X_n)\le C n^{-\frac{p-2}{4p}}$ for all $n\ge1$.
\begin{proof}
We follow the proof of~\cite[Lemma 4.7]{AntMel19}. Define the piecewise constant process 
$V'_n(t)=n^{-1/2}\sum_{j=n-\lfloor nt\rfloor}^{n-k-1}\psi\circ F_j$, $t\in[0,1]$, where 
$\psi=\int_0^1w\circ F_s\diff s=m+\chi\circ F_1-\chi$ from equation~\eqref{eq:primary_dec}, and $k=k_n(t)$ is the random variable from~\eqref{eq:Xn_Kubilius}. So,
\begin{equation}\label{eq:estimate_lem_Kub}
\begin{split}\textstyle
h\circ \widehat{W}_n(t)-\sigma X_n(t)&\textstyle=n^{-\frac12}\bigl(
\int_{n-\lfloor nt\rfloor}^nw\circ F_s\diff s-\sum_{j=1}^km\circ F_{n-j}
\bigr)+R_n(t)\\
&=n^{-\frac12}(
\psi_n-\psi_{n-\lfloor nt\rfloor}-(m_n-m_{n-k})
)+R_n(t)\\
&=n^{-\frac12}(
\psi_{n-k}-\psi_{n-\lfloor nt\rfloor}+\chi\circ F_n-\chi\circ F_{n-k}
)+R_n(t)\\
&=V'_n(t) + n^{-\frac12}(\chi\circ F_n-\chi\circ F_{n-k}
)+R_n(t),
\end{split}
\end{equation}
for every $t\in[0,1]$, where
$\bigl|\sup_{t\in[0,1]}|R_n(t)|\bigr|_p\le n^{-\frac12}\bigl(|w|_\infty+\bigl|\max_{1\le k\le n}|m\circ F_{k-1}|\bigr|_p\bigr)$. 
Reasoning as in the proof of Lemma~\ref{lem:estimate_Xn_flow}, we get 
$\bigl|\sup_{t\in[0,1]}|R_n(t)|\bigr|_p\ll  n^{-\frac{p-2}{2p}}\|w\|_\eta$. Using~\eqref{eq:estimate_chi_bis}, 
\[\begin{split}\SMALL
n^{-\frac12}\bigl|\sup_{t\in[0,1]}|\chi\circ F_n-
\chi\circ F_{n-k_n(t)}|\bigr|_p&\SMALL=
n^{-\frac12}\bigl|\sup_{t\in[0,1]}|\chi\circ F_{k_n(t)}-
\chi|\bigr|_p\\
&\SMALL=n^{-\frac12}\bigl|\max_{1\le k\le n}|\chi\circ F_k-
\chi|\bigr|_p\ll n^{-\frac{p-2}{2p}}.
\end{split}\]
By Propositions~\ref{proposition:estimates_nodes} and~\ref{proposition:Z_n}, and by Cauchy-Schwarz,
\begin{align*}\SMALL
\bigl|\sup_{t\in[0,1]}|V'(t)|\bigr|_{p-1}&\SMALL\le
n^{-\frac12}\bigl|Z_n(n^{-\frac12}\sup_{t\in[0,1]}|
\lfloor nt\rfloor -k_n(t)|+3)\bigr|_{p-1}\\
&\SMALL\le n^{-\frac12} |Z_n|_{2(p-1)}(
n^{-\frac12}\bigl|\sup_{t\in[0,1]}|
\lfloor nt\rfloor -k_n(t)|\bigr|_{2(p-1)}+3)\\
&\SMALL\ll n^{-\frac12} |Z_n|_{2(p-1)}\ll
n^{-(\frac14-\frac{1}{4(p-1)})}=n^{-\frac14\frac{p-2}{p-1}}.
\end{align*}
Applying these estimates to~\eqref{eq:estimate_lem_Kub}, 
$|\sup_{t\in[0,1]}|h\circ \widehat{W}_n(t)-\sigma X_n(t)||_{p-1}\ll
n^{-\frac14\frac{p-2}{p-1}}$. Finish by applying~\cite[Proposition 4.5(b)]{AntMel19} with $q=p-1$.
\end{proof}
\end{lemma}

\begin{lemma}\label{lem:Kubilius_applied}
Let $B\colon[0,1]\to\R$ be a standard Brownian motion. 
There exists $C>0$ such that 
$\Pi(X_n,B)\le C n^{-\frac{p-2}{4p}}$ for all $n\ge1$.
\begin{proof}
This is identical to the proof of of~\cite[Lemma 4.3]{AntMel19}, adapting the notation and applying~\cite[Theorem 1]{Kub94} of Kubilius. We remark the importance of~\eqref{eq:estimate_Kubilius} to finish this proof. (See also~\cite[Theorem 3.59]{Pav23}).
\end{proof}
\end{lemma}

\begin{proof}[\textbf{Proof of Theorem~\ref{theorem:rates_N=1_flow} ($p\in (2,\infty)$)}]
This proof is identical to the proof of~\cite[Theorem~2.2]{AntMel19} with the appropriate notational changes. We write it to keep this work self-contained. 
Let $\widehat{W}_n$ be from~\eqref{eq:Wn'}.
Since $\widehat{W}_n(0)=0$, Proposition~\ref{proposition:final_inequality} yields that
$\Pi(\widehat{W}_n,W)\le 2
\Pi(h\circ \widehat{W}_n,W).$ Using that $W=_d\sigma B$, we get
\[
\Pi(\widehat{W}_n,W)\ll\Pi(h\circ \widehat{W}_n,\sigma B)
\ll \Pi(h\circ \widehat{W}_n, \sigma X_n)+ \Pi(\sigma X_n,\sigma B).
\]
Conclude by Lemmas~\ref{lem:Xn_Kubilius} and~\ref{lem:Kubilius_applied}.
\end{proof}

\section{Nonuniformly hyperbolic flows}\label{sec:sinai_flow}
In the previous sections we showed how to derive rates in the WIP for certain non-invertible systems using martingale techniques and Gordin's method. However, these approaches are not applicable directly in the invertible setting. Therefore, in this section we illustrate how the theorems outlined in Section~\ref{sec:setup_and_main_results} remain valid for a class of invertible systems.

Following the arguments of~\cite[Remark 6.2(a)]{AntMel19}, Theorems~\ref{theorem:rate_wass_maps} and~\ref{theorem:rates_N=1_maps} are still valid for a class of maps that are nonuniformly hyperbolic in the sense of Young~\cite{You98,You99}, namely they display polynomial tails with inducing time in $L^p$, and have an exponential contraction along stable leaves. Such conditions are satisfied if the Young tower displays exponential tails.
However, the matter of passing the rates from semiflows to flows is more delicate than in the discrete setting.  In the following, we show how Theorems~\ref{theorem:rate_wass_flow} and~\ref{theorem:rates_N=1_flow} are satisfied for a class of nonuniformly hyperbolic flows that display an exponential contraction along stable leaves.

\subsection{Setup}\label{subsec:setup}

Let $(M,d)$ be a bounded metric space  and let $\Psi_t\colon M\to M$, $t\in\R$, be a flow, so $\Psi_0=\Id$ and $\Psi_{t+s}=\Psi_t\circ\Psi_s$, $t,s\in\R$. As in Subsection~\ref{subsec:nonunif_exp_semiflow}, we assume Lipschitz continuous dependence on initial conditions, on initial conditions~\eqref{eq:continuous_dependence} and  Lipschitz continuity in time~\eqref{eq:Lipschitz}. 
Our main assumption is that $\Psi_t$ can be modelled by a suspension with Hölder continuous first return time, where the induced map is nonuniformly hyperbolic in the sense of Young~\cite{You98}.
Here follows a detailed list of our assumptions.

\begin{itemize}

\item  There exists $X\subset M$ and a function $r\colon X\to[1,\infty)$, $r(x)=\inf\{t>0:\Psi_tx\in X\}$, such that $T\colon X\to X$, $Tx=\Psi_{r(x)}x$ is nonuniformly hyperbolic. We suppose that there is $\eta_0\in(0,1]$ such that $|r|_{\eta_0}=\sup_{x\ne x'}|r(x)-r(x')|/d(x,x')^{\eta_0}<\infty$.

\item There is a Borel probability measure $m$ on $X$, which is preserved by $T$.

\item There exists a measurable $Y\subset X$, $m(Y)>0$, with an (at most countable) measurable partition $\{Y_j\}_{j\ge1}$ and a function $\tau\colon Y\to\Z_+$ that is constant on partition elements, such that $T^{\tau(y)}y\in Y$ for all $y\in Y$. We define $F\colon Y\to Y$ as 
$F y =  T^{\tau(y)}y$.

\item Let $p\in[2,\infty]$. We have an $F$-invariant probability measure $\mu$ on $Y$ such that  
$\tau\in L^p(Y,\mu)$. 

\item There is a partition of $Y$ consisting of stable leaves $\mathcal{W}^s$ of $F$, which is a refinement of~$\{Y_j\}_{j\ge1}$. For $y\in Y$, let $\mathcal{W}^s(y)$ denote the stable leaf containing $y$. The stable leaves are invariant under $F$, that is $F(\mathcal{W}^s(y))\subset \mathcal{W}^s(Fy)$. 

\item We have the quotient space $\overline{Y}=Y/\sim$, where $y\sim y'$ if $y\in \mathcal{W}^s(y')$, with projection $\overline{\pi}\colon Y\to \overline{Y}$. We can also quotient the map $F$ into $\overline{F}\colon \overline{Y}\to \overline{Y}$ with invariant probability measure $\overline{\mu}=\overline{\pi}_*\mu$. 

\item There is a partition $\{\overline{Y}_j\}_{j\ge1}$ of $\overline{Y}$ such that $Y_j=\overline{\pi}^{-1}(\overline{Y}_j)$. The separation time $s(y,y')$ is defined as the infimum of $n\in\N_0$ such that $\overline{F}^ny$ and $\overline{F}^ny'$ belong to different partition elements.
The function $\overline{F}$ is a full-branch Gibbs-Markov map and it separates trajectories:  $s(y,y')<\infty$ if and only if $y\ne y'$. We extend the separation time to $y,y'\in Y$ by $s(y,y')=s(\overline{\pi}(y),\overline{\pi}(y'))$. 

\item Let $\wt{Y}\subset Y$ be such that it intersects every stable leaf only once, that is for every $y\in Y$ the set $\mathcal{W}^s(y)\cap\wt{Y}$ contains exactly one element, and for $j\ge1$ let $\wt{Y}_j=\wt{Y}\cap Y_j$. Let $\pi\colon Y\to \wt{Y}$ be the associated projection, so that $\pi y\in\mathcal{W}^s(y)\cap \wt{Y}$ and  $\pi y=\pi y'$ if and only if $y'\in\mathcal{W}^s(y)$.
\end{itemize}

For $n\ge0$ and $y\in Y$, we let $\beta_n(y)=N$ to be the unique integer such that 
\begin{equation}\label{eq:lap_numer_f}
\SMALL\sum_{j=0}^{N-1}\tau\circ F^j(y)\le n<
\sum_{j=0}^{N}\tau\circ F^j(y).    
\end{equation}
To conclude our list of assumptions on $T$, we include some estimates.

\begin{itemize}
\item Suppose there exist~$C>0$ and $\gamma\in(0,1)$ such that
\begin{equation}\label{eq:first_assumption}\SMALL
d(T^ny,T^ny')\le C\bigl(\gamma^nd(y,y')+\gamma^{s(y,y')-\beta_n(y)}\bigr),
\end{equation}
for any $n\ge0$ and $y,y'\in Y$.

\item Finally, we assume that
\begin{equation}\label{eq:second_assumption}
d(T^ny,T^ny')\le C\gamma^{s(y,y')-\beta_n(y)},
\end{equation}
for all $n\ge0$ and $y,y'\in \wt Y$.
\end{itemize}
In particular, by~\eqref{eq:first_assumption} we have contraction of $T$ along stable leaves:
\begin{equation}\label{eq:exp_contr}
d(T^ny,T^ny')\le C\gamma^nd(y,y'),
\end{equation}
for all $n\ge0$ and $y,y'\in Y$ with $y'\in\mathcal{W}^s(y)$.

\begin{definition}
A flow $\Psi_t \colon M \to M$, $t \in \mathbb{R}$, satisfying all the assumptions above is called a \textit{nonuniformly hyperbolic flow}.    
\end{definition}

\begin{remark}
Condition~\eqref{eq:exp_contr} is essential for our purposes and was assumed by~\cite{Mel07,Mel18}. It is satisfied by flows modelled by Young towers with exponential tails and by some classes of intermittent solenoids. Nevertheless, we mention there are slow mixing models of interests where~\eqref{eq:exp_contr} is not satisfied, such as Bunimovich flowers~\cite{Bun73}.
\end{remark}

Define  $\phi\colon Y\to[1,\infty)$ as 
$\phi(y)=\sum_{j=0}^{\tau(y)-1}r(T^jy)$. 
By $\phi\le|r|_\infty\tau$ and $\tau\in L^p(Y)$, we get that $\phi\in L^p(Y)$.
Define 
$Y^\phi=\{(y,u)\in Y\times[0,\infty):u\in [0,\phi(y)]\}/\sim$
where $(y,\phi(y))\sim (Fy,0)$.
The suspension flow $F_t:Y^\phi\to Y^\phi$ is given by
$F_t(y,u)=(y,u+t)$ computed modulo identifications. 
The projection 
$\pi_M\colon Y^\phi\to M$, $\pi_M(y,u)=\Psi_uy$, is a semiconjugacy from $F_t$ to $\Psi_t$. 
We define the ergodic $F_t$-invariant probability measure
$\mu^\phi=(\mu\times{\rm Lebesgue})/\bar\phi$,
where $\bar\phi=\int_Y\phi\,d\mu$. Then, $\nu=(\pi_M)_*\mu^\phi$ is an ergodic $\Psi_t$-invariant probability measure on $M$.\vspace{1ex}

We define the space of Hölder functions $C^\eta(M,\R^N)$ with norm~$\|\cdot\|_\eta$ similarly to Definition~\ref{def:Hölder}.
Let $v\in C^\eta(M,\R^N)$ with mean zero and define the sequence 
\[
W_n(t)=\frac{1}{\sqrt{n}}\int_0^{nt}v\circ\Psi_s\diff s,
\]
for $n\ge1$ and $t\in[0,1]$. Every $W_n$ is a random element in $C([0,1],\R^N)$ defined on the probability space $(M,\nu)$.
If $\Psi_t$ is a nonuniformly hyperbolic flow and $p\in[2,\infty]$, then the WIP is satisfied for Hölder observables. This follows by passing the WIP for maps~\cite{Melvar16} to the flow using~\cite{KelMel16}. So, there exists a centred Brownian motion $W\colon[0,1]\to\R^N$ such that $W_n\to_w W$. We state now the last theorem of this paper.

\begin{theorem}\label{theorem:rates_NUH_flows}
Let $v\in C^\eta(M,\R^N)$ with mean zero and let $p\in(2,\infty]$. The rates of convergence for the WIP in Wasserstein and Prokhorov metrics of Theorems~\ref{theorem:rate_wass_flow} and~\ref{theorem:rates_N=1_flow} are also valid for nonuniformly hyperbolic flows, with the same conditions on $p$ and $N\ge1$.
\end{theorem}

The remainder of this section is dedicated to proving Theorem~\ref{theorem:rates_NUH_flows}. Let us outline our approach. Firstly, we show how to model a nonuniformly hyperbolic flow by a suspension with a roof function that is constant along stable leaves, and which projects into a Gibbs-Markov semiflow. 
Secondly, we establish that a lift of a Hölder observable $v$ is cohomologous to a function $h$ that only depends on future coordinates. This function projects to an observable $\overline{h}$  that fits within the appropriate functional space described Section~\ref{sec:new_dec}.
We obtain the rates for the WIP for $\overline{h}$ following the approach of Section~\ref{sec:continuous_time_rates}. Finally, these rates are passed to the WIP for the function $v$ through standard arguments.

We start by presenting a couple of estimates. We fix $\gamma_1=\gamma^{\eta_0}\in(0,1)$, where~$\eta_0\in(0,1]$ is the Hölder exponent for $r\colon X\to [1,\infty)$ and $\gamma\in(0,1)$ is from equations~\eqref{eq:first_assumption} and~\eqref{eq:second_assumption}.
For the return time function $r\colon X\to[1,\infty)$ defined above and $k\ge1$, henceforth we write $r_k=\sum_{j=0}^{k-1}r\circ T^j$.

\begin{proposition}\label{proposition:phi_NUH_flows}
There is $C>0$ such that, for all $j\ge1$, $y,y'\in \wt{Y}_j$, and $0\le n\le \tau(y)$,
\[\SMALL
|r_n(y)-r_n(y')|\le C(\inf_{Y_j}\phi)\gamma_1^{s(y,y')}.
\]
In particular,
\[\SMALL
|\phi(y)-\phi(y')|\le C(\inf_{Y_j}\phi)\gamma_1^{s(y,y')}.
\]
\begin{proof}
    Let us prove the estimate as done in~\cite[Proposition~7.4]{BalButMel19}. 
For $y,y'\in\wt Y_j$ and $0\le j\le\tau(y)-1$, note that $\beta_j(y)=0$. Hence, equation~\eqref{eq:second_assumption} gives a $C>0$ such that
\begin{equation}\label{eq:derived_estimates}\SMALL
d(T^jy,T^jy')\le C\gamma^{s(y,y')}.
\end{equation}
So,
\[\SMALL
|r_n(y)-r_n(y')|\le\sum_{j=0}^{n-1}|r(T^jy)-r(T^jy')|\le
|r|_{\eta_0} \sum_{j=0}^{n-1    }d(T^jy,T^jy')^{\eta_0}\ll n\gamma_1^{s(y,y')}.
\]
The first estimate follows because $n\le\tau(y)\le\inf_{Y_j}\phi$, whereas the inequality for $\phi$ follows by taking $n=\tau(y)$.
\end{proof}
\end{proposition}

\begin{proposition}\label{proposition:regularity}
For every $v\in C^\eta(M,\R^N)$ there exists $C>0$ such that,  
\[\SMALL
|v(\Psi_uy)-v(\Psi_sy')|\le C\bigl((\inf_{Y_j}\phi)(\gamma_1^{\eta})^{s(y,y')}+|u-s|^\eta\bigr),
\]
for all $(y,u),(y',s)\in Y^\phi$ with $y,y'\in \wt{Y}_j$.
\begin{proof} 
We have that
\[
|v(\Psi_uy)-v(\Psi_sy')|\le |v|_\eta d(\Psi_uy,\Psi_sy')^\eta\le |v|_\eta\{
d(\Psi_uy,\Psi_uy')^\eta+d(\Psi_uy',\Psi_sy')^\eta\}.
\]
We use Lipschitz continuity~\eqref{eq:Lipschitz} for the last term:
$d(\Psi_uy',\Psi_sy')^\eta\le L^\eta|u-s|^\eta$.

Let $n=n(y,u)\ge0$ be the unique integer such that $r_n(y)\le u<r_{n+1}(y)$.
Reasoning as in Proposition~\ref{proposition:v_in_F},
\[
d(\Psi_uy,\Psi_uy')\ll
d(T^ny,T^ny')+ L|r_n(y)-r_n(y')|.
\]
By $u\le\phi(y)$, we have $n\le\tau(y)$ and we can apply equation~\eqref{eq:derived_estimates}. Combining this with Proposition~\ref{proposition:phi_NUH_flows},  we conclude that
\[\SMALL
d(\Psi_uy',\Psi_uy')^\eta\ll \bigl(\gamma_1^{s(y,y')}\bigr)^\eta+
\bigl((\inf_{Y_j}\phi)\gamma_1^{s(y,y')}\bigr)^\eta\le(\inf_{Y_j}\phi)(\gamma_1^\eta)^{s(y,y')}.\tag*{\qedhere}
\]
\end{proof}
    \end{proposition}

\begin{remark}
Our assumptions are naturally satisfied when the flow $\Psi_t$ is uniformly hyperbolic. In such a case, we have $Y=X$, the return function $\tau$ is constantly $1$, and $\phi=r\in L^\infty(Y)$ is Hölder continuous. Moreover, an observable $v\in C^\eta(M)$ lifts directly to a Lipschitz observable with respect to a two-sided symbolic metric.
\end{remark}

\subsection{Reduction to a roof function constant along stable leaves}\label{subsec:reduction}

We recall here standard arguments from~\cite{Bow75,Sin72} in order to quotient the suspension $Y^\phi$ along stable leaves. For a more complete and general exposition of our setting, we refer to~\cite{BalButMel19}, where the authors do not assume the exponential contraction~\eqref{eq:exp_contr}. We remark that our methods rely strongly on such a requirement. 

We introduce the Young tower~\cite{You99} 
$\Delta=\{(y,\ell)\in Y\times\Z:0\le\ell\le\tau(y)-1\}$ with the tower map $f\colon\Delta\to\Delta$ 
\[
f(y,\ell)=
\begin{cases}
(y,\ell+1),&\ell\le\tau(y)-2\\
(Fy,0),&\ell=\tau(y)-1
\end{cases}.
\]
Let $\pi_X\colon\Delta\to X$ be the projection $\pi_X(y,\ell)=T^\ell y$, hence $\pi_X\circ f=T\circ\pi_X$. Using that $\pi y$ and $y$ belong to the same partition element for any $y\in Y$, and that $\tau$ is constant along partition elements, we can define $\pi\colon \Delta\to\Delta$ as $\pi(y,\ell)=(\pi y,\ell)$.

We define $\widetilde{\chi}\colon\Delta\to\R$ as
$\widetilde{\chi}=
\sum_{n=0}^\infty \bigl(r\circ\pi_X\circ f^n\circ\pi-r\circ\pi_X\circ f^n\bigr)$. 
Using that $\pi_X\circ f^n=T^n\circ\pi_X$ for any $n\ge0$, we get
\begin{align*}\SMALL
|\wt\chi(y,\ell)|&\SMALL\le\sum_{n=0}^\infty
\bigl|r\circ\pi_X(f^n(\pi y,\ell))-r\circ\pi_X(f^n(y,\ell))\bigr|\\
&\SMALL=\sum_{n=0}^\infty
\bigl|r(T^{n+\ell}\pi y)-r(T^{n+\ell} y)\bigr|\\
&\le\SMALL|r|_{\eta_0}\sum_{n=0}^\infty  d\bigl(T^{n+\ell}\pi y,T^{n+\ell} y\bigr)^{\eta_0}.
\end{align*}
By~\eqref{eq:exp_contr}, the series for $\widetilde{\chi}$ converges absolutely on $\Delta$, and the function $\widetilde{\chi}$ is bounded. 
For~$y\in Y$, we write $\wt \chi_Y(y)=\wt\chi(y,0)$, so naturally $|\wt\chi_Y|_\infty\le|\wt\chi|_\infty$. Note that 
\begin{equation}\label{eq:chi_Y}\SMALL
\wt\chi_Y(y)=\sum_{n=0}^\infty\bigl(r(T^n\pi y)-r(T^n y)\bigr),
\end{equation}
and $\chi_Y(y)=0$  for any $y\in\wt Y$.

For the remaining of this subsection, we fix $\gamma_1=\gamma^{\eta_0}$ and $\theta=\gamma_1^{1/2}$, where~$\eta_0\in(0,1]$ is the Hölder exponent for $r\colon X\to [1,\infty)$ and $\gamma\in(0,1)$ is from equations~\eqref{eq:first_assumption} and~\eqref{eq:second_assumption}.

\begin{proposition}\label{proposition:chi}  
There exists $C>0$ such that
\[\SMALL
|\wt\chi_Y(y)-\wt\chi_Y(y')|\le C\bigl(d(y,y')^{\eta_0}+\theta^{s(y,y')}\bigr),
\]
for all $y,y'\in Y$.
\begin{proof}
This proof follows~\cite[Lemma 8.4]{BalButMel19}. 
Suppose $y,y'\in Y$ and write $N=\lfloor s(y,y')/2)\rfloor$.
By~\eqref{eq:chi_Y},
\[\SMALL
|\wt\chi_Y(y)-\wt\chi_Y(y')|\le A(y,y')+A(\pi y,\pi y')+B(y)+B(y'),
\]
where 
\[\begin{split}
&\SMALL A(y,y')=\sum_{n=0}^{N-1}|r(T^ny)-r(T^ny')|,\\
&\SMALL B(y)=\sum_{n=N}^{\infty}|r(T^ny)-r(T^n(\pi y))|.
\end{split}
\]

By equation~\eqref{eq:exp_contr}, for all $y\in Y$
\[\SMALL
B(y)\le |r|_{\eta_0}\sum_{n=N}^\infty d(T^ny,T^n\pi y)^{\eta_0}
\ll\sum_{n=N}^\infty (\gamma^{\eta_0})^n.
\]
By $\gamma_1=\gamma^{\eta_0}$ and
$
\sum_{n=N}^\infty \gamma_1^n
= \gamma_1^N/(1-\gamma_1),
$
it follows that $B(y),B(y')\ll \gamma_1^N\ll\theta^{s(y,y')}$.
\vspace{1ex}

Let us focus on $A$. Using~\eqref{eq:first_assumption}, we get
\begin{align*}\SMALL
A(y,y')
& \le |r|_{\eta_0}\sum_{n=0}^{N-1} d(T^ny,T^ny')^{\eta_0}
\le C|r|_{\eta_0}\sum_{n=0}^{N-1} (\gamma_1^n\,d(y,y')^{\eta_0} +\gamma_1^{s(y,y')-\beta_n(y)})
\\
& \le C|r|_{\eta_0}\sum_{n=0}^{N-1} (\gamma_1^n\,d(y,y')^{\eta_0} +\gamma_1^{s(y,y')-n})
 \\ & 
\ll d(y,y')^{\eta_0}+\gamma_1^{s(y,y')-N}\ll d(y,y')^{\eta_0}+ \theta^{s(y,y')}.
\end{align*}
Similarly, using~\eqref{eq:second_assumption},
\[
A(\pi y,\pi y') \ll  \theta^{s(\pi y,\pi y')} = \theta^{s(y,y')}.
\]
The statement is proved by combining the estimates for $A$ and~$B$.
\end{proof}
\end{proposition}

We recall the notation $r_k=\sum_{j=0}^{k-1}r\circ T^j$ for $k\ge1$.
Since $\inf r\ge1$, there exists~$k\ge1$ such that $\inf r_{k}\ge 4|\widetilde{\chi}|_\infty+1$.
We take without loss $k=1$ (otherwise in the following we could substitute $r$ with $r_k$ and $T$ with~$T^k$). This implies that $\inf \phi\ge 4|\widetilde{\chi}|_\infty+1$, for $\phi\in L^p(Y)$ the roof function defined in Subsection~\ref{subsec:setup}.
For $w=r\circ\pi_X$, we define $\wt{r}\colon \Delta\to\R$ as
\begin{equation}\label{eq:tilde_r}
\wt{r}=w+\wt{\chi}\circ f-\wt{\chi},    
\end{equation}
so that $\inf\wt{r}\ge\inf r-2|\wt{\chi}|_\infty\ge1$.
A calculation gives
\begin{equation}\label{eq:tilde_r_consant}
\begin{split}
\wt{\chi}\circ f-\wt{\chi}&=\sum_{n=0}^\infty 
\bigl(w\circ f
^{n+1}-w\circ f^n\circ\pi\circ f-
w\circ f^n+w\circ f^n\circ\pi\bigr)\\ 
&=\sum_{n=0}^\infty 
\bigl(w\circ f^{n+1}\circ\pi-w\circ f^n\circ\pi\circ f\bigr)-w+w\circ\pi=H-w+w\circ\pi.
\end{split}
\end{equation}
Since $\pi\circ f=\pi\circ f\circ\pi$, we have that $\wt r=w\circ\pi+H$ is  constant along stable leaves, that in this context means $\wt r\circ\pi=\wt r$.
\vspace{1ex}

We let $\wt{\phi}\colon Y\to\R$
be
$\wt{\phi}(y)=\sum_{\ell=0}^{\tau(y)-1}\wt r(y,\ell)$, so $\inf\wt \phi\ge1$.
So, 
$\wt{\phi}(y)=\wt{\phi}(\pi y)$ for all $y\in Y$, and so
$\wt{\phi}$ is constant along stable leaves.
A calculation gives that
\begin{equation}\label{eq:tilde_phi_coboudary}
\wt \phi=\phi+\wt\chi_Y\circ F-\wt\chi_Y.    
\end{equation}
It follows that
$\int_Y\widetilde{\phi}\diff \mu=\int_Y \phi\diff \mu$, and so $\wt{\phi}$ is an integrable roof function. Moreover, by $\phi\in L^p(Y)$ and
$|\wt{\phi}|_p\le|\phi|_p+2|\wt{\chi}_Y|_\infty$, we have that
$\wt{\phi}\in L^p(Y)$. 

\begin{proposition}\label{proposition:phi_final}
There  exists $C>0$ such that 
\[\SMALL
|\wt{\phi}(y)-\wt{\phi}(y')|\le C(\inf_{Y_j}\wt{\phi})\theta^{s(y,y')}.
\]
for any $j\ge1$ and $y,y'\in Y_j$.
\begin{proof}

We follow the proof of~\cite[Proposition~6.1]{BalButMel19}. Since 
$\inf \phi\ge 4|\widetilde{\chi}|_\infty+1$, we have that 
\begin{equation*}\label{eq:phi_wtphi}\SMALL
\inf_{Y_j}\phi\le \inf_{Y_j}\wt{\phi}+2|\wt{\chi}|_\infty\le 
\inf_{Y_j}\wt{\phi}+\frac12\inf\phi\le
\inf_{Y_j}\wt{\phi}+\frac12\inf_{Y_j}\phi,
\end{equation*}
giving that $\inf_{Y_j}\phi\le2\inf_{Y_j}\wt{\phi}$.

Since $\wt{\phi}$ is constant along stable leaves, we can assume without loss that $y,y'\in \wt{Y}_j$. 
Using that $\chi_Y=0$ on $\wt Y$, we get from~\eqref{eq:tilde_phi_coboudary}
\[\SMALL
|\wt{\phi}(y)-\wt{\phi}(y')|\le|\phi(y)-\phi(y')|+
|\wt\chi_Y(Fy)-\wt{\chi}_Y(Fy')|.
\] 
By Proposition~\ref{proposition:phi_NUH_flows},
\[\SMALL
|\phi(y)-\phi(y')|\ll(\inf_{Y_j}\phi)\theta^{s(y,y')}\le 2
(\inf_{Y_j}\wt{\phi})\theta^{s(y,y')}.
\]
To bound the term for $\wt\chi_Y$, we note that $\tau(y)=\tau(y')$ and $\beta_{\tau(y)}(y)=1$. So,~\eqref{eq:second_assumption} yields that
\[
d(Fy,Fy')=d(T^{\tau(y)}y,T^{\tau(y)}y')\ll \gamma^{s(y,y')-1}\ll\gamma^{s(y,y')}.
\]
By $s(Fy,Fy')=s(y,y')-1$, Proposition~\ref{proposition:chi} implies that
\[\SMALL
|\wt\chi_Y(Fy)-\wt{\chi}_Y(Fy')|\ll d(Fy,Fy')^{\eta_0}+\theta^{s(Fy,Fy')}
\ll \gamma_1^{s(y,y')}+\theta^{s(y,y')}\ll\theta^{s(y,y')}.    
\]
Combine the previous estimates to conclude the proof.
\end{proof}
\end{proposition}

Let $Y^{\wt{\phi}}$ be the suspension over the map $F\colon Y\to Y$ with roof function $\wt{\phi}$, and write $\wt{F}_t\colon Y^{\wt{\phi}}\to Y^{\wt{\phi}}$ for the associated suspension flow. We define the probability measure 
$\mu^{\wt{\phi}}=(\mu\times{\rm Lebesgue})/\int\wt\phi\diff\mu$.
Let $g\colon Y^{\wt{\phi}}\to Y^\phi$ be  
$g(y,u)=(y,u+\wt\chi_Y(y))$, following the identifications on $Y^   \phi$.
So,     
\[\SMALL
g(y,\wt{\phi}(y))=(y,\wt{\phi}(y)+\wt\chi_Y(y))=(y,\phi(y)+\wt\chi_Y(Fy))
=(Fy,\wt\chi_Y(Fy))=g(Fy,0).
\]
Hence, the function $g$ respects the identifications on $Y^{\wt{\phi}}$, and so it is well-defined. Note that $F_t\circ g=g\circ \wt{F}_t$, where $F_t$ is the suspension flow on $Y^\phi$. Let $\wt{\pi}_M=\pi_M\circ g$; since both $\pi_M$ and $g$ are measure-preserving, so it is $\wt \pi_M$. By $\Psi_t\circ \pi_M=\pi_M\circ F_t$, we get that 
$\Psi_t\circ  \wt{\pi}_M=\wt{\pi}_M\circ\wt{F}_t$, and so $\wt{\pi}_M$ is a semiconjugacy between $(Y^{\wt \phi},\wt{F}_t,\mu^{\wt\phi})$ and $(M,\Psi_t,\nu)$. 

\begin{lemma}\label{lem:G_n}
There is $C>0$ such that 
\[
\bigl|\wt{\chi}(f^n(y,0))-\wt{\chi}(f^n(y',0))\bigr|\le C\gamma_1^n,
\]
for all $n\ge1$ and $y,y'\in Y$ such that $y'\in\mathcal W^s(y)$.
\end{lemma}

\begin{proof}
For $n\ge1$ and $y\in Y$, let $N=\beta_n(y)$ satisfy~\eqref{eq:lap_numer_f}. Writing $\tau_k=\sum_{j=0}^{k-1}\tau\circ F^j$ for $k\ge1$, we can describe the iterations of the tower map as
\[
f^n(y,0)=\bigl(F^Ny,n-\tau_N(y)\bigr),
\]
where $0\le n-\tau_N(y)\le \tau(F^Ny)-1$. 
If $y$ and $y'$ belong to the same stable leaf, then $\tau(y)=\tau(y')$ and $\pi F^Ny=\pi F^Ny'$.
Hence, 
\begin{equation}\label{eq:proj_f^n}
\pi(f^n(y,0))=\bigl(\pi F^Ny,n-\tau_N(y)\bigr)=\bigl(\pi F^Ny',n-\tau_N(y')\bigr)=\pi(f^n(y',0)).
\end{equation}

Using the definition of $\wt \chi$ and $\pi_X\circ f=T\circ \pi_X$, we can write
\begin{align*}
\wt \chi(f^n(y,0))&=
\sum_{j=0}^\infty\Bigl(r\circ\pi_X\bigl(f^j(\pi(f^n(y,0)))\bigr)-
r\circ\pi_X\bigl(f^j(f^n(y,0))\bigr)\Bigr)\\
&=\sum_{j=0}^\infty\Bigl(r\circ\pi_X\bigl(f^j(\pi(f^n(y,0)))\bigr)-
r(T^{j+n}y)\bigr)\Bigr).
\end{align*} 
So, by the identity~\eqref{eq:proj_f^n},
\[
G_n=\wt \chi(f^n(y,0))-\wt \chi(f^n(y',0))=
\sum_{j=0}^\infty\bigl(r(T^{j+n}y')-r(T^{j+n}y)\bigr),
\]
for any $y,y'\in Y$ such that $y'\in\mathcal W^s(y)$.
By~\eqref{eq:exp_contr},
\begin{align*}
|G_n|&\le
|r|_{\eta_0}\sum_{j=0}^\infty 
d\bigl(T^{j+n}y',T^{j+n}y\bigr)^{\eta_0}
\ll\sum_{j=0}^{\infty}\gamma_1^{j+n}
=\gamma_1^n/(1-\gamma_1)\ll\gamma_1^n.
\end{align*}
This finishes the proof.
\end{proof}

We can now show that the exponential contraction along stable leaves of the map $F$ can be lifted to the suspension flow $\wt F_t\colon Y^{\wt \phi}\to  Y^{\wt \phi}$.
\begin{proposition}\label{proposition:exp_decay_flow}
There exist $C>0$ and $\gamma_2\in(0,1)$ for which
\[
d\bigl(\wt{\pi}_M\circ\wt{F}_t(y,0),\wt{\pi}_M\circ\wt{F}_t(y',0)\bigr)
\le C\gamma_2^t,
\]
for all $t\ge0$ and $y,y'\in Y$ such that $y'\in\mathcal{W}^s(y)$.
\begin{proof}
For $k\ge1$, and $\wt r\colon\Delta\to[1,\infty)$ from~\eqref{eq:tilde_r}, we write 
$S_k\wt r=\sum_{j=0}^{k-1}\wt r\circ f^j$.
For $y\in Y$ and $t\ge0$, we let $n=n(y,t)\ge0$ such that $S_n\wt{r}(y,0)\le t<S_{n+1}\wt{r}(y,0)$. Since $\wt{r}$ is constant along stable leaves, it follows that  
$t=S_n\wt{r}(y,0)+E(y)=S_n\wt{r}(y',0)+E(y)$, where $E(y)\le|\wt{r}|_\infty$. Note that
\begin{align*}
\wt\pi_M \circ \tilde F_t (y,0) &= \pi_M \circ g \circ \wt F_t(y,0) = \pi_M \circ F_t\circ g(y,0) \\
&= \pi_M \circ F_t \circ (y, \wt\chi_Y(y)) =\pi_M\circ F_{t+\wt\chi_Y(y)}(y,0)\\
&=\Psi_{t+\wt\chi_Y(y)}(y).   
\end{align*}
So, by applying~\eqref{eq:continuous_dependence}, 
\[\begin{split}
d(\wt{\pi}_M\circ\wt{F}_t(y,0),\wt{\pi}_M\circ\wt{F}_t(y',0))
&=d(\Psi_{t+\wt{\chi}_Y(y)}y,\Psi_{t+\wt{\chi}_Y(y')}y')\\
&\ll d(\Psi_{S_n\wt{r}(y,0)+\wt\chi_Y(y)}y,\Psi_{S_n\wt{r}(y',0)+\wt\chi_Y(y')}y').
\end{split}
\]

Using the identities $\pi_X\circ f=T\circ\pi_X$ and $S_n\wt r=S_n(r\circ\pi_X)+\wt\chi\circ f^n-\wt{\chi}$, we get that
\[
S_n\wt{r}(y,0)+\wt\chi_Y(y)=
r_n(y)+\wt{\chi}(f^n(y,0)).
\] 
So,
\[
d\bigl(\wt{\pi}_M\circ\wt{F}_t(y,0),\wt{\pi}_M\circ\wt{F}_t(y',0)\bigr)\ll
d(\Psi_{r_n(y)+\wt{\chi}(f^n(y,0))}y,\Psi_{r_n(y')+\wt{\chi}(f^n(y',0))}y').
\]

Let $G_n=\wt \chi(f^n(y,0))-\wt \chi(f^n(y',0))$.
By $|\wt \chi|_\infty<\infty$, Equation~\eqref{eq:continuous_dependence} gives that 
\[\begin{split}
d(\Psi_{r_n(y)+\wt{\chi}(f^n(y,0))}y,\Psi_{r_n(y')+\wt{\chi}(f^n(y',0))}&y')\ll
d(\Psi_{r_n(y)+G_n}y,\Psi_{r_n(y')} y')\\
&\le d(\Psi_{r_n(y)+G_n}y,\Psi_{r_n(y)}y) + d(\Psi_{r_n(y)}y,\Psi_{r_n(y')} y')\\
&=d(\Psi_{r_n(y)+G_n}y,\Psi_{r_n(y)}y) + d(T^ny,T^ny').
\end{split}
\]
Therefore, using~\eqref{eq:Lipschitz} it follows that
\[
d\bigl(\wt{\pi}_M\circ\wt{F}_t(y,0),\wt{\pi}_M\circ\wt{F}_t(y',0)\bigr)
\ll |G_n|+d(T^ny,T^ny').
\]
Applying Lemma~\ref{lem:G_n} and~\eqref{eq:exp_contr}, we get that
$|G_n|+d(T^ny,T^ny')\ll \gamma_1^n$. 
Since $n+1\ge t/|\wt{r}|_\infty$ uniformly in $y\in Y$, the proof is finished letting $\gamma_2=\gamma_1^{1/|\wt{r}|_\infty}$.
\end{proof}
\end{proposition}

\begin{proposition}\label{proposition:regularity_final}
 Let $v\in C^\eta(M,\R^N)$. There exists $C>0$ such that
 \begin{equation*}
 |v\circ \wt{\pi}_M(y,u)-v\circ \wt{\pi}_M(y',s)|\le C\bigl((\inf_{Y_j}\wt \phi)(\theta^{\eta})^{s(y,y')}+|u-s|^\eta \bigr),
 \end{equation*}
 for all $(y,u), (y',s)\in Y^{\wt{\phi}}$ such that $y,y'\in \wt Y_j$.
\begin{proof}
Using $\wt\chi_Y=0$ on $\wt Y$ and applying Proposition~\ref{proposition:regularity},
\[
\begin{split}
|v\circ \wt{\pi}_M(y,u)-v\circ \wt{\pi}_M(y',s)|&=|v(\Psi_{u+\wt\chi_Y(y)}y)-v(\Psi_{s+\wt\chi_Y(y')}y')|\\
&=|v(\Psi_uy)-v(\Psi_sy')|\\
&\SMALL\ll(\inf_{Y_j}\phi)(\gamma_1^\eta)^{s(y,y')}+|u-s|^\eta.
\end{split}
\]
The inequality $\inf_{Y_j} \phi \le 2 \inf_{Y_j} \wt{\phi}$ follows from the proof of Proposition~\ref{proposition:phi_final}, hence concluding our estimate.
\end{proof}
\end{proposition}

\begin{remark}\label{rem:reduction}
By what we have seen in the current subsection, we can assume without loss that $\Psi_t\colon M\to M$ is modelled by a suspension flow 
$F_t\colon Y^\phi\to Y^\phi$
with a roof function $\phi\in L^p(Y)$, $p\ge2$, which is  constant along stable leaves and satisfies the condition of Proposition~\ref{proposition:phi_final}.
By Proposition~\ref{proposition:exp_decay_flow}, we can suppose that the projections of $F_t(y,0)$ and $F_t(y',0)$ into $M$ contract exponentially for any $y,y'\in Y$ such that $y'\in\mathcal W^s(y)$.
Finally, we can assume that Hölder observables on $M$ lifted to $Y^\phi$ satisfy the regularity condition of Proposition~\ref{proposition:regularity_final}.
\end{remark}

\subsection{Reduction to an observable constant along stable leaves}\label{subsec:sinai_flows} 

In this subsection we present an adaptation of~\cite[Theorem 5]{MelTor02} and~\cite[Theorem 7.1]{AraMelVar15} to our family of nonuniformly hyperbolic flows.

Following Remark~\ref{rem:reduction}, we consider $\phi\colon Y \to [1,\infty)$ satisfying $\phi(\pi y)=\phi(y)$ for all $y\in Y$, where $\pi\colon Y\to\wt{Y}$ is a projection to a chosen set $\wt{Y}\subset Y$ that intersects every stable leaf only once. Hence, if $(y,u)\in Y^\phi$ then $(\pi y,u)$ lies in $Y^\phi$ as well, and we can define $\pi\colon Y^\phi\to Y^\phi$ as $\pi(y,u)=(\pi y,u)$. 
Moreover, $\phi\in L^p(Y)$ for some $p\in[2,\infty]$ and there exist $C>0$ and $\gamma\in(0,1)$ such that 
\begin{equation}\label{eq:phi_final}\SMALL
|\phi(y)-\phi(y')|\le C(\inf_{Y_j}\phi)\gamma^{s(y,y')},
\end{equation}
for all $j\ge1$ and any $y,y'\in Y_j$. We let $F\colon Y\to Y$ be a map as defined in Subsection~\ref{subsec:setup}, and let  $F_t\colon Y^\phi\to Y^\phi$, $t\in\R$, be the suspension flow of $F$ on $Y^\phi$. 

As stated in Remark~\ref{rem:reduction}, there is a map $\pi_M\colon Y^\phi\to M$ such that 
\begin{enumerate}[(I)]
\item\label{item_regularity} For every $v\in C^\eta(M,\R^N)$ there exists $C>0$ such that
\begin{equation*}\SMALL
\bigl|v\circ\pi_M(y,u)-v\circ\pi_M(y',s)\bigr|\le C\bigl((\inf_{Y_j}\phi)\gamma^{s(y,y')}+|u-s|^\eta\bigr),
\end{equation*}
for all $(y,u), (y',s)\in Y^\phi$ such that $y,y'\in \wt{Y}_j$;
\item\label{item:exp_contr_flow} 
There exists $C>0$ such that 
\[
d\bigl(\pi_M\circ F_t(y,0),\pi_M\circ F_t(y',0)\bigr)
\le C\gamma^t,
\]
for all $t\ge0$ and $y,y'\in Y$ such that $y'\in\mathcal{W}^s(y)$.
\end{enumerate}

Let $\overline{Y}$ denote the quotient of $Y$ by the partition into stable leaves.
Since $\phi$ is constant along stable leaves, we can project it to
$\overline{\phi}\colon\overline{Y}\to[1,\infty)$ that is integrable with respect to $\overline{\mu}=\overline{\pi}_*\mu$. Let $\overline{Y}^{\overline{\phi}}$ be the suspension over the map $\overline{F}\colon\overline Y\to \overline Y$ with roof function $\overline{\phi}$, suspension semiflow 
$\overline{F}_t\colon\overline{Y}^{\overline{\phi}}\to \overline{Y}^{\overline{\phi}}$, $t\ge0$, and ergodic measure 
$\overline{\mu}^{\overline{\phi}}=(\overline{\mu}\times\text{Leb})/\int_Y \overline{\phi}\diff\overline{\mu}$. 
We have by~\eqref{eq:phi_final} that for all $j\ge1$ and  $y,y'\in\overline{Y}_j$,
\[\SMALL
|\overline{\phi}(y)-\overline{\phi}(y')|\le C(\inf_{\overline{Y}_j}\overline{\phi})\gamma^{s(y,y')}.
\]
\begin{remark}\label{rem:Gibbs-Markov_semiflow}
Since $\gamma\in(0,1)$, we know that $d_\gamma(y,y')=\gamma^{s(y,y')}$ is a metric on $\overline{Y}$ such that 
$d_\gamma(y,y')\le \gamma^{-1}d_\gamma(\overline{F}y,\overline{F}y')$. Hence, the function $\overline{\phi}$ satisfies equation~\eqref{eq:phi} with metric $d_\gamma$ and $\eta=1$. As in Section~\ref{sec:new_dec} and in~\cite[Definition~2.2]{BalButMel19} we say that $\overline{F}_t\colon\overline{Y}^{\overline{\phi}}\to \overline{Y}^{\overline{\phi}}$ is a Gibbs-Markov semiflow.
\end{remark}

\begin{proposition}\label{proposition:sinai_flows}
Let $v\in  C^\eta(M,\R^N)$. There exist bounded functions
$h,\chi\colon Y^{\phi}\to\R^N$ such that
\[
v\circ\pi_M=h+\chi-\chi\circ F_1.
\]
The function $h$ is constant along stable leaves (that is $h=h\circ\pi$). Hence, $h$ projects to an observable $\overline{h}\colon\overline{Y}^{\overline{\phi}}\to\R^N$. 

\begin{proof}

For $w=v\circ\pi_M$, write 
$\chi=\sum_{n=0}^\infty\{w\circ F_n-w\circ F_n\circ\pi\}$ and  define~$h=w+\chi\circ F_1-\chi$.
A calculation similar to~\eqref{eq:tilde_r_consant} (replacing $f$ with $F_1$) gives
\begin{equation}\label{eq:H}
\begin{split}
\chi\circ F_1-\chi
&=\sum_{n=0}^\infty
\bigl\{
w\circ F_{n+1}\circ\pi-w\circ F_n\circ\pi\circ F_1
\bigr\}-w+w\circ\pi
=H-w+w\circ\pi.
\end{split}
\end{equation}
Since $\pi\circ F_1=\pi\circ F_1\circ\pi$, we have that $h=w\circ\pi+H$ is  constant along stable leaves.

Using point~\ref{item:exp_contr_flow},
\[\begin{split}
\bigl|w( F_n(y,u))-w( F_n(\pi(y,u)))\bigr|&\le 
|v|_\eta d\bigl(\pi_M\circ F_{n+u}(y,0),\pi_M\circ F_{n+u}(\pi y,0)\bigr)^\eta\\
&\ll (\gamma^\eta)^{n+u}\le (\gamma^\eta)^n,    
\end{split}
\]
for all $n\ge0$. 
So, the series for~$\chi$ converges absolutely on  $Y^\phi$ and $\chi$ is bounded.
\end{proof}
\end{proposition}

\begin{definition}[Function space on $\overline{Y}^{\overline{\phi}}$]
Let $\theta\in(0,1)$ and $N\ge1$. For $j\ge1$, define  $\overline{Y}^{\overline{\phi}}_j=\{(y,u)\in Y^{\overline{\phi}}:y\in \overline{Y}_j\}$.
We denote with $\mathcal{H}^\theta(\overline{Y}^{\overline{\phi}},\R^N)$  the space of bounded observables  $g:\overline{Y}^{\overline{\phi}}\to\R^N$ such that
\[ 
|g|_{\theta}=\sup_{j\ge1}\sup_{(y,u),(y',u)\in \overline{Y}^{\overline{\phi}}_j,\ y\neq y'} \frac{|g(y,u)-g(y',u)|}{(\inf_{\overline{Y}_j}\overline{\phi})\theta^{s(y,y')}}<\infty.
\] 
\end{definition}

\begin{remark}\label{rem:correct_function_space}
As in Remark~\ref{rem:Gibbs-Markov_semiflow}, we write $d_\theta(y,y')=\theta^{s(y,y')}$ for $\theta\in(0,1)$. Any function $g\in\mathcal{H}^\theta(\overline{Y}^{\overline{\phi}},\R^N)$ satisfies
\[\SMALL
|g(y,u)-g(y',u)|\le \theta^{-1}|w|_\theta (\inf_{\overline{Y}_j}\overline{\phi})d_\theta(\overline{F}y,\overline{F}y'),
\]
for all $j\ge1$ and $(y,u),(y',u)\in\overline{Y}^{\overline{\phi}}_j$. Such a condition is the one
of Definition~\ref{def:functions_on_Y^phi}. Therefore, we can apply our results from Section~\ref{sec:new_dec} to any observable in $\mathcal{H}^\theta(\overline{Y}^{\overline\phi},\R^N)$.
\end{remark}

Using the notation
$\phi_m=\sum_{j=0}^{m-1}\phi\circ F^j$, we define for 
$y\in Y$ and $t\ge0$ the \textit{lap number}  
$M_t(y)=m$ to be the unique integer such that 
$\phi_m(y)\le t<\phi_{m+1}(y)$. 
Since $\phi(y)=\phi(\pi y)$, we have also that $M_t(y)=M_t(\pi y)$.
The suspension flow $F_t\colon Y^\phi\to Y^\phi$ can be written as 
\begin{equation}\label{eq:suspension_flow}
F_t(y,u)=\bigl(F^{M_{t+u}(y)}y,t+u-\phi_{M_{t+u}(y)}(y)\bigr).
\end{equation}

The next result shows that, if $v\in C^\eta(M,\R^N)$, then the function $\overline{h}$ defined in Proposition~\ref{proposition:sinai_flows} lies in 
$\mathcal{H}^{\theta}(\overline{Y}^{\overline{\phi}},\R^N)$ 
for some parameter $\theta\in(0,1)$. By Remarks~\ref{rem:Gibbs-Markov_semiflow} and~\ref{rem:correct_function_space}, such a regularity is sufficient to apply our results on the rates of convergence for semiflows.

\begin{proposition}\label{proposition:final_prop}
Let $v\in C^\eta(M,\R^N)$ and let  $\overline{h}\colon\overline{Y}^{\overline{\phi}}\to\R^N$ be from Proposition~\ref{proposition:sinai_flows}.
There exists $\theta\in(0,1)$ such that 
$\overline{h}\in\mathcal{H}^{\theta}(\overline{Y}^{\overline{\phi}},\R^N)$.
\begin{proof}
We follow the proof of~\cite[Theorem 7.1]{AraMelVar15} and use the same notation of Proposition~\ref{proposition:sinai_flows}.  
Let $C>0$ and $\gamma\in(0,1)$ be the constants in equation~\eqref{eq:phi_final} and let $\theta=\gamma^{\eta/2}$. Since $h$ is constant along stable leaves, 
it is sufficient to prove that there exists $K>0$ such that $|h(y,u)-h(y',u)|\le K (\inf_{Y_j}\phi)\theta^{s(y,y')}$ for all $(y,u),(y',u)\in Y^\phi$ such that $y,y'\in\wt{Y}_j$.

For all $j\ge1$, define $L_j=\gamma(1-\gamma)/(C\inf_{Y_j}\phi)$ and fix $y,y'\in \wt{Y}_j$ such that $\theta^{s(y,y')}\le L_j/2$ (the case $\theta^{s(y,y')}> L_j/2$ is trivial). Write $w=v\circ\pi_M$ and let  $u\in[0,\min\{\phi(y),\phi(y')\}]$. We deal separately with the two addends of  $h=w\circ\pi + H$, where $H$ is from equation~\eqref{eq:H}.   For the first one, we see by point~\ref{item_regularity} 
 that 
 \[\SMALL
 |w(\pi y,u)-w(\pi y',u)|\ll (\inf_{Y_j}\phi)\gamma^{s(\pi y, \pi y')}\le (\inf_{Y_j}\phi)\theta^{s(y, y')}.
 \]

Let us deal with $H$.
For $\widehat{N}=\lfloor s(y,y')/2\rfloor$, we write 
\[\SMALL
|H(y,u)-H(y',u)|\le A_1(y,y')+A_2(y,y')+B(y)+B(y'),
\]
where 
\begin{align*}
&\SMALL A_1(y,y')=\sum_{n=0}^{\widehat{N}-1}|w( F_{n+1}(\pi y,u))-w( F_{n+1}(\pi y',u))|,\nonumber\\
&\SMALL A_2(y,y')=\sum_{n=0}^{\widehat{N}-1}|w( F_n(\pi( F_1( y,u))))-w( F_n(\pi( F_1( y',u))))|,\nonumber\\
&\SMALL B(y)=\sum_{n=\widehat{N}}^{\infty}|w( F_{n+1}(\pi(y,u)))-w( F_n(\pi ( F_1(y,u))))|.\label{eq:B}
\end{align*}

Let us focus on $B$. By $\inf\phi\ge1$, we have
$M_{1+u}(y)=M_{1+u}(\pi y)=m\in\{0,1\}$. Using~\eqref{eq:suspension_flow},
\begin{equation*}\label{eq:no_lap}
 F_1(\pi(y,u))=(F^m\pi y, u')\qquad\text{and}\qquad
\pi ( F_1(y,u))=(\pi F^m y, u'),
\end{equation*}
where $u'$ is either $u$ or $u+1-\phi(y)$ (depending on $m$).
Since $F^m\pi y$ and $\pi F^my$ belong to the same stable leaf,
point~\ref{item:exp_contr_flow} yields,
\begin{align*}\label{eq:decay_B}\SMALL
B(y)&\SMALL=\sum_{n=\widehat{N}}^{\infty}
\bigl|v\circ\pi_M( F_n(\pi F^m y, u'))-v\circ\pi_M( F_n(F^m\pi y, u'))\bigr|\\
&\SMALL\le|v|_\eta \sum_{n=\widehat{N}}^\infty 
d\bigl(\pi_M\circ F_n(\pi F^m y, u')-\pi_M\circ F_n(F^m\pi y, u')\bigr)^\eta\\
&\SMALL\ll \sum_{n=\widehat{N}}^\infty (\gamma^\eta)^{n+u'}\le \sum_{n=\widehat{N}}^\infty (\gamma^\eta)^{n}=
(\gamma^\eta)^{\widehat{N}}/(1-\gamma^\eta).
\end{align*}
By our choice of $\widehat{N}$, we get $B(y)\ll(\gamma^\eta)^{\widehat{N}}\ll(\gamma^{\eta/2})^{s(y,y')}=\theta^{s(y,y')}$.
\vspace{1ex}

Le us deal with $A_1$ and $A_2$. Note that for all 
$n=0,\dots, \widehat{N}+1$, 
equation~\eqref{eq:phi_final} yields
\begin{equation*}\begin{split}\label{eq:same_laps}
|\phi_n(y)-\phi_n(y')|&\le C(\inf_{Y_j}\phi)\sum_{i=0}^{n-1}
\gamma^{s( F^iy, F^iy')}\le C
(\inf_{Y_j}\phi)\sum_{i=0}^{\widehat{N}}\gamma^{s(y,y')-i}\\
&\le \frac{C(\inf_{Y_j}\phi)}{1-\gamma}\gamma^{s(y,y')-\widehat{N}}\le
\frac{C(\inf_{Y_j}\phi)}{(1-\gamma)\gamma}\gamma^{s(y,y')/2}\le
L_j^{-1}\theta^{s(y,y')}\le\frac12\le\frac{\inf\phi}{2}.
\end{split}
\end{equation*}
Hence, for 
$n=0,\dots,\widehat{N}$, the intervals $[\phi_n(y),\phi_{n+1}(y)]$ and $[\phi_n(y'),\phi_{n+1}(y')]$ have the initial points and endpoints closer than $(\inf \phi)/2$.
Let~$k=M_{n+u}(y)$ and $k'=M_{n+u}(y')$; by $\inf \phi\ge1$ it follows that $k,k'\le n\le \widehat{N}$. Moreover, by
\[
n+u \in [\phi_k(y),\phi_{k+1}(y)]\cap [\phi_{k'}(y'),\phi_{k'+1}(y')]
\]
we have that $|k-k'|\le1$.

For $n=1,\dots,\widehat{N}$, let $a_1(y,y')=|w( F_n(\pi y,u))-w( F_n(\pi y',u))|$. By~\eqref{eq:suspension_flow},
\[\begin{split}\SMALL
a_1(y,y')=
|w( F^k\pi y,n+u-\phi_k(y))-w( F^{k'}\pi y',n+u-\phi_{k'}( y'))|.
\end{split}
\]
For $\ell=0,\dots,\widehat{N}-1$, let $a_2(y,y')=|w( F_\ell(\pi( F_1( y,u))))-w( F_\ell(\pi( F_1( y',u))))|$ and define
\[\begin{split}
    &k_1=M_{1+u}(y)\qquad k_2=M_{\ell+u+1-\phi_{k_1}(y)}( F^{k_1}y)\\
    &k'_1=M_{1+u}(y')\qquad k'_2=M_{\ell+u+1-\phi_{k'_1}(y')}( F^{k'_1}y').
\end{split}
\]
Writing $n=\ell+1\in\{1,\dots,\widehat{N}\}$, we have by~\eqref{eq:suspension_flow} that
\[
a_2(y,y')=|w( F^{k_2}\pi F^{k_1}y,n+u-\phi_k(y))-w( F^{k'_2}\pi F^{k'_1}y',n+u-\phi_{k'}(y'))|.
 \]
Note that $k_1+k_2=k=M_{n+u}(y)$ and $k'_1+k'_2=k'=M_{n+u}(y')$. 
We assume without loss that $k\ge k'$.

We claim for $i=1,2$ that
\begin{equation}\label{eq:global_aim}\SMALL
a_i(y,y')\ll (\inf_{Y_j}\phi)\gamma^{s(F^k y, F^k y')}+|\phi_k(y)-\phi_k(y')|^\eta.
\end{equation}
Assuming the claim and by equation~\eqref{eq:phi_final}, we can bound
\begin{equation*}\label{eq:k=h}
\begin{split}\SMALL
a_i(y,y') &\SMALL\ll
(\inf_{Y_j}\phi)\bigl\{\gamma^{s(y,y')-k}+
\sum_{\ell=0}^{k-1}(\gamma^\eta)^{s( F^\ell y, F^\ell y')}\bigr\}\\
&\SMALL= (\inf_{Y_j}\phi)\bigl\{\gamma^{s(y,y')-k}+
\sum_{\ell=0}^{k-1}(\gamma^\eta)^{s(y,y')-\ell}\bigr\}\\
&\SMALL=
(\inf_{Y_j}\phi)\bigl(1+\gamma^\eta/(1-\gamma^\eta)\bigr)(\gamma^\eta)^{s(y,y')-k}.
\end{split}
\end{equation*} 
By $k\le n$,we have $a_i(y,y')\ll(\inf_{Y_j}\phi)(\gamma^\eta)^{s(y,y')-n}$ and hence
\begin{equation*}
(\inf_{Y_j}\phi)^{-1}A_i(y,y')\ll \sum_{n=1}^{\widehat{N}}(\gamma^\eta)^{s(y,y')-n}\le
\frac{(\gamma^\eta)^{s(y,y')-\widehat{N}}}{1-\gamma^\eta}\ll
(\gamma^\eta)^{s(y,y')/2}=
\theta^{s(y,y')}.
\end{equation*}
The main statement is proven by combining the estimates for $A_1,A_2$ and $B$.\vspace{1ex}

Let us show the claim for $a_1(y,y')$. 
By $|k-k'|\le 1$, we have to deal with two cases. If $k=k'$, equation~\eqref{eq:global_aim} follows from point~\ref{item_regularity}. 
If $k=k'+1$, using again point~\ref{item_regularity} and the identifications on $ Y^\phi$,
\begin{equation}\label{eq:k=h-1}
\begin{split}
a_1(y,y')&\le 
|w( F^{k'}\pi y',n+u-\phi_{k'}(y'))-w( F^{k'}\pi   y',\phi( F^{k'} y'))|\\
&\quad +|w( F^{k'+1}\pi y',0)-w( F^{k'+1}\pi y,n+u-\phi_{k'+1}(y))|\\
&\SMALL\le C(\inf_{Y_j}\phi)
\gamma^{s( F^{k'+1}\pi y, F^{k'+1}\pi y')}\\
&\quad+
(\phi( F^{k'} y')-n-u+\phi_{k'}(y'))^\eta+(n+u-\phi_{k'+1}(y))^\eta.
\end{split}
\end{equation}
Note that on the last line of~\eqref{eq:k=h-1} the quantities inside the brackets are positive. Hence, equation~\eqref{eq:global_aim} follows by the general inequality $\alpha^\eta+\beta^\eta\le2(\alpha+\beta)^\eta$ for all $\alpha,\beta\ge 0$.

To deal with  $a_2(y,y')$, we use the identity $s(F^{k_2}\pi F^{k_1}y,F^{k'_2}\pi F^{k'_1}y')=s(F^ky,F^{k'}y')$, and the fact that $\phi$ is constant on stable leaves. 
Hence, the claim is proven analogously to $a_1(y,y')$: if $k=k'$, equation~\eqref{eq:global_aim} follows by point~\ref{item_regularity}; if $k=k'+1$, we conclude by   reasoning as in~\eqref{eq:k=h-1}.
\end{proof}
\end{proposition}

\begin{proof}[\textbf{Proof of Theorem \ref{theorem:rates_NUH_flows}}]
For $v\in C^\eta(M,\R^N)$, let $w=v\circ\pi_M$.
By Proposition~\ref{proposition:sinai_flows}, there exist functions $h,\chi\colon Y^\phi\to\R^N$ such that $w=h+\chi-\chi\circ F_1$ with $\chi\in L^\infty(Y^\phi,\R^N)$. Let 
$\overline{h}\colon\overline{Y}^{\overline{\phi}}\to\R^N$ be the projection of $h$ via $\overline{\pi}$.
For $t\in[0,1]$, we define the sequence 
\[\SMALL
\overline{W}_n(t)=n^{-1/2}\int_0^{nt}\overline{h}\circ\overline{F}_s\diff s,
\]
on the probability space $(\overline{Y}^{\overline{\phi}},\overline{\mu}^{\overline{\phi}})$. The first step of the proof is to show that
\[\mathcal W(W_n,\overline{W}_n)\ll n^{-1/2}
\qquad\text{and}\qquad
\Pi(W_n,\overline{W}_n)\ll n^{-1/2}.
\]

Define on the space $(Y^\phi,\mu^\phi)$ the sequences
\[
\SMALL W_n'(t)=n^{-1/2}\int_0^{nt}w\circ F_s\diff s\qquad\text{and}\qquad
\SMALL W_n''(t)=n^{-1/2}\int_0^{nt}h\circ F_s\diff s.
\]
Note that $W_n'(t)=W_n(t)\circ\pi_M$ and  $\overline{W}_n=W''_n\circ\overline{\pi}$. Since $\pi_M$ and $\overline{\pi}$ are measure-preserving, we have $W_n=_dW_n'$ and $\overline{W}_n=_dW''_n$. 
For $x>1$, we get by a change of variables
\[\SMALL
 \int_0^x(\chi-\chi\circ F_1)\circ F_s\diff s=
 \int_0^x\chi\circ F_s\diff s-
 \int_1^{x+1}\chi\circ F_{s}\diff s=
 \int_0^1\chi\circ F_s\diff s-\int_x^{x+1}\chi\circ F_s\diff s,
\]
which implies that 
$|\int_0^x(\chi-\chi\circ F_1)\circ F_s\diff s|\le 2|\chi|_\infty$, for all $x>0$.
So,
\[\SMALL
\bigl|\sup_{t\in[0,1]}|W'(t)-W''(t)|\bigr|_\infty=
n^{-1/2}\bigl|\sup_{t\in[0,1]}\bigl|\int_0^{nt}(\chi-\chi\circ F_1)\circ F_s\bigr|\bigr|_\infty\le 2 n^{-1/2}|\chi|_\infty.
\]

For the Wasserstein metric, we get
\[\SMALL
\mathcal{W}(W_n,\overline{W}_n)=
\mathcal{W}(W'_n,W''_n)
\le \mathbb{E}[\sup_{t\in[0,1]}|W'_n(t)-
W''_n(t)|]\ll n^{-1/2}.
\]
Since $\Pi$ is bounded by the infinity norm, 
$\Pi(W_n,\overline{W}_n)=\Pi(W'_n,W''_n)\ll n^{-1/2}$.
Hence, we are left to estimate the rates for $\overline{W}_n$.

Using Remarks~\ref{rem:Gibbs-Markov_semiflow} and~\ref{rem:correct_function_space}, we see that $\overline{F}_t$ is a Gibbs-Markov semiflow and the observable $\overline{h}$ belongs to the functional space of Definition~\ref{def:functions_on_Y^phi}.
So, the sequence $\overline{W}_n$ coincides with $\widehat{W}_n$ of~\eqref{eq:Wn'}.   Hence, the rates for
$\mathcal{W}(\overline{W}_n,W)$ and $\Pi(\overline{W}_n,W)$ can be deduced as in the proofs of Theorem~\ref{theorem:rate_wass_flow} and~\ref{theorem:rates_N=1_flow} (which can be found in Section~\ref{sec:continuous_time_rates}).
\end{proof}

\subsection*{Acknowledgements}
The author expresses all his gratitude to his PhD supervisor, Prof. Ian Melbourne, for the thorough guidance in the completion of his thesis~\cite{Pav23}, which made this paper possible. The author also thanks Alexey Korepanov for his valuable help and discussions on the topic. A special thanks goes to Nicholas Fleming Vázquez for sharing reference~\cite{CunDedMer20}, which is essential to prove the rates in 1-Wasserstein independently of the dimension.
I am grateful to the anonymous referees for their detailed feedback and corrections, which helped improve both the mathematical clarity and the overall readability of this manuscript.

\subsection*{Statements and declarations}
\noindent \textbf{Funding.}
This research was supported by ISF research grant 3056/21.

\noindent \textbf{Conflict of interest.} The author declares that he has no financial interests.

\noindent \textbf{Data availability.} There is no associated data.

\end{document}